\documentclass[a4paper]{article}
\usepackage{endnotes}
\let\footnote=\endnote

\usepackage{amsmath,amsthm}
\usepackage{amssymb}
\usepackage{subcaption}
\usepackage{a4wide}

\usepackage{graphicx}
\usepackage{tikz}
\usetikzlibrary{plotmarks}
\usepackage{pgfplots}
\tikzset{mark size=2.5}
\usepackage{enumitem}
\usepackage[pdftex,colorlinks=true,urlcolor=blue,citecolor=blue,pdfstartview=FitH]{hyperref}
\usepackage{subcaption}

\usepackage{multirow}
\usepackage{bigstrut}

\definecolor{col4}{rgb}{0.8, 0.01568627450980392, 0.}
\definecolor{col5}{rgb}{0.9372549019607843, 0.6274509803921569, 0.1}
\definecolor{col2}{rgb}{0.9686274509803922, 0.87, 0.}
\definecolor{col3}{rgb}{0.3176470588235294, 0.58, 0.0784313725490196}
\definecolor{col1}{rgb}{0.22745098039215686, 0.23921568627450981, 
  0.43}
\definecolor{col6}{rgb}{0.5019607843137255, 0.16, 0.4470588235294118}

\newtheorem{lemma}{Lemma}
\newtheorem{theorem}{Theorem}
\newtheorem{proposition}{Proposition}

\theoremstyle{definition}	
\newtheorem{remark}[theorem]{Remark}

\newtheorem{algorithm}{Algorithm}
\newtheorem{conjecture}{Conjecture}
\newtheorem{insight}{Managerial insight}

\def \l {\lambda}
\def \m {\mu}
\def \d {\delta}
\def \b {\beta}
\def \g {\gamma}
\def \a {\alpha}
\def \iy {\infty}
\def \E {\mathbb{E}}
\def \P {\mathbb{P}}

\def \QY {\hat Q_1^b(t)}
\def \QH {\hat Q_1^h(t)}
\def \CY {\hat Q_2^b(t)}
\def \CH {\hat Q_2^h(t)}

\def \qy {q_1^b}
\def \qh {q_1^h}
\def \cy {q_2^b}
\def \ch {q_2^h}

\def \Poi {{\rm Pois}}

\def \less[#1] {\stackrel{(#1)}{\leq}}
\def \ra[#1] {\stackrel{(#1)}{\Rightarrow}}
\def \ed {\stackrel{d}{=}}

\begin{document}

\title{The restricted Erlang-R Queue: Finite-size effects in service systems with returning customers}

\author{Johan S.H. van Leeuwaarden\footnote{
Department of Mathematics and Computer Science, Eindhoven
University of Technology, P.O. Box 513, 5600 MB Eindhoven, The
Netherlands (\{j.s.h.v.leeuwaarden,b.w.j.mathijsen,f.sloothaak\}@tue.nl)} \and 
Britt W.J. Mathijsen\footnotemark[1] \and 
Fiona Sloothaak\footnotemark[1] \and 
Galit B. Yom-Tov \footnote{Faculty of Industrial Engineering and Management, Technion - Israel Institute of Technology, 32000 Haifa, Israel, (gality@ie.technion.ac.il)}
}

\maketitle

\begin{abstract}Motivated by health care systems with repeated services that have both personnel (nurse and physician) and space (beds) constraints, we study a restricted version of the Erlang-R model. The space restriction policies we consider are blocking or holding in a pre-entrant queue. We develop many-server approximations for the system performance measures when either policy applies, and explore the connection between them.
We show that capacity allocation of both resources should be determined simultaneously, and derive the methodology to determine it explicitly. We show that the system dynamics is captured by the fraction of needy time in the network, and that returning customers should be accounted for both in steady-state and time-varying conditions. We demonstrate the application of our policies in two case studies of resource allocation in hospitals.
 \end{abstract}

\section{Introduction}

Because service systems are stochastic in nature, it is common practice to use queueing theory for performance analyses and workforce planning. 
Traditionally, systems are modeled after a single station queue, such as the $M/M/s$ (Erlang-C), $M/M/s/s$ (Erlang-B) or $M/M/s+M$ (Erlang-A) models, and fluid and diffusion approximations are used to provide insights into the process dynamics. 
However, single station models often fail to capture the more intricate dynamics of service networks. 
In a health care context, think of flows of patients in a hospital from one medical ward to another \cite{Armony2015}, within the Emergency Department (ED) between different stages of treatment \cite{Junfei2015}, or between medical facilities \cite{zychlinski2016bedblocking}.  
Queueing networks can capture the dependency between several service stages and several resources.  More specifically, we are interested in the ubiquitous feature, particularly present in health care environments, that customers during their stay in the system might require a specific resource multiple times, e.g.\ physicians and nurses who treat patients several times during their stay in the medical wards \cite{Jennings2011} or the ED \cite{YomTov2014}, while multiple resources are limited (e.g.\ physicians, nurses and beds).

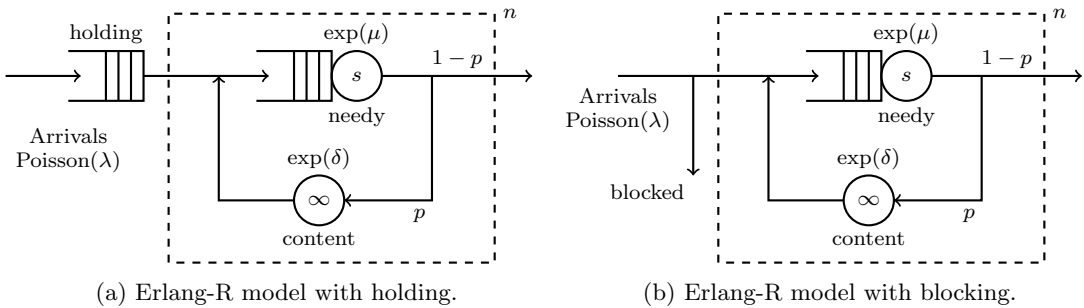
\begin{figure} 
\centering
\begin{subfigure}{0.48\textwidth}
\begin{tikzpicture}[scale = 0.66]
\footnotesize
\draw [thick, ->] (-2.75,4.5) -- (-1.25,4.5);
\draw [thick] (-1.5,5) -- (0,5) -- (0,4) -- (-1.5,4);
\draw [thick] (-0.25,4) -- (-0.25,5);
\draw [thick] (-0.5,4) -- (-0.5,5);
\draw [thick] (-0.75,4) -- (-0.75,5);
\draw [thick, ->] (0,4.5) -- (2.5,4.5);
\draw [thick] (2.25,5) -- (3.75,5) -- (3.75,4) -- (2.25,4);
\draw [thick] (3,4) -- (3,5);
\draw [thick] (3.25,4) -- (3.25,5);
\draw [thick] (3.5,4) -- (3.5,5); 
\draw [thick] (4.25,4.5) circle [radius=0.5] node {$s$} node[below= 0.3 cm] {needy} node[above=0.3cm] {$\exp(\mu)$};
\draw [thick, ->] (4.75,4.5) -- (7.75,4.5);
\draw [thick,->] (5.75,4.5) -- (5.75,2) -- (4,2);
\draw [thick] (3.5,2) circle [radius=0.5] node {$\infty$} node[below=0.3cm] {content} node[above=0.3cm] {$\exp(\delta)$};
\draw [thick,->] (3,2) -- (1.5,2) -- (1.5,4.5);
\node [above] at (-1.5,3) {\footnotesize Arrivals};
\node [above] at (-1.5,2.4) {\footnotesize Poisson($\lambda$)};
\node [below] at (5.5,2) {\footnotesize $p$};
\node [above] at (6.25,4.5) {\footnotesize $1-p$};
\draw [thick, dashed] (0.5,0.75) rectangle (7,5.75) node[right] {\footnotesize $n$}; 
\node [above] at (-0.75,5) {\footnotesize holding};
\end{tikzpicture}
\caption{Erlang-R model with holding.}
\label{fig:Erlang_R_holding}
\end{subfigure}
\begin{subfigure}{0.48\textwidth}
\begin{tikzpicture}[scale = 0.66]
\footnotesize
\draw [thick, ->] (-1.5,4.5) -- (2.5,4.5);
\draw [thick, ->] (0,4.5) -- (0,2.5) node[below left] {blocked};
\draw [thick] (2.25,5) -- (3.75,5) -- (3.75,4) -- (2.25,4);
\draw [thick] (3,4) -- (3,5);
\draw [thick] (3.25,4) -- (3.25,5);
\draw [thick] (3.5,4) -- (3.5,5); 
\draw [thick] (4.25,4.5) circle [radius=0.5] node {$s$} node[below= 0.3 cm] {needy} node[above=0.3cm] {$\exp(\mu)$};
\draw [thick, ->] (4.75,4.5) -- (7.75,4.5);
\draw [thick,->] (5.75,4.5) -- (5.75,2) -- (4,2);
\draw [thick] (3.5,2) circle [radius=0.5] node {$\infty$} node[below=0.3cm] {content} node[above=0.3cm] {$\exp(\delta)$};
\draw [thick,->] (3,2) -- (1.5,2) -- (1.5,4.5);
\node [above] at (-1.5,3.8) {\footnotesize Arrivals};
\node [above] at (-1.5,3.2) {\footnotesize Poisson($\lambda$)};
\node [below] at (5.5,2) {\footnotesize $p$};
\node [above] at (6.25,4.5) {\footnotesize $1-p$};
\draw [thick, dashed] (0.5,0.75) rectangle (7,5.75) node[right] {\footnotesize $n$};
\end{tikzpicture}
\caption{Erlang-R model with blocking.}
\label{fig:Erlang_R_blocking}
\end{subfigure}
\caption{Restricted Erlang-R models with maximally $n$ customers in system.}
\label{fig:Erlang_R_model}
\end{figure}

An often ignored yet essential feature of service networks concerns the restriction on the number of customers that can reside in the system simultaneously. 
Call centers, for instance, only have a finite number of trunk lines \cite{Khudyakov2010}, whereas in health care facilities, the number of beds poses a constraint on the number of patients that can be admitted at the same time to a medical unit. 
In such situations, two or more resources (personnel, trunk lines, beds) are restricted. 
In this paper, we investigate the influence of such multiple restrictions on the network dynamics and the required staffing policies.\\
\\*
\noindent
\textbf{The restricted Erlang-R model.}
The canonical model for service networks with returns is the Erlang-R model \cite{YomTov2014} in which customers, during their stay in the system receive a random number of services from the same pool of servers. Yom-Tov \& Mandelbaum~\cite{YomTov2014} showed that such a simple network model can be used to determine staffing in an Israeli ED both in stable and time-varying conditions. Nevertheless, empirical studies report that some countries, such as the US, use a different operational mode that apply strict restrictions on entering the ED \cite{EDexperiment}. In typical US EDs, a patient will not enter the ED until both a bed and a physician are available to treat her. Those restrictions can be either physical (beds) restrictions or managerial ones---for instance by imposing a patient-to-physician ratio. In this work, we extend the Erlang-R model by enforcing a constraint on the maximum number of available places inside the service facility. Our model hence incorporates two kinds of resource constraints: servers that provide the actual service and the maximum available places inside the service system. Both affect the system in a highly interdependent way. The model, presented in Figure \ref{fig:Erlang_R_model}, assumes $s$ servers and a maximum capacity of $n$ concurrent places. We assume that customers arrive according to a Poisson process with rate $\lambda(t)$. In case a new arrival finds $n$ or more customers already present, we consider two options---either she waits outside the service facility in a holding queue until a vacant space becomes available (Figure \ref{fig:Erlang_R_holding}) or she is blocked (Figure \ref{fig:Erlang_R_blocking}), such as is the case when patients are sent to an alternative facility. Once customers are admitted, they require assistance from one of the $s$ servers for an exponentially distributed duration with mean $1/\mu$. Then, with probability $1-p$, customers leave the system or, with probability $p$, return to service again after an exponentially distributed time with mean $1/\delta$. Following Jennings \& de V\'ericourt~\cite{Jennings2011} and Yom-Tov \& Mandelbaum~\cite{YomTov2014}, we call patients {\it needy} when they require attention from one of the servers and {\it content} when they are in the delayed return phase. In addition, we call customers {\it holding} when they are waiting outside the facility for an available space. We assume that the arrival process, the needy times and content times are mutually independent. 
In the holding queue and the needy queue, we apply the First-Come-First-Served (FCFS) discipline.



As mentioned, we consider two versions of the finite-capacity constraint. 
The first version is called \emph{Erlang-R with holding}, in which customers wait for an available space in the system. 
The second version is called \emph{Erlang-R with blocking}, in which customers meeting a full system are blocked.
Naturally, intermediate scenarios can be constructed in which a proportion of the total arrival volume of customers indeed leaves upon finding a full system, while the rest joins the holding room. 
While this paper focuses on the two extreme cases, straightforward adaptions can fit these intermediate scenarios. \\
\\*
\noindent
\textbf{Examples of restricted Erlang-R.}
As noted before, an ED operated in the US can be modeled using a restricted Erlang-R model. Another health care example are Medical Units (MU) in a hospital.
Such units specialize in specific types of illnesses (cadriatric, oncology, etc.) and have limited resources such as nurses and beds. If the unit is full, new patients are either allocated to an alternative medical unit, i.e.\ blocked, or wait for an available bed. Both policies are problematic in terms of quality-of-care, because the personnel in the alternative unit (or the ED) may be less knowledgable and waiting in the ED was shown to increase mortality. 
Moreover, ED waiting may reduce available capacity for treating ED patients \cite{Carmen2016,israelit}, hence endangering both the delayed patient as well as others. Both the number of personnel (nurses and physicians) and the number of beds impact service dynamics and quality-of-care. Research so far looked at the capacity allocation of those resources separately. Green \& Yankovic~\cite{GY2011} and Jennings \& de V\'ericourt \cite{Jennings2008} looked at nurse staffing in medical units, while Bekker et al.~\cite{Bekker2009b} looked into bed allocation. The unified model we suggest enables us to capture the dependency between those two decisions, and its impact on other medical units in the hospital.
At the same time, we capture the two most commonly used modes of operation---blocking and holding of new patients.  

\subsection{Contributions}

Our main goal is to provide staffing policies for the restricted Erlang-R models that ensures high resource utilization, while at the same time maintains a good quality-of-care. This goal relates to the philosophy of the Quality-and-Efficiency-Driven (QED) regime known in many-server asymptotic theory. 
We discuss the main ideas behind this regime further in \S 2.
In this paper, we obtain asymptotic results for the Erlang-R model with blocking in the QED regime (\S\ref{sec:QED_limit_block}). Following \cite{Jennings2008}, we employ a two-fold QED staffing policy: $s=R_1 +\beta \sqrt{R_1}$  for the number of servers and $n=R_1/r+\gamma \sqrt{R_1/r}$ for the number of customers in the system, where   $\beta$ and $\gamma$ are constants, $R_1$ is the offered load of the servers and $r$ is the fraction of time a customer spends in the needy state. We establish limiting expressions for performance measures, such as the probability of delay and blocking, in the form of explicit functions that depend solely on $\beta$ and $\gamma$. 
In deriving these limit results, we use the available product-form solution for the stationary distribution.

Likewise, we pursue QED performance for the Erlang-R model with holding.
However, a direct analytic approach is obstructed by the absence of product-form solutions. We provide two solutions for establishing QED behavior. 
First, we provide stochastic performance bounds that stay meaningful in the QED regime (\S \ref{sec:bounds}), which demonstrate the non-degenerate behavior of the two-fold scaling in the large-system limit. Second, we develop a heuristic method that quantifies the difference between the scalable holding model and the blocking model (\S \ref{sec:QED_limit_holding}). This is based on the following unique approach: initially blocked patients in the blocking model are seen as if they reattempt to get access after a some delay. The behavior of this retrial model then resembles the Erlang-R model with blocking without retrials, yet with an increased arrival rate. The increase in arrival rate turns out to be the solution of a fixed-point equation. Using our results on the asymptotic behavior of the model with blocking in the QED regime, we then obtain approximative QED performance measures for the model with holding. Finally, we use these QED results to develop algorithms for dimensioning and time-varying staffing (\S \ref{sec:dimensioning_block}). 


Using the approximations developed, numerical analysis and simulation we provide the following managerial insights:
\begin{itemize}
\item We show that all resource allocation of personnel and beds should be synchronized in order to avoid waste, and that the QED scaling provides an efficient, flexible, and easy to implement methodology to do so.
\item We conclude that \emph{reentrant} customers in a restricted network are more significant than in an open network (Erlang-R). In contrast to the open model, in which returning customers need to be accounted for only in time-varying systems, the restricted Erlang-R model requires explicit consideration of returning customers under stationary conditions as well. 
\item We show that the influence of the network structure on the system dynamics crucially depend on the fraction of time a patient spends being needy during her stay in the system.
We then explore the influence of $r$ on operational decisions.
\item Combining the theoretical results, we explore the implication of managerial decisions in designing an MU (\S \ref{sec:dimensioning}) and an ED (\S \ref{sec:case_study}). In \S \ref{sec:dimensioning}, we show that enabling customers to hold in ED before entering an MU, requires more resources both in the MU and the ED. In Section \ref{sec:case_study}, we compare the pros and cons of imposing strict constraints on entering an ED.  We find that size restrictions have the ability to improve the quality-of-service of the processes within the facility, at the expense of a slight increase in pre-entrant wait and server efficiency levels.
\item Finally, we show that restricting the number of admitted customers protects those customers with complicated demand consisting of relatively many retrials/interruptions (\S \ref{subsec:num_visit}).
\end{itemize}

When dealing with time-varying patient demands, one requires some modifications in the QED results in order to obtain stable performance at every moment in time. We transform the two-fold staffing policy into a time-varying one based on the Modified Offered Load (MOL) method (\S \ref{sec:case_study}). This method approximates the offered load at the needy station at each point in time via a corresponding system with ample resources. By noting that the latter system coincides with the Erlang-R model with ample servers, we adopt the offered load approximations given in~\cite{YomTov2014}. Numerical experiments justify this method and we use this approach in our case study \S \ref{sec:case_study}.  

\section{Literature review}
%
%
Due to increasing demand and tightening budgets in health care, there is a growing need for efficient workforce management \cite{Green2008}. Personnel (nurse and physician) expenditure is one of the biggest factors in hospital costs \cite{Kazahaya2005}, and inadequate nursing levels have been mentioned as a significant factor in medical errors and ED overcrowding. In order to establish appropriate nursing levels, a staffing policy requires assessment of a wide range of variables, such as differing nurse expertise and patient acuity during the day. Current methods, such as the minimum nurse-to-patient ratios, are often too inflexible to capture those varying conditions. The American Hospital Association (AHA) and others call for dynamic staffing policies that can deal with the complex and evolving nature of health care \cite{AHA2007}.
Workforce management in health care systems has been studied extensively; see \cite{Denton2013,Hall2006,Hall2012} for overviews.
In recent years it has become apparent that queueing models can be helpful in developing staffing and routing recommendations, not just for large-scale service systems, but also for the small and complicated health care systems.

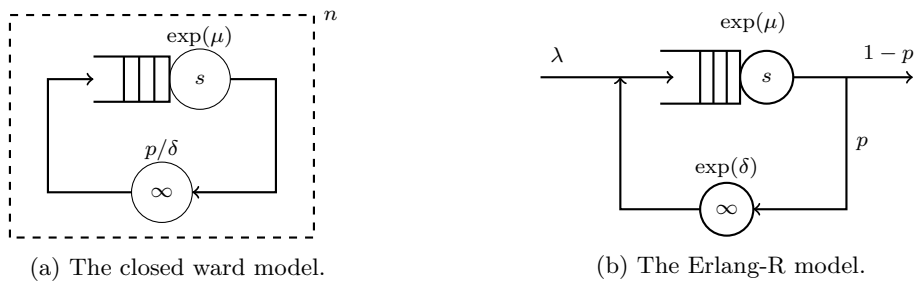
\begin{figure} \centering

\begin{subfigure}{0.48\textwidth}
\centering
\begin{tikzpicture}[scale=1]
\draw [dashed, thick] (-0.5,-0.1) rectangle (3.5,2.85) node[right] {\footnotesize $n$};
\draw [thick,->] (1.1,0.5) -- (0,0.5) -- (0,2) -- (0.6,2);
\draw [thick,->] (2.4,2) -- (3,2) -- (3,0.5) -- (1.9,0.5);
\draw [thick] (0.6,1.7) -- (1.6,1.7) -- (1.6,2.3) -- (0.6,2.3);
\draw [thick] (1,1.7) -- (1,2.3);
\draw [thick] (1.2,1.7) -- (1.2,2.3);
\draw [thick] (1.4,1.7) -- (1.4,2.3);
\draw (2.4,2) -- (3,2) -- (3,0.5) -- (1.9,0.5);
\draw (1.5,0.5) circle [radius=0.4] node[above=0.3cm] {\footnotesize $p/\d$} ;
\draw (2,2) circle [radius=0.4] node[above=0.3cm] {\footnotesize exp($\mu$)};
\node at (1.5,0.5) {\footnotesize $\infty$};
\node at (2,2) {\footnotesize $s$};
\end{tikzpicture}
\caption{The closed ward model.}
\label{fig:Jennings}
\end{subfigure}
\begin{subfigure}{0.48\textwidth}
\centering
\begin{tikzpicture}[scale=0.7]
\draw [thick, ->] (0,4.5) node[above=0.3cm,right] {\footnotesize $\l$} -- (2.5,4.5);
\draw [thick] (2.25,5) -- (3.75,5) -- (3.75,4) -- (2.25,4);
\draw [thick] (3,4) -- (3,5);
\draw [thick] (3.25,4) -- (3.25,5);
\draw [thick] (3.5,4) -- (3.5,5);
\draw [thick] (4.25,4.5) circle [radius=0.5];
\draw [thick, ->] (4.75,4.5) -- node[above=0.3cm,right] {\footnotesize $1-p$} (7,4.5);
\draw [thick,->] (5.75,4.5) -- node[right] {\footnotesize $p$} (5.75,2) -- (4,2);
\draw [thick] (3.5,2) circle [radius=0.5];
\draw [thick,->] (3,2) -- (1.5,2) -- (1.5,4.5);
\node at (4.25,4.5) {\footnotesize $s$};
\node at (3.5,2) {\footnotesize $\infty$};
\node [above] at (4,5.15) {\footnotesize exp($\mu$)};
\node [above] at (3.5,2.4) {\footnotesize exp($\delta$)};
\end{tikzpicture}
\caption{The Erlang-R model.}
\label{fig:ErlangR}
\end{subfigure}
\caption{Two basic models with interrupted services.}
\end{figure}

The first to try such an approach through queueing models were Green \& Savin~\cite{Green2006,Green2008} who used the single station stationary Erlang-C model to set staffing levels in EDs and panel sizes for clinics. Using a similar approach, \cite{Bekker2009a} used Erlang-B model to determine bed allocation for medical wards. 
The first to observe the significant impact of interrupted services in a health care setting were Jennings \& de V\'ericourt~\cite{Jennings2008,Jennings2011}. Motivated by the need to set nurse-to-patient ratios for internal wards, they considered a closed queueing system with $s$ nurses and $n$ beds, whuch we will refer to as the \emph{closed ward} model. This is essentially the Erlang-C model with the additional restriction that a finite population of the $n$ patients requires care. In their model, all beds are always occupied, and patients alternate between two phases: the needy phase where patients require service of a nurse and the content phase where they do not; see Figure \ref{fig:Jennings}. The system dynamics of restricted Erlang-R model are equivalent to those of the closed ward model of \cite{Jennings2008} if the holding queue would never be empty. 

Campello et al.~\cite{Campello2016} analyzed a similar operational decision, referred to as ED case management, which determines the maximal number of patients a physician should handle in parallel. They also used queueing networks and analyzed the stationary distribution. Note that in practice such decision is not only affected by operational measurements such as waiting times, but also by psychological constraints that limit physician capability to manage multiple tasks (patients) in parallel. Diwas~\cite{diwas} provided empirical evidence that physicians should not treat more than 6-7 patients at the same time. Therefore, many hospitals in the US restrict entrance to EDs even if beds are available if physicians are overloaded. 
We too consider such constraints, and analyze their impact on performance. We take a different approach than \cite{Campello2016}; instead of analyzing numerically steady-state distributions, we develop many-server approximations that can produce insight into the system dynamics, and can be incorporated into time-varying staffing procedures (see \S \ref{sec:case_study}).  


The model in~\cite{Jennings2008,Jennings2011} was developed for modeling internal dynamics within an internal ward. However, in the ED, beds are not constantly occupied and the utilization level depends on the flow of patients that arrive from outside the system.
Yom-Tov \& Mandelbaum~\cite{YomTov2014} highlight the interrupted services while accounting for the transient nature of patient's arrival process, and introduced the Erlang-R model as a model for an ED. The Erlang-R model is an open two-station queueing network that has the same layout as the restricted Erlang-R model, except that all patients find a bed available upon arrival, see Figure \ref{fig:ErlangR}. In both models patients experience the interrupted services, but the Erlang-R model has no further restrictions on the bed capacity, hence neglecting the finite-size effects. Yom-Tov \& Mandelbaum~\cite{YomTov2014} showed, using a simulator tailored to an Israeli ED, that the complicated small ED dynamics can be captured using the relatively simple Erlang-R model, and hence, its recommendations can be implemented in ED workforce management.
Although the feature of interrupted services is present in many systems, it is particularly important for modeling EDs, because the duration of the interruption is typically much longer than the time patients require care from a nurse. This explains why the Erlang-R model is considered to be the canonical model for EDs. The restricted Erlang-R model with holding/blocking thus extends the Erlang-R model with finite-size constraints which, like interrupted services, are expected to have a decisive impact on performance.
 \\
 \\*
 \noindent
\textbf{Quality-and-Efficiency-Driven regime.} 
The Quality-and-Efficiency driven (QED) regime for many-server systems, also known as the Halfin-Whitt regime, adheres to a square-root staffing rule, which is best explained for the Erlang-C model. This system can be characterized completely by the staffing level $s$ and the offered load $R=\lambda/\mu$, which is the average workload pressure per time unit on the system. Moreover, let $\rho=R/s$ denote the server utilization level. Since exact analysis provides little qualitative insight in performance with respect to the system parameters, we resort to asymptotic analysis. This in turn provides approximations for the true system behavior. In the QED regime \cite{Borst2004,HalfinWhitt1981}, the utilization level is driven to unity in accordance to $(1-\rho)\sqrt{s} \rightarrow \beta$ as $s \rightarrow \infty$, for some fixed parameter $\beta >0$. This gives rise to the square-root dimensioning rule $s=R+\beta \sqrt{R}$, which prescribes that the number of nurses $s$ exceeds the minimally required offered load $R$, but only by a relatively small amount $\beta \sqrt{R}$. As $s$ grows large, the probability of delay tends to a non-degenerate function that only depends on the parameter $\beta$. This function is strictly decreasing in $\beta>0$ with range $(0,1)$. Consequently, any targeted delay probability can be achieved by adjusting $\beta$. Moreover, the mean delay is of order $1/\sqrt{s}$ and hence asymptotically negligible.
In this paper, we take the same approach, but determine capacity for nurses and beds simultaneously in such a way that both the probability to wait for a nurse and the probability to wait for a bed are non-degenerate.
Moreover, the utilization of both nurses and beds goes to unity as the size of the system increases. 
While the QED regime gives precise limits when the system size ($s$ and $n$ in our case) goes to infinity, it is by now well known that the asymptotic behavior kicks in quickly, so that QED limits serve as sharp approximations, already for small systems. See \cite{Janssen2011,erlanga,Leeuwaarden2011,Leeuwaarden2012,Gamarnik2013,Sanders2016} for various works that provide theoretical support for this fast relaxation.

\noindent 
In recent years, a new approach for developing approximations of performance measures of queueing systems  with finite-size constraints in the QED regime was suggested by van Leeuwaarden et al.~\cite{Leeuwaarden2015,Leeuwaarden2016}. This approach characterizes the asymptotic dynamics of these analytically intractable systems with finite-size restrictions through more tractable ones. We adopt this approach for developing the approximations for the holding model in \S\ref{sec:QED_limit_holding}.


\section{Models and performance measures}
\label{sec:modeldescription}

\subsection{Three-dimensional Markov process}
\label{sec:Markov_process}

Since in the restricted Erlang-R model described the arrival process is taken Poisson, and all service and content times are assumed independent and exponential, the system can be characterized in terms of a Markov process.
Let $Q(t) = (H(t),Q_1(t),Q_2(t))$ represent the number of patients in the \emph{holding}, \emph{needy} and \emph{content} state at time $t$, respectively. 
In both variants, $n$ is the maximum number of patients admitted to system, we have $Q_1(t)+ Q_2(t)\leq n$ for all $t\geq 0$. 
Due to the absence of holding patients in the Erlang-R model with blocking,  $H(t)=0$ is enforced in this case, whereas $H(t)$ has unbounded support in the model with holding. 
This distinction requires us to explore the stationary distribution of the two variants separately. 
Before doing so, we introduce some additional notation. 
We define 
\begin{equation}
 R_1 := \frac{\l}{(1-p)\mu}\, \qquad R_2 := \frac{p\l}{(1-p)\d}, 
 \label{eq:R1_R2}
 \end{equation}
where $R_1$ and $R_2$ can be interpreted as the offered workload brought towards the needy queue, and the content (infinite-server) queue, respectively. 
Furthermore, we define 
\begin{equation}
r:= \frac{\d}{\d+p\m},
\label{eq:delta}
\end{equation} 
which is the fraction of time a patient spends in the needy state (in case she experienced no wait during her sojourn). 

\subsubsection{Erlang-R model with blocking.}

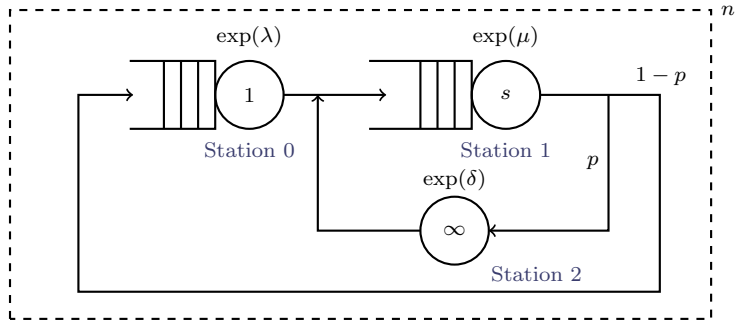
\begin{figure}[t]
\centering
\begin{tikzpicture}[scale = 0.9]
\draw [thick] (-1.25,5) -- (0,5) -- (0,4) -- (-1.25,4);
\draw [thick] (0.5,4.5) circle [radius = 0.5] node {\footnotesize 1} node[above=0.5cm] {\footnotesize exp$(\l)$} 
	node[below =0.5cm] {\footnotesize \color{col1} Station 0} ;
\draw [thick] (-0.25,4) -- (-0.25,5);
\draw [thick] (-0.5,4) -- (-0.5,5);
\draw [thick] (-0.75,4) -- (-0.75,5);
\draw [thick, ->] (1,4.5) -- (2.5,4.5);
\draw [thick] (2.25,5) -- (3.75,5) -- (3.75,4) -- (2.25,4);
\draw [thick] (3,4) -- (3,5);
\draw [thick] (3.25,4) -- (3.25,5);
\draw [thick] (3.5,4) -- (3.5,5);
\draw [thick] (4.25,4.5) circle [radius=0.5] node {\footnotesize $s$} node[above=0.5cm] {\footnotesize exp$(\mu)$}
node[below = 0.5cm] {\footnotesize \color{col1} Station 1} ;
\draw [thick, ->] (4.75,4.5) -- node[right=0.8cm,above] {\footnotesize $1-p$} (6.5,4.5) -- (6.5,1.6) -- (-2,1.6) -- (-2,4.5) -- (-1.2,4.5);
\draw [thick,->] (5.75,4.5) -- node[left] {\footnotesize $p$} (5.75,2.5) -- (4,2.5);
\draw [thick] (3.5,2.5) circle [radius=0.5] node {\footnotesize $\infty$} node[above=0.4cm] {\footnotesize exp$(\d)$}
node[below right = 0.35cm] {\footnotesize \color{col1} Station 2} ;;
\draw [thick,->] (3,2.5) -- (1.5,2.5) -- (1.5,4.5);

\draw [thick, dashed] (-3,1.2) rectangle (7.25,5.75) node[right] {\footnotesize $n$};
\end{tikzpicture}
\caption{The Erlang-R model with blocking, viewed as a closed Jackson network.}
\label{fig:ErlangR_blocking}
\end{figure}

In case of the blocking model, $Q(t)$ reduces to a finite-state Markov process $Q(t) = (Q_1(t),Q_2(t))$, where $Q_1(t)+Q_2(t)\leq n$ for all $t\geq 0$. 
In fact, this is equivalent to the closed Jackson network depicted in Figure \ref{fig:ErlangR_blocking} with finite population $n$. 
Station 1 in Figure \ref{fig:ErlangR_blocking} is an $M/M/s$ queue with service rate $\mu$, modeling the number of needy patients, $Q_1(t)$. 
Station 2 models the number of content patients, $Q_2(t)$, and can therefore be represented as an infinite-server queue with service rate $\d$. 
A patient can enter the unit only if $Q_1(t)+Q_2(t)<n$. 
Station 0---a single-server queue---moderates this as it only produces output at rate $\l$ in case its queue length is positive, i.e.\ if $n-Q_1(t)-Q_2(t)>0$.

Observe that because patients finding a full network are blocked, the number of patients in the system cannot grow beyond $n$. 
Hence, the system is stable for all parameter settings, and hence a steady-state distribution exists. Moreover, the simplification of the model with blocking allows us to express the steady-state distribution of the system in explicit product-form. 
Let $\pi_b(j,k)$ denote the steady-state probabilities of having $j$ needy and $k$ content patients in the system. Then, 
\begin{equation}\label{eq:pih(i,j)}
\pi_b(j,k) = \left\{
\begin{array}{ll}
\pi_0\,\frac{1}{\kappa(j)}\,\frac{1}{k!}\cdot R_1^j\cdot R_2^k, & ~~~\text{if }j+k \leq n,\\
0, & ~~~\text{else,}
\end{array}\right.
\end{equation}
where 
\begin{equation*}
\kappa(j) := \left\{
\begin{array}{ll}
j! , & ~~\text{if }j \leq s,\\
s!\, s^{j-s}, &~~ \text{else,}
\end{array}\right.
\end{equation*}
and $
\pi_0^{-1} = \sum_{j+k\leq n} \frac{1}{\kappa(j)}\,\frac{1}{k!}\cdot R_1^j\cdot R_2^k.
$

\subsubsection{Erlang-R with holding.}\label{ref:modelsoft}
The Erlang-R model with holding does not lead to a Jackson network with an elegant product-form solution for the steady-state distribution, because the holding queue cannot be modeled as a station that is independent from the other queues in the system.
However, we are able to describe the system as a two-dimensional Markov process without loss of information. 
To see this, define $N:= \{N(t),\,t\geq 0\}$ with $N(t) := H(t)+Q_1(t) + Q_2(t)$, the total number of patients in the system (including the holding queue).  
Using the restriction $Q_1(t)+Q_2(t) \leq n$ together with the fact that no bed is left vacant if a patient is waiting in the holding queue, this yields
\begin{equation*}
H(t) = \left(N(t) - n\right)^+, \quad t\geq 0,
\end{equation*}
where $(\cdot)^+ := \max\{0,\cdot\}$.
For the same reason, $Q_2(t) = N(t) - Q_1(t)$ if $H(t)=0$, and $Q_2(t) = n-Q_1(t)$ otherwise. 
In other words,
\begin{equation*}
Q_2(t) = \min\{N(t),n\} - Q_1(t), \quad t \geq 0.
\end{equation*}
Therefore, we can express the state of all three queues in the Erlang-R model with holding using a two-dimensional Markov process $X:= \{X(t),\,t\geq 0\}$, where 
\begin{equation*}
X(t) :=\left( N(t), Q_1(t) \right).
\end{equation*}
The process $X$ lives on the semi-infinite strip
 \begin{equation*}
X(t) \in \left\{\,(i,j)\, |\, j \leq \min\{i,n\}, i\in \mathbb{N}_0, j \in \{0,1,\ldots,n\}\, \right\},
\end{equation*}
and belongs to the class of Quasi-Birth-Death (QBD) processes.
The reader is referred to Appendix~\ref{app:QBDdescription} in the e-companion for a detailed description of this process, in terms of its transition diagram and generator matrix.

Contrary to the model with blocking, the system with holding \emph{can} become unstable in case capacity is insufficient to satisfy patient demand. 

\begin{proposition}\label{prop:StabilityCondition}
The Erlang-R model with holding is stable if and only if 
\begin{equation}
\frac{\lambda}{(1-p)\mu s} < \frac{ \sum_{i=0}^s \frac{i}{s}\binom{n}{i} \left(\frac{\d}{p\mu}\right)^i + \sum_{i=s+1}^n \binom{n}{i} \frac{i!}{s!} s^{s-i} \left(\frac{\d}{p\m}\right)^i}
{ \sum_{i=0}^s \binom{n}{i} \left(\frac{\d}{p\mu}\right)^i + \sum_{i=s+1}^n \binom{n}{i} \frac{i!}{s!} s^{s-i} \left(\frac{\d}{p\m}\right)^i}
=: \rho_{\max}(s,n).
\label{eq:StabilityCondition} 
\end{equation}
\end{proposition}
The proof is given in Appendix~\ref{app:stability} and follows from the general theory for QBD processes. 

\begin{figure}
\centering
\begin{subfigure}{0.48\textwidth}
\centering
\begin{tikzpicture}[scale = 0.75]
\begin{axis}[
	xmin = 0,
	xmax = 24,
	ymin = 0,
	ymax = 20,	
	grid = none, 
	axis line style={->},
	axis lines = left,
	xlabel = $s$,
	ylabel = {$R_{\rm max}(s,n)$},
	legend cell align=left,
	legend pos = north west
]

\addplot[col1,thick,mark=*,mark repeat = 2] table[x=s,y=n20] {r025_n_fixed.txt};
\addplot[col1,opacity = 0.75,thick,mark=x,mark repeat = 2] table[x=s,y=n40] {r025_n_fixed.txt};
\addplot[col1,opacity = 0.5,thick,mark=o,mark repeat = 2] table[x=s,y=n60] {r025_n_fixed.txt};
\addplot[col1,opacity = 0.25,thick,mark=square,mark repeat = 2] table[x=s,y=n80] {r025_n_fixed.txt};

\legend{$n=20$,$n=40$,$n=60$,$n=80$}
\end{axis}
\end{tikzpicture}
\caption{$R_{\rm max}$ as a function of $s$.}
\label{fig:Rmax_a}
\end{subfigure}
\begin{subfigure}{0.48\textwidth}
\centering
\begin{tikzpicture}[scale = 0.75]
\begin{axis}[
	xmin = 0,
	xmax = 100,
	ymin = 0,
	ymax = 20,	
	grid = none, 
	axis line style={->},
	axis lines = left,
	xlabel = $n$,
	ylabel = {$R_{\rm max}(s,n)$},
	legend cell align=left,
	legend pos = north west
]

\addplot[col1,thick,mark=*,mark repeat = 6] table[x=n,y=s5] {r025_s_fixed.txt};
\addplot[col1,opacity=0.75,thick,mark=x,mark repeat = 6] table[x=n,y=s10] {r025_s_fixed.txt};
\addplot[col1,opacity=0.5,thick,mark=o,mark repeat = 6] table[x=n,y=s15] {r025_s_fixed.txt};
\addplot[col1,opacity=0.25,thick,mark=square,mark repeat = 6] table[x=n,y=s20] {r025_s_fixed.txt};

\legend{$s=5$,$s=10$,$s=15$,$s=20$}
\end{axis}
\end{tikzpicture}
\caption{$R_{\rm max}$ as a function of $n$.}
\label{fig:Rmax_b}
\end{subfigure}
\caption{The maximum achievable workload in the restricted Erlang-R model with holding for $r=0.25$.}
\label{fig:Rmax}
\end{figure}
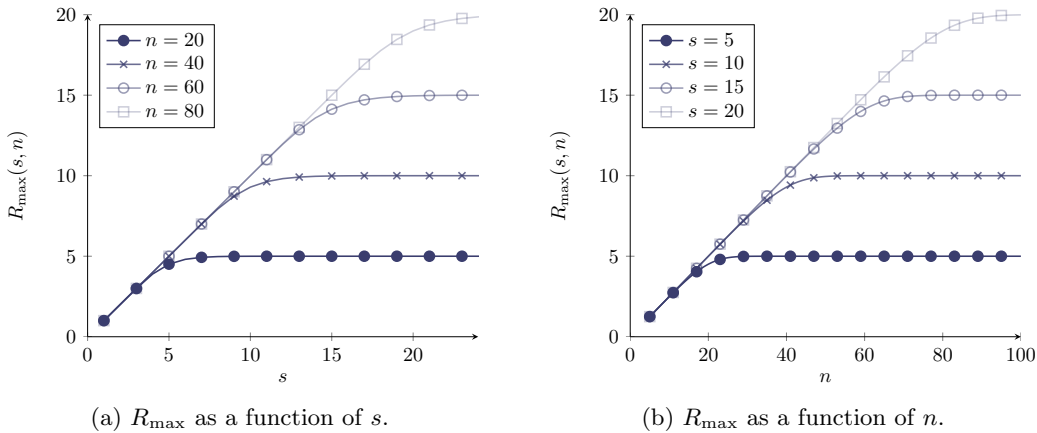

Observe that $\rho_{\max}(s,n)$ poses an upper bound on the occupancy level of the servers in the holding model, which is clearly smaller than 1 for all $s$ and $n$. 
In addition, this implies that the maximum workload $R_{\max}(s,n) := s\cdot\rho_{\max}(s,n)$ the system is able to handle is strictly less than $s$.
If we compare this to the open Erlang-R model, in which the maximal attainable workload equals $s$, we observe the effect of finite-size constraints on operational performance. 
Figure \ref{fig:Rmax} shows the influence of both $s$ and $n$ on the maximum feasible workload in case $r=0.25$. 
From these graphs, note that if $s\ll rn$, $R_{\max}$ grows almost linearly with $s$. 
Furthermore, $R_{\rm max}(s,n)$ is increasing in $n$ for $s$ fixed.
A logical practical consequence is that a larger number of beds allows for a larger patient volume to enter the ED with the same number of nurses.
Moreover, $R_{\rm max}(s,n)$ is increasing in $s$, but as in Figure \ref{fig:Rmax_a}, adding an extra nurse does not increase the stability region in case $n$ is too tight. 
 Conversely, adding extra beds does not increase $R_{\rm max}(s,n)$ if the number of nurses does not allow for an increase in offered load, see Figure \ref{fig:Rmax_b}. 
Additionally, it is easily verified that $R_{\rm max}(s,n)$ is upper bounded by both $s$ and $R_{\rm max}(n,n) = rn$. Therefore, a careful balance is called for between servers (nurses) and beds, so that resources will be efficiently utilized. We observe that when the ratio $s/n\approx r$, the system is better balanced. 
We will propose an appropriate balance between resources by defining a synchronized QED capacity recommendation for both servers and beds in \S \ref{sec:QED_scaling}. 

\begin{insight}
Both the number of nurses $s$ and the number of beds $n$ play a restricting role on the maximum demand the system can handle. Therefore, they should be balanced together.
\end{insight} 

Provided that the system is stable, the stationary distribution of the QBD process $X$ can be obtained numerically by the matrix geometric method \cite{Neuts1981}.
Subsequently, we can derive the stationary distribution of the original $Q(t)$, denoted by $\pi_h(\cdot,\cdot,\cdot)$.

\subsection{Performance measures}
\label{sec:performance_metrics}
In this work, we concentrate on five performance measures that are central to our analysis. 
In the definitions that follow, we present expressions for these measures in terms of a general three-dimensional measure $\pi$, which one can replace by either $\pi_b$ or $\pi_h$, depending on the scenario considered. 
In the remainder of this work, we will augment the measures related to the Erlang-R model with blocking and holding by the superscript $b$ and $h$, respectively\footnote{In line with $H(t)=0$, we use $\pi_b(i,j,k) = \pi_b(j,k)$ if $i=0$, with $\pi_b(j,k)$ as in \eqref{eq:pih(i,j)}, and $\pi_b(i,j,k) = 0$ otherwise, when considering the model with blocking.}.

As relevant performance measures, we consider the probability of holding (blocking) at entering the system, the probability of delay at the needy queue, expected waiting time for a nurse, utilization of nurses and utilization of beds:
\begin{equation}
\P({\rm hold}) = \sum_{i=0}^\iy \sum_{j=0}^n \pi(i,j,n-j), \qquad
\P({\rm delay}) \approx \sum_{i=0}^\infty\sum_{j=s}^{n}\sum_{k=0}^{n-j} \pi(i,j,k), 
\label{eq:delay_probability}
\end{equation}
\begin{equation}
\label{eq:EW_exact}
\E [W] \approx \sum_{i=0}^\infty\sum_{j=s}^{n}\sum_{k=0}^{n-j} \frac{\max\{0,j-s+1\}}{\mu}\,\pi(i,j,k), 
\end{equation}
\begin{equation}
\label{eq:utilization}
\rho_s = \frac{1}{s}\,\sum_{i=0}^\infty \sum_{j=0}^n \sum_{k=0}^{n-j} \min\{j,s\} \pi(i,j,k), \qquad
\rho_n = \frac{1}{n}\,\sum_{i=0}^\infty \sum_{j=0}^n \sum_{k=0}^{n-j} \min\{i,n\} \pi(i,j,k).
\end{equation}

It should be stressed that the above expression for the delay probability and the expected waiting time for a nurse is not exact. For the blocking model one can use the Arrival Theorem \cite{ChanYao}, whereby the exact expression uses $n-1$ instead of $n$; for the holding model, the arrival process to the needy queue, which consists of both external arrivals and content patients becoming needy, is not Poisson. Therefore, we cannot use the PASTA argument for the holding model. However, for both models,  as we will be studying the system as $s$ and $n$ become large, this approximation error will become negligible. 

\subsection{Stochastic bounds}
\label{sec:bounds}


Although the two variants of the Erlang-R model differ with respect to the admission policy, and require different mathematical treatment, we would like to be able to capture their relative performance.  
We substantiate the intuition that the holding room leads to more patients in the ED, in the following result. 

\begin{proposition}\label{thm:stochasticordering}
Let $Q_1^b$, $Q_2^b$, $Q_1^h$, $Q_2^h$ denote the nurse and content queue length processes in the Erlang-R model with blocking and holding, respectively.
Let $H(0) = 0$,  $Q_1^b(0)=Q_1^h(0)$ and $Q_2^b(0)=Q_2^b(0)$. For all $t\geq 0$, 
\begin{align}
Q_1^b(t) + Q_2^b(t) &\preceq_{\rm st} Q_1^h(t) + Q_2^h(t) \preceq_{\rm st} n ,\\
Q_2^b(t) &\preceq_{\rm st} Q_2^h(t),\\
Q_1^b(t) &\preceq_{\rm st} Q_1^h(t) + H(t),
\end{align}
where $X\preceq_{\rm st} Y$ implies $\P(X\geq k) \leq \P(Y\geq k)$ for all $k\geq 0$.
\end{proposition}

The proof of Proposition \ref{thm:stochasticordering} uses sample path coupling and can be found in Appendix \ref{app:stochastic_ordering}.
Note that as an immediate consequence, we have
\[ \P^b( {\rm hold }) = \lim_{t\to\iy} \P\left( Q_1^b(t)+Q_2^b(t) \geq n \right) \leq \lim_{t\to\iy} \P\left( Q_1^h(t) + Q_2^h(t) \geq n \right) = \P^h( {\rm hold }) \]
and by similar reasoning $\rho^b_n \leq  \rho_n^h$.
\begin{insight}
\label{in:order}
Under similar offered load and capacity constraints, utilization levels for the nurses in the Erlang-R model with blocking are lower than in the Erlang-R model with holding. Moreover, the total number of waiting patients in the setting with holding is stochastically larger than in the setting with blocking, and in the open Erlang-R model. 
\end{insight}
\noindent
We further discuss the differences between both models in \S\ref{sec:dimensioning} and \S\ref{sec:analysis}.

\section{Two-fold QED regime}
\label{sec:QED_scaling}

We do not want to waste capacity of either servers or beds without getting significant advantage in term of performance. 
We therefore take an asymptotic approach that lets the external arrival rate $\l$ grow to infinity, while scaling $s$ and $n$ accordingly.
In doing so, we intend to establish QED-type system behavior, i.e.\ high occupancy levels of both nurses and beds and good quality-of-service.

\subsection{Two-fold scaling rule}

In order to identify the scaling of $s$ and $n$ as $\l\to\infty$, we draw inspiration from the two-fold scaling rule in \cite{Jennings2008} and \cite{Khudyakov2010}, which follows the celebrated square-root staffing principle discussed in \S 2.
This principle suggests that, in the most general setting, capacity should be equal to the expected offered load entering the system, let us say $R$, plus an additional variability hedge that is proportional to $\sqrt{R}$.
In the restricted Erlang-R model, we have two capacity sources, namely $s$ and $n$, which experience different relevant amount of works.

The offered load the servers in the needy queue experience is given by $R_1$, as in the regular Erlang-R model; 
it does not change due to the finite-size effects, since all patients are served eventually. Hence, we only need to account for the interrupted services. It follows that the appropriate staffing rule for the nurses in the QED regime remains $s=R_1+\beta \sqrt{R_1}$ for some constant $\beta >0$.

To establish the bed capacity level, we need to reflect on the load offered to the beds. Observe that beds remain occupied both in needy and content states. This suggests that $R_n:=R_1+R_2=R_1/r$, with $R_1$ and $R_2$ as in \eqref{eq:R1_R2} and $r$ the expected fraction of time a patient spends at the nurse station, defined in \eqref{eq:delta}.
%
As a result, the appropriate staffing rule is $n=R_n+\gamma \sqrt{R_n}$ for some constant $\gamma>0$. In conclusion, the two-fold QED scaling rule is given by
\begin{equation}\label{eq:twofoldscaling}
\begin{array}{ll}
s &= R_1 + \beta \sqrt{R_1} + o(\sqrt{R_1}) \\
n &= \frac{R_1}{r}+\gamma \sqrt{\frac{R_1}{r}} + o(\sqrt{R_1})
\end{array}
\end{equation}
with $\beta,\gamma>0$ constants and $R_1:=\lambda/((1-p)\mu)$.

Recall that we saw in Figure \ref{fig:Rmax} that resources seem efficiently utilized if $s/n\approx r$. 
Scaling \eqref{eq:twofoldscaling} is in line with this reasoning since 
\[
\frac{s}{n} = r\left(1+ \frac{\beta - \gamma\sqrt{r}}{\sqrt{R_1}}+ O(1/R_1) \right) . 
\]

\begin{remark}
In \cite{Jennings2008}, a similar scaling regime is considered, which only relates $s$ and $n$ through a square-root scaling, namely the regime $s = r n + \hat\g\sqrt{n}$,
which is equivalent to the second relation in \eqref{eq:twofoldscaling} if $\hat\g = \b\sqrt{r} - \g r$.
Due to the absence of external arrivals in this closed system, they let the number of beds $n$ approach infinity as opposed to $\l$ in our settings. 
Nevertheless, this results in the same asymptotic regime.
\end{remark}

Before turning to asymptotic expressions for the performance measures concerning the Erlang-R model with blocking/holding, we conduct a few numerical experiments to confirm that the scaling in \eqref{eq:twofoldscaling} indeed leads to desired QED behavior. 

In Figure \ref{fig:sample_paths} we plotted the sample paths of the three-dimensional queue length process of the holding model in which $\b$ and $\g$ are fixed, and $R_1$ is increased. 
Observe that the needy queue length $Q_1(t)$, plotted in orange in Figure \ref{fig:sample_paths}, fluctuates around the values $s$, and stabilizes for larger values of $R_1$. 
This naturally implies that the server (nurses) utilization approaches 100\%, while the number of patients waiting is $O(\sqrt{R_1})$.
Furthermore, we see that the percentage of occupied beds also tends to 100\%, while the holding queue length remains small. 
The holding queue is of much smaller order than $R_1$, which implies that the holding time of a patient becomes negligible as $R_1\to\iy$. 
From these empirical findings we deduce that  under scaling \eqref{eq:twofoldscaling} the restricted Erlang-R model exhibits QED behavior on two levels: Outside the facility while waiting for an available bed, and inside the facility while waiting for attention of a nurse.

\begin{figure}[h] \centering
\centering
\begin{subfigure}{0.32\textwidth}
\begin{tikzpicture}[scale=0.59]

\begin{axis}[
	xmin = 0,
	xmax = 200,
	ymin = 0.0,
	ymax = 28,
	ytick = {0,5,10,15,20,25},	
	grid = none, 
	axis line style={->},
	axis lines = left,
	xlabel = $t$,
	legend cell align=left,
	legend pos = north east
]

\addplot[very thick,col1] file {R5_holding.txt};
\addplot[very thick,col5] file {R5_service.txt};
\addplot[very thick,col3] file {R5_total.txt};
\addplot[very thick,dashed] coordinates {
	(0,7)
	(200,7)
	};
\addplot[very thick,dashed] coordinates {
	(0,24)
	(200,24)
	};

\end{axis}

\end{tikzpicture}
\caption{$R_1=5$}
\end{subfigure}
\begin{subfigure}{0.32\textwidth}
\begin{tikzpicture}[scale=0.59]

\begin{axis}[
	xmin = 0,
	xmax = 200,
	ymin = 0.0,
	ymax = 128.333,
	ytick = {0,20,40,60,80,100,120},	
	grid = none, 
	axis line style={->},
	axis lines = left,
	xlabel = $t$,
	legend cell align=left,
	legend pos = north east
]

\addplot[very thick,col1] file {R25_holding.txt};
\addplot[very thick,col5] file {R25_service.txt};
\addplot[very thick,col3] file {R25_total.txt};
\addplot[very thick,dashed] coordinates {
	(0,30)
	(200,30)
	};
\addplot[very thick,dashed] coordinates {
	(0,110)
	(200,110)
	};
\end{axis}

\end{tikzpicture}
\caption{$R_1=25$}
\end{subfigure}
\begin{subfigure}{0.32\textwidth}
\begin{tikzpicture}[scale=0.59]

\begin{axis}[
	xmin = 0,
	xmax = 200,
	ymin = 0.0,
	ymax = 490,
	ytick = {0,100,200,300,400},	
	grid = none, 
	axis line style={->},
	axis lines = left,
	xlabel = $t$,
	legend cell align=left,
	legend pos = north east
]
 	
\addplot[very thick,col1] file {R100_holding.txt};
\addplot[very thick,col5] file {R100_service.txt};
\addplot[very thick,col3] file {R100_total.txt};
\addplot[very thick,dashed] coordinates {
	(0,110)
	(200,110)
	};
\addplot[very thick,dashed] coordinates {
	(0,420)
	(200,420)
	};
	
\end{axis}

\end{tikzpicture}
\caption{$R_1=100$}
\end{subfigure}
\caption{Sample paths of $H(t)$ (blue), $Q_1(t)$ (orange) and $Q_1(t)+Q_2(t)$ (green) of the Erlang-R model with holding with parameters $\mu = 1$, $\d=0.25$, $p=0.75$ and $\b=\g=1$. The staffing levels $s$ and $n$ are depicted by the dashed lines.}
\label{fig:sample_paths}
\end{figure}
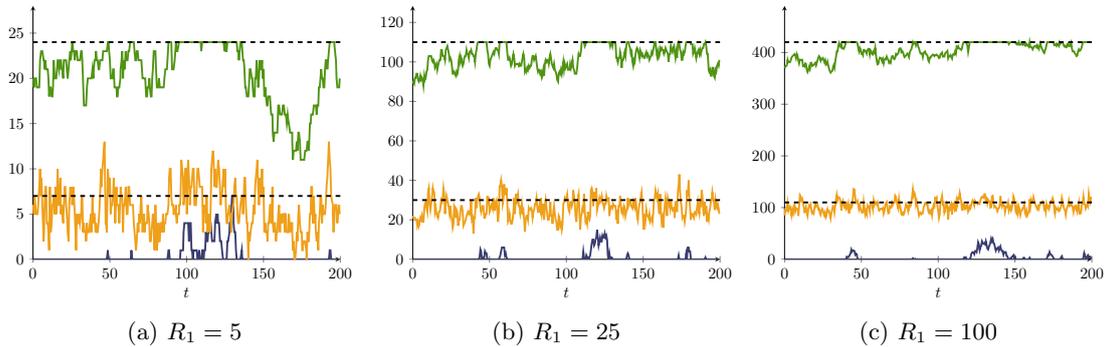

\begin{figure}[h] \centering
\centering
\begin{subfigure}{0.32\textwidth}

\tikzset{mark size=3}
\begin{tikzpicture}[scale=0.6]

\begin{axis}[
	xmin = 0,
	xmax = 145,
	ymin = 0.0,
	ymax = 0.7,
	ytick = {0,0.1,...,0.7},
	xlabel = $\l$,
	grid = both, 
	axis line style={->},
	axis lines = left,
	legend cell align=left,
	legend pos = south east
]

\addplot[very thick,col1,mark=*] file {delayProbErlangH.txt};
\addplot[very thick,col1,opacity=0.75,mark=x] file {delayProbYomTov.txt};
\addplot[very thick,col1,opacity=0.5,mark=square] file {delayProbJennings.txt};
\small
\legend{Holding,Blocking,Closed ward};
\end{axis}

\end{tikzpicture}
\caption{Delay probability nurse}
\end{subfigure}
\begin{subfigure}{0.32\textwidth}
\tikzset{mark size=3}
\begin{tikzpicture}[scale=0.6]

\begin{axis}[
	xmin = 0,
	xmax = 145,
	ymin = 0.0,
	ymax = 0.3,
	ytick = {0,{0.05},0.1,0.15,0.2,0.25,3},	
	grid = both, 
	axis line style={->},
	axis lines = left,
	xlabel = $\l$,
	legend cell align=left,
	legend pos = north east,
	y tick label style={
        /pgf/number format/.cd,
        fixed,
        fixed zerofill,
        precision=2,
        /tikz/.cd
    }
]

\addplot[very thick,col1,mark=*] file {EWErlangH.txt};
\addplot[very thick,col1,opacity=0.75,mark=x] file {EWYomTov.txt};
\addplot[very thick,col1,opacity=0.5,mark=square] file {EWJennings.txt};
\small
\legend{Holding,Blocking,Closed ward};
\end{axis}

\end{tikzpicture}
\caption{Expected wait}
\end{subfigure}
\begin{subfigure}{0.32\textwidth}
\tikzset{mark size=3}
\begin{tikzpicture}[scale=0.6]

\begin{axis}[
	xmin = 0,
	xmax = 145,
	ymin = 0.7,
	ymax = 1.02,
	grid = both, 
	axis line style={->},
	axis lines = left,
	xlabel = $\l$,
	legend cell align=left,
	legend pos = south east
]

\addplot[very thick,col1,mark=*] file {rhoErlangH.txt};
\addplot[very thick,col1,opacity=0.75,mark=x] file {rhoYomTov.txt};
\addplot[very thick,col1,opacity=0.5,mark=square] file {rhoJennings.txt};
\small
\legend{Holding,Blocking,Closed ward};
\end{axis}

\end{tikzpicture}
\caption{Nurse utilization}
\end{subfigure}
\caption{Asymptotic behavior of the restricted Erlang-R model with holding and blocking, and the Closed ward model for $\m=1$, $\d = 0.2$, $p=0.8$ and $\b=\g=0.5$.}
\label{fig:empiricalAsymptotics}
\end{figure}
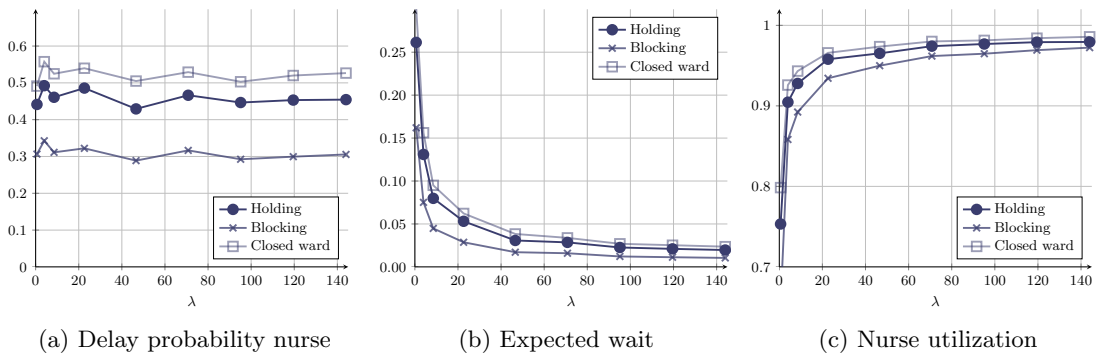

\tikzset{mark size=2.5}

We also check how the Erlang-R model with blocking or holding and the closed ward model of \cite{Jennings2008} relate under scaling \eqref{eq:twofoldscaling}. 
In Figure~\ref{fig:empiricalAsymptotics}, we plot the performance measures, obtained through simulation, for the three models in which we fix $\b=\g=0.5$ and vary the arrival rate $\l$. 
First, we see that $\P({\rm delay})$ stabilizes as $\l\to\iy$ in all three models under scaling \eqref{eq:twofoldscaling}, and the delay probability of the model with holding lies in between the other two. 
Second, note that the expected waiting time for a nurse in all models converges to 0 as $\l$ increases. 
In fact, the rate of decay is similar in all three models.  
We observe that $\rho_s$ approaches unity in all three models, and the rate of convergence seems again comparable. 
Finally, and most importantly, we notice an ordering between the three models.
Namely, in all performance measures considered in Figure \ref{fig:empiricalAsymptotics}, Erlang-R with holding appears to be upper bounded by the closed ward and lower bounded by the Erlang-R with blocking. 
In a multitude of parameter settings of $(\b,\g)$, we have seen the same ordering, leading to the following conjecture: 
\begin{conjecture}\label{conj:stochorder}
Let $Q^b_1(\iy)$, $Q_1^h(\iy)$ and $Q_1^c(\iy)$ denote the stationary number of needy patients in the Erlang-R model with blocking, holding and the closed ward, respectively. Then,
\begin{equation}
Q_1^b(\iy) \preceq_{\rm st} Q_1^h(\iy) \preceq_{\rm st} Q_1^c(\iy).
\end{equation}
\end{conjecture}
Observe that Conjecture \ref{conj:stochorder} poses a stronger statement than the third assertion in Proposition \ref{thm:stochasticordering}.

\subsection{QED limits for Erlang-R with blocking}
\label{sec:QED_limit_block}
We now continue our analysis by examining the limiting behavior under scaling \eqref{eq:twofoldscaling}.
We first derive QED limits for some performance measures of the Erlang-R model with blocking.
Using the explicit expressions for the blocking model in \eqref{eq:pih(i,j)}, we derive the limiting values of the relevant performance measures defined in \S\ref{sec:performance_metrics} in terms of $\b$ and $\g$. 

\begin{theorem}\label{thm:limits_YT}
Let $s$ and $n$ scale as in \eqref{eq:twofoldscaling} with ${-}\infty<\b<\infty,\,\g>0$ as $\l\to\infty$. Then, if $\b \neq 0$,
\begin{align}
\label{eq:yt_limit_delay}
g^b(\b,\g)
&:= \lim_{\l\to\iy} \P^b({\rm delay}) = 
\left(1 + 
\frac{ \b \, \int_{-\iy}^\b \Phi\left(\frac{\g-t\sqrt{r}}{\sqrt{1-r}}\right)\, d\Phi(t) }
{\phi(\b)\Phi(\eta) -  \phi(\sqrt{\b^2+\eta^2}){\rm e}^{\tfrac{1}{2} \omega^2} \Phi(\omega)}
\right)^{-1},\\
\label{eq:yt_limit_block}
f^b(\b,\g) 
&:= \lim_{\l\to\iy} \sqrt{R_1}\cdot\P^b({\rm block}) = 
\frac{
\sqrt{r}\phi(\g)\Phi(-\omega\sqrt{r}) + \phi(\sqrt{\b^2+\eta^2})\,{\rm e}^{\frac{1}{2} \omega^2} \Phi(\omega)
}{
\int_{-\iy}^\b \Phi\left(\frac{\g-t\sqrt{r}}{\sqrt{1-r}}\right)\, d\Phi(t) +
\frac{\phi(\b)\Phi(\eta)}{\b} -  \frac{\phi(\sqrt{\b^2+\eta^2})}{\b}{\rm e}^{\tfrac{1}{2} \omega^2} \Phi(\omega)
 },\\
\label{eq:yt_limit_Edelay}
h^b(\b,\g) &:= \lim_{\l\to\iy} \sqrt{R_1}\cdot\E[W]
=
\frac{
\frac{\phi(\b)\Phi(\eta)}{\b^2} +
 \left(\frac{\b}{r}-\frac{\g}{\sqrt{r}}-\frac{1}{\b}\right)\,\frac{\phi(\sqrt{\eta^2+\b^2})}{\b}\, {\rm e}^{\tfrac{1}{2}\omega^2}\, \Phi(\omega) 
 - \sqrt{\frac{1-r}{r}}\,\frac{\phi(\b)\phi(\eta)}{\b}
}{
\int_{-\iy}^\b \Phi\left(\frac{\g-t\sqrt{r}}{\sqrt{1-r}}\right)\, d\Phi(t) +
\frac{\phi(\b)\Phi(\eta)}{\b} -  \frac{\phi(\sqrt{\b^2+\eta^2})}{\b}{\rm e}^{\tfrac{1}{2} \omega^2} \Phi(\omega)
 },
\end{align}
and if $\b=0$,
\begin{align}
\label{eq:yt_limit_delay_beta0}
g^b_0(\g) 
&:= \lim_{\l\to\iy} \P^b({\rm delay})=
\left(1+
\frac{
\int_{-\iy}^0 \Phi\left(\frac{\g-t\sqrt{r}}{\sqrt{1-r}}\right)\, d\Phi(t)
}{
\sqrt{\frac{1-r}{r}} \frac{1}{\sqrt{2\pi}}\,\left(\eta \,\Phi(\eta) + \phi(\eta) \right)
}
\right)^{-1},\\
\label{eq:yt_limit_block_beta0}
f^b_0(\g) 
&:= \lim_{\l\to\iy} \sqrt{R_1}\cdot\P^b({\rm block}) = 
\frac{
\sqrt{r}\,\phi(\g)\Phi(-\omega\sqrt{r}) + \frac{1}{\sqrt{2\pi}} \Phi(\eta)
}{
\int_{-\iy}^\b \Phi\left(\frac{\g-t\sqrt{r}}{\sqrt{1-r}}\right)\, d\Phi(t) +
\sqrt{\frac{1-r}{r}} \frac{1}{\sqrt{2\pi}}\,\left(\eta \,\Phi(\eta) + \phi(\eta) \right)
 },\\
\label{eq:yt_limit_Edelay_beta0}
h_0^b(\g) &:= \lim_{\l\to\iy} \sqrt{R_1}\cdot\E[W] 
= \frac{1}{2\mu}\, \frac{ \left( \gamma^2/r+1\right) \Phi(\eta) + \eta \phi(\eta) }
{ \frac{r}{1-r} \sqrt{2\pi} \int_{-\infty}^0 \Phi\left(\frac{\g-t\sqrt{r}}{\sqrt{1-r}}\right)\, d\Phi(t) + \sqrt{\frac{r}{1-r}} \left(\eta \Phi(\eta)+\phi(\eta)\right)},
\end{align}
where $\eta = \frac{\g - \b\sqrt{r}}{\sqrt{1-r}}$ and $\omega := \frac{\g - \b/\sqrt{r}}{\sqrt{1-r}}$.
\end{theorem}
The proof of Theorem \ref{thm:limits_YT} can be found in Appendix \ref{app:proof_block}.

Theorem \ref{thm:limits_YT} proves that the scaling \eqref{eq:twofoldscaling} results in QED behavior: the probability of waiting in Equations \eqref{eq:yt_limit_delay} and \eqref{eq:yt_limit_delay_beta0} converges to a limit that is strictly between 0 and 1. 
Notice that all limits in Theorem \ref{thm:limits_YT} are functions of three parameters: $\beta$ and $\gamma$, which are decision variables, and the fraction of needy time $r$, which is dictated by the physics of the system. 

Theorem \ref{thm:limits_YT} also shows that the probability of blocking (Equations \eqref{eq:yt_limit_block} and \eqref{eq:yt_limit_block_beta0}) is of order $1/\sqrt{R_1}$. 
For example, assume that the fraction of needy time $r$ is $0.5$ and the system is large (100 servers). 
Using Figure \ref{fig:pdelay_pblock}, we observe that, by choosing the pair $\gamma = 1$ and $\beta = 0.245$, we actually aim at a probability of getting served immediately to be 40\%. At the same time, the probability of getting immediately a bed is 97\%. Thus, our QED policy puts more emphasis on preventing blocking than waiting.
\pgfplotsset{scaled y ticks=false}
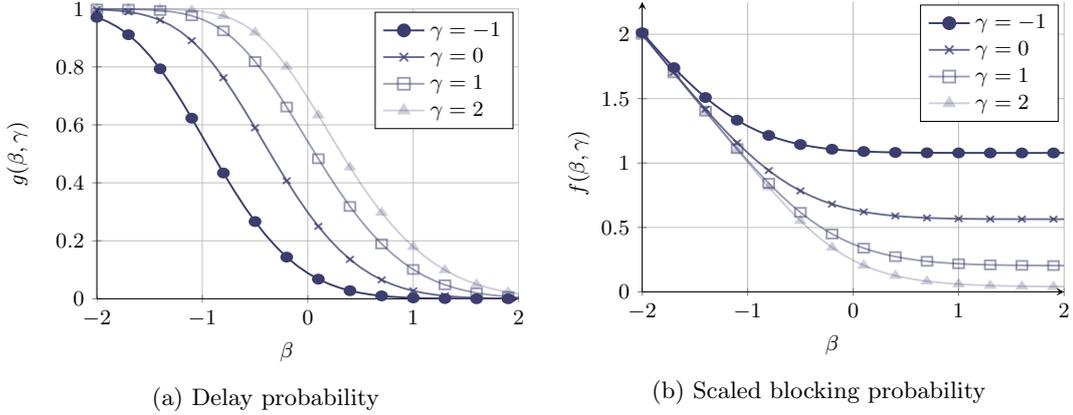
\begin{figure}
\centering
\begin{subfigure}{0.48\textwidth}
\centering
\begin{tikzpicture}[scale = 0.9]
\small
\begin{axis}[
	xmin = -2,
	xmax = 2,
	ymin = 0,
	ymax = 1,
	grid = both, 
	axis line style={->},
	axis lines = left,
	tick label style={/pgf/number format/fixed},
	xlabel = $\beta$,
	ylabel = {$g(\beta,\gamma)$},
	y label style = {at = {(axis cs: -1.95,0.5)}}, 
	xscale = 0.9,
	yscale = 0.75,
	legend cell align=left,
	legend style = {at = {(axis cs: 1.95,0.99)},anchor = north east}
]

\addplot[thick,col1,mark repeat = 3, mark = *] table[x=beta,y=g_min1] {limit_probabilities_delay.txt};
\addplot[thick,col1,opacity=0.75,mark repeat = 3, mark = x] table[x=beta,y=g_0] {limit_probabilities_delay.txt};
\addplot[thick,col1,opacity=0.5,mark repeat = 3, mark = square] table[x=beta,y=g_1] {limit_probabilities_delay.txt};
\addplot[thick,col1,opacity=0.25,mark repeat = 3, mark = triangle*] table[x=beta,y=g_2] {limit_probabilities_delay.txt};

\legend{{$\gamma = -1$},{$\gamma = 0$},{$\gamma = 1$},{$\gamma = 2$}};
\end{axis}
\end{tikzpicture}
\caption{Delay probability}
\label{fig:pdelay_pblock_a}
\end{subfigure}
\begin{subfigure}{0.48\textwidth}
\centering
\begin{tikzpicture}[scale = 0.9]
\small
\begin{axis}[
	xmin = -2,
	xmax = 2,
	ymin = 0,
	ymax = 2.25,
	grid = both, 
	axis line style={->},
	axis lines = left,
	xlabel = $\beta$,
	tick label style={/pgf/number format/fixed},
	ylabel = {$f(\beta,\gamma)$},
	y label style = {at = {(axis cs: -1.8,1)}}, 
	xscale = 0.9,
	yscale = 0.75,
	legend cell align=left,
	legend style = {at = {(axis cs: 1.95,2.23)},anchor = north east}
]
\addplot[thick,col1,mark repeat = 3, mark = *] table[x=beta,y=g_min1] {limit_probabilities_block.txt};
\addplot[thick,col1,opacity=0.75,mark repeat = 3, mark = x] table[x=beta,y=g_0] {limit_probabilities_block.txt};
\addplot[thick,col1,opacity=0.5,mark repeat = 3, mark = square] table[x=beta,y=g_1] {limit_probabilities_block.txt};
\addplot[thick,col1,opacity=0.25,mark repeat = 3, mark = triangle*] table[x=beta,y=g_2] {limit_probabilities_block.txt};


\legend{{$\gamma = -1$},{$\gamma = 0$},{$\gamma = 1$},{$\gamma = 2$}};
\end{axis}
\end{tikzpicture}
\caption{Scaled blocking probability}
\label{fig:pdelay_pblock_b}
\end{subfigure}
\caption{Asymptotic delay and scaled blocking probability for $r=0.5$ based on Theorem \ref{thm:limits_YT}. }
\label{fig:pdelay_pblock}
\end{figure}

Theorem \ref{thm:limits_YT} further shows that the expected waiting (Equations \eqref{eq:yt_limit_Edelay} and \eqref{eq:yt_limit_Edelay_beta0})
is of order $1/\sqrt{R_1}$ too and hence vanishes in the large-system limit.

We see from Theorem \ref{thm:limits_YT} that achieving target service levels is always an interplay between $\beta$ and $\gamma$.
Figure \ref{fig:pdelay_pblock_a} shows for instance that in order to keep $\P({\rm delay})\in (0.25,0.75)$,  choosing $\gamma=-1$ requires $\beta$ to stay within the range $[-1.4,-0.5]$, while $\gamma=1$ corresponds to values of $\beta$ in $[-0.4,0.5]$.

While the two-fold scaling rule in \eqref{eq:twofoldscaling} automatically captures the right dimensioning ratio as the system scales up, Theorem \ref{thm:limits_YT} shows that the parameters $\beta$ and $\gamma$ provide a means to fine-tune the performance.
Figure \ref{fig:pdelay_pblock_b} confirms how adding nurses, i.e.~increasing $\beta$, does not improve the blocking probability if the number of beds, i.e.~$\gamma$, is too tight.
This is in accordance with our previous observations in Figure \ref{fig:Rmax} for the exact steady-state distribution.

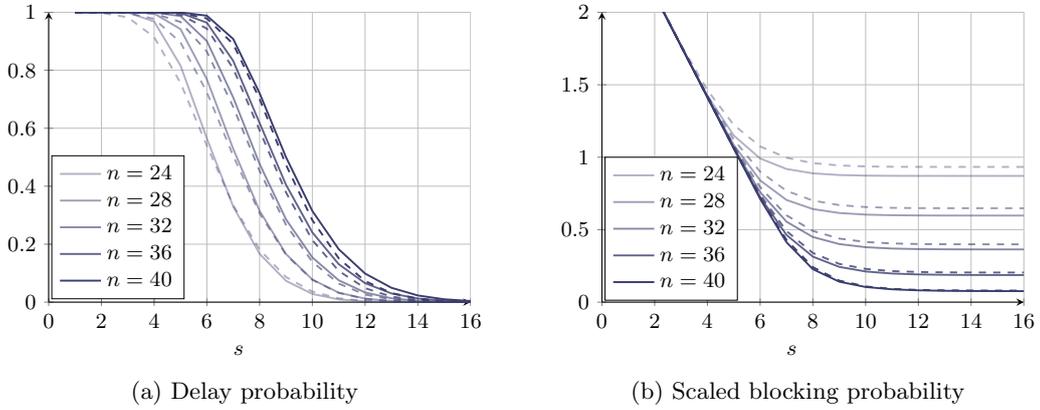
\begin{figure}
\centering
\begin{subfigure}{0.48\textwidth}
\centering
\begin{tikzpicture}[scale = 0.9]
\small
\begin{axis}[
	xmin = 0,
	xmax = 16,
	ymin = 0,
	ymax = 1,
	grid = both, 
	axis line style={->},
	axis lines = left,
	xlabel = $s$,
	xscale = 0.9,
	yscale = 0.75,
	legend cell align=left,
	legend style = {at = {(axis cs: 0.1,0.01)},anchor = south west}
]
 

\addplot[thick,col1,opacity=0.4] table[x=s,y=n24] {accuracy_pdelay_case2.txt};
\addplot[thick,col1,opacity=0.5] table[x=s,y=n28] {accuracy_pdelay_case2.txt};
\addplot[thick,col1,opacity=0.6] table[x=s,y=n32] {accuracy_pdelay_case2.txt};
\addplot[thick,col1,opacity=0.8] table[x=s,y=n36] {accuracy_pdelay_case2.txt};
\addplot[thick,col1,opacity=1] table[x=s,y=n40] {accuracy_pdelay_case2.txt};

\addplot[thick,col1,opacity=0.4,dashed] table[x=s,y=approx_n24] {accuracy_pdelay_case2.txt};
\addplot[thick,col1,opacity=0.5,dashed] table[x=s,y=approx_n28] {accuracy_pdelay_case2.txt};
\addplot[thick,col1,opacity=0.6,dashed] table[x=s,y=approx_n32] {accuracy_pdelay_case2.txt};
\addplot[thick,col1,opacity=0.8,dashed] table[x=s,y=approx_n36] {accuracy_pdelay_case2.txt};
\addplot[thick,col1,opacity=1,dashed] table[x=s,y=approx_n40] {accuracy_pdelay_case2.txt};

\legend{{$n=24$},{$n=28$},{$n=32$},{$n=36$},{$n=40$}};
\end{axis}
\end{tikzpicture}
\caption{Delay probability}
\end{subfigure}
\begin{subfigure}{0.48\textwidth}
\centering
\begin{tikzpicture}[scale = 0.9]
\small
\begin{axis}[
	xmin = 0,
	xmax = 16,
	ymin = 0,
	ymax = 2,
	grid = both, 
	axis line style={->},
	axis lines = left,
	xlabel = $s$,
	xscale = 0.9,
	yscale = 0.75,
	legend cell align=left,
	legend style = {at = {(axis cs: 0.1,0.01)},anchor = south west}
]

\addplot[thick,col1,opacity=0.4] table[x=s,y=n24] {accuracy_pblock_case2.txt};
\addplot[thick,col1,opacity=0.5] table[x=s,y=n28] {accuracy_pblock_case2.txt};
\addplot[thick,col1,opacity=0.6] table[x=s,y=n32] {accuracy_pblock_case2.txt};
\addplot[thick,col1,opacity=0.8] table[x=s,y=n36] {accuracy_pblock_case2.txt};
\addplot[thick,col1,opacity=1] table[x=s,y=n40] {accuracy_pblock_case2.txt};

\addplot[thick,col1,opacity=0.4,dashed] table[x=s,y=approx_n24] {accuracy_pblock_case2.txt};
\addplot[thick,col1,opacity=0.5,dashed] table[x=s,y=approx_n28] {accuracy_pblock_case2.txt};
\addplot[thick,col1,opacity=0.6,dashed] table[x=s,y=approx_n32] {accuracy_pblock_case2.txt};
\addplot[thick,col1,opacity=0.8,dashed] table[x=s,y=approx_n36] {accuracy_pblock_case2.txt};
\addplot[thick,col1,opacity=1,dashed] table[x=s,y=approx_n40] {accuracy_pblock_case2.txt};

\legend{{$n=24$},{$n=28$},{$n=32$},{$n=36$},{$n=40$}};
\end{axis}
\end{tikzpicture}
\caption{Scaled blocking probability}
\end{subfigure}
%
\caption{Comparison of exact performance measures (solid) against asymptotic approximations (dashed) with $\beta=(s-R_1)/\sqrt{R_1}$ and $\gamma=(n-R_1/r)/\sqrt{R_1/r}$ for $\lambda = 2$, $\mu=1$, $\delta=0.25$ and $p=0.75$.}
\label{fig:accuracy_blocking}
\end{figure}

To test the accuracy of the asymptotic results in Theorem \ref{thm:limits_YT} as approximations in a realistic setting, we plot in Figure \ref{fig:accuracy_blocking} the exact probability of delay and blocking for an Erlang-R model with $R=8$ and $r=0.25$, as a function of $s$. The exact probabilities are given by Equation 
\eqref{eq:delay_probability}, and their respective asymptotic approximations are based on Theorem \ref{thm:limits_YT}. 
Despite the realistic moderate size of the system ($R=8$), we see that the QED approximations are remarkably accurate for many settings $(s,n)$. 
This fast relaxation is in line with observations made earlier in the QED literature \cite{Borst2004,Janssen2011}.

In Appendix \ref{app:accuracy_blocking}, we furthermore compare the asymptotic delay and blocking probability in three additional scenarios. In Tables \ref{tab:numerics_case1}--\ref{tab:numerics_case3} we compute the exact probabilities of delay and blocking through the explicit forms in \eqref{eq:delay_probability} for various values of the offered-load $R_1$, which are omitted here due to space constraints.
The numerical results show that $g^b(\b,\g)$, $f^b(\b,\g)$ and $h^b(\b,\g)$ provide accurate approximations to $\P({\rm delay})$, $\sqrt{R_1}\P({\rm block})$ and $\sqrt{R_1}\,\E[W]$ in pre-limit systems.
The quality of the approximations increases with $R_1$.  
Naturally, fluctuations occur for relatively small values of $R_1$, because $s$ and $n$ need to be rounded to an integer. 

\subsection{QED limits for Erlang-R with holding}
\label{sec:QED_limit_holding}

As explained in \S\ref{sec:QED_scaling}, the model with holding has no product-form steady-state distribution, which makes it hard (if not impossible) to obtain QED limits.
Instead, we derive QED approximations by exploiting a connection with the blocking model. 

We first prove that under scaling \eqref{eq:twofoldscaling}, the upper bound on the utilization level of the nurses needed to achieve stability in the model with holding, as given in Proposition \ref{prop:StabilityCondition}, converges to unity as $R\to\infty$. 
This facilitates high utilization levels of both nurses and beds, a key characteristic of the QED regime.

\begin{proposition}\label{prop:stability_convergence}
Let $s$ and $n$ scale with $R_1\to\infty$ as in \eqref{eq:twofoldscaling}. Then, for $\lambda\to\infty$,
\[
\rho_{\max}(s,n) \to 1.
\]
\end{proposition}

The proof can be found in Appendix \ref{app:proof_stability_convergence}. 
Combining Proposition \ref{prop:stability_convergence} with Proposition \ref{prop:StabilityCondition} shows that indeed the scaling we use results in a highly utilized system.

As observed before, the nature of the two variants of the model is similar up to the fact that a fraction of the patients is deferred on arrive in the setting with blocking, whereas all the arriving patients are eventually admitted into the system in the holding model. 
This implies that, given $s$ and $n$, the nurses face an increased workload in case of a holding room. 
In fact, Theorem \ref{thm:limits_YT} shows that the blocking probability is of order $1/\sqrt{R_1}$, yielding a volume of blocked patients of order $\sqrt{R_1}$ in setting with blocking. 
Accordingly, if $R^b = R_1$ and $R^h$ denote the nominal load arriving to the nurses in the model with blocking and holding, respectively, we can argue that 
\[R^h = R^b + \a \sqrt{R^b} + o(\sqrt{R^b}),\] 
for some $\a>0$. 
Notice that this additional load is of the same order as the safety staffing in the blocking model staffing rule \eqref{eq:twofoldscaling}. 
As $s$ and $n$ remain unchanged, we rewrite \eqref{eq:twofoldscaling} in terms of $R^h$,
\begin{align}
s &= R^h + (\b-\a)\sqrt{R^h} + o(\sqrt{R^h}), \nonumber \\
n &= \frac{R^h}{r} + \left(\g-\a/\sqrt{r}\right)\sqrt{\frac{R^h}{r}} + o(\sqrt{R^h}),
\label{eq:fixed_point_scaling}
\end{align} 
where we have used $R^b = O(R^h)$.
Observe that the square-root principle prevails also after this substitution, albeit with different hedging parameters.
We therefore heuristically argue that the holding model under scaling \eqref{eq:twofoldscaling} with parameters $\b$ and $\g$ mimics the blocking model with parameters $\b-\a$ and $\g-\a/\sqrt{r}$, respectively. 

Observe, however, that we have not yet specified the value of $\a$. 
By definition, $\a\sqrt{R^b}$ is the expected volume of patients that would be rejected in the model with blocking, that is, $R^h$ times the probability of not being admitted to the system directly. 
By the construction in \eqref{eq:fixed_point_scaling}, this volume asymptotically equals $R^h \cdot \P^b({\rm block})$, with parameters $\b-\a$ and $\g-\a/\sqrt{r}$, which by Theorem \ref{thm:limits_YT} is approximated by
\[f^b\left(\b-\a,\gamma-\a/\sqrt{r}\right) / \sqrt{R^h}\]
as $R^h$ grows large.
In conclusion, $\a$ is characterized as the solution of the fixed-point equation
\begin{equation}
 \label{eq:fixedpoint}
 \a = f^h\left(\b-\a,\gamma-\a/\sqrt{r}\right),
 \end{equation} 
 and as a result, we are able to approximate the delay probability in the Erlang-R model with holding as
 \begin{equation}
 \P^h({\rm delay}) \approx g^b(\b-\a,\g-\a/\sqrt{r}) =: g^h(\b,\g).
 \label{eq:fixed_point_Pwait}
 \end{equation}
Likewise, the scaled the mean waiting time for a server can be approximated by
 \begin{equation}
 \sqrt{R_1} \cdot \E[W] \approx h^b(\b-\a,\g-\a/\sqrt{r}) =: h^h(\b,\g).
 \label{eq:fixed_point_Ewait}
 \end{equation}

This also implies that the holding queue is $O(\sqrt{R_1})$.
Subsequently, we argue that the expected holding time (pre-entering wait) under the QED policy is $O(1/\sqrt{R_1})$ and hence asymptotically negligible. 
We justify this claim numerically in \S\ref{sec:analysis}. 

\begin{remark}
\label{rem:holding_limit}
Notice that in the reasoning leading to \eqref{eq:fixedpoint}, we implicitly assumed that the additional volume $\a\sqrt{R_1}$ is an independent Poisson process, which is obviously not the case. Therefore, \eqref{eq:fixed_point_Pwait}--\eqref{eq:fixed_point_Ewait} are approximations for pre-limit systems that are not asymptotically correct as $\lambda\to\iy$.
Nevertheless, our heuristic approach seems to performs well as we confirm numerically next. 
\end{remark} 
 
\begin{figure}[ht]
\centering
\begin{subfigure}{0.48\textwidth}
\centering
\begin{tikzpicture}[scale = 0.9]
\small
\begin{axis}[
	xmin = 0,
	xmax = 16,
	ymin = 0,
	ymax = 1,
	grid = both, 
	axis line style={->},
	axis lines = left,
	xlabel = $s$,
	xscale = 0.9,
	yscale = 0.75,
	legend cell align=left,
	legend style = {at = {(axis cs: 0.1,0.01)},anchor = south west}
]
 
\addplot[ultra thick,col1,opacity=1] table[x=s,y=sim32] {accuracy_pdelay_holding_case2.txt};
\addplot[very thick,col1,opacity=0.7] table[x=s,y=sim36] {accuracy_pdelay_holding_case2.txt};
\addplot[thick,col1,opacity=0.5] table[x=s,y=sim40] {accuracy_pdelay_holding_case2.txt};

\addplot[ultra thick,col1,opacity=1,dashed] table[x=s,y=approx32] {accuracy_pdelay_holding_case2.txt};
\addplot[very thick,col1,opacity=0.7,dashed] table[x=s,y=approx36] {accuracy_pdelay_holding_case2.txt};
\addplot[thick,col1,opacity=0.5,dashed] table[x=s,y=approx40] {accuracy_pdelay_holding_case2.txt};

\legend{{$n=32$},{$n=36$},{$n=40$}};
\end{axis}
\end{tikzpicture}
\caption{Delay probability}
\label{fig:accuracy_holding_a}
\end{subfigure}
\begin{subfigure}{0.48\textwidth}
\centering
\begin{tikzpicture}[scale = 0.9]
\small
\begin{axis}[
	xmin = 0,
	xmax = 16,
	ymin = 0,
	ymax = 1,
	grid = both, 
	axis line style={->},
	axis lines = left,
	xlabel = $s$,
	xscale = 0.9,
	yscale = 0.75,
	legend cell align=left,
	legend style = {at = {(axis cs: 0.1,0.01)},anchor = south west}
]
 
\addplot[very thick,col1,opacity=1] table[x=s,y=sim32] {accuracy_holding_probability.txt};
\addplot[very thick,col1,opacity=0.7] table[x=s,y=sim36] {accuracy_holding_probability.txt};
\addplot[thick,col1,opacity=0.5] table[x=s,y=sim40] {accuracy_holding_probability.txt};

\legend{{$n=32$},{$n=36$},{$n=40$}};
\end{axis}
\end{tikzpicture}
\caption{Holding probability}
\label{fig:accuracy_holding_b}
\end{subfigure}
\caption{Comparison of simulated delay probability (solid) against asymptotic approximations (dashed) with $\beta = (s-R_1)/\sqrt{R_1}$ and $\gamma = (n-R_1/r)/\sqrt{R_1/r}$ for $\lambda = 2$, $\mu=1$, $\delta=0.25$ and $p=0.75$.}
\label{fig:accuracy_holding}
\end{figure}
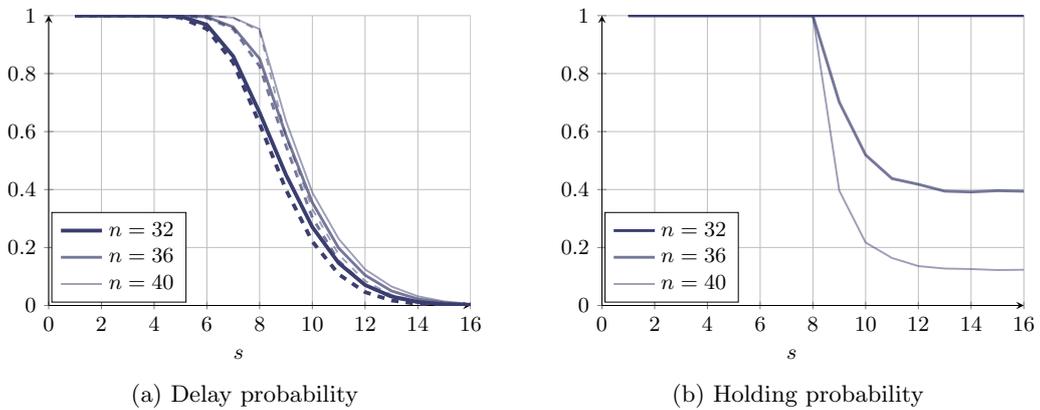

In Figure \ref{fig:accuracy_holding}, we repeat the numerical experiments of Figure \ref{fig:accuracy_blocking} for the model with holding. 
Since the heuristic does not provide an approximation for the holding probability, Figure \ref{fig:accuracy_holding_b} only plots the simulated holding probabilities. 
Those are provided to better understand the implication of operational decision.
Recall that the holding system is only stable (i.e. $\P({\rm hold})<1$) if both $s>R_1=8$ and $n > R_1/r  = 32$.
We nevertheless included the boundary case $n=32$ for illustrative purposes. 
The graphs in Figure \ref{fig:accuracy_holding} show that the heuristic captures the congestion levels well, even for this moderate-size system.
 
To see how this heuristic approach performs under different settings, and particularly if $R_1\to \iy$, we compare in Appendix \ref{app:accuracy_holding} the approximated delay probability in the Erlang-R model with holding as solution of the fixed-point procedure to the outcomes of simulation experiments. 
We omit the tables here due to space constraints.
%
We conclude from these tables that the approximation is good. As $R_1$ increases, the simulated values move closer to the heuristic approximation. Extensive numerical experiments suggest that load is slightly underestimated in the limit. 
The best results in terms of accuracy are attained for small $r$. 
This suggests that the quality of the heuristic method improves as $r$ gets smaller. 
These are exactly the parameter settings for which this model is relevant. 


\section{Dimensioning}
\label{sec:dimensioning}

We will now use the accurate asymptotic approximations of the previous section to define a procedure that determines resource capacity in the restricted Erlang-R models. 
That is, we aim to set the number of nurses $s$ and the number of beds $n$, such that a preset performance level is achieved. 
We take the probability of delay at the needy queue and the probability of blocking/holding at the pre-entrant queue as the target performance objectives.

\subsection{Capacity setting for Erlang-R with blocking}
\label{sec:dimensioning_block}
In the setting with blocking, we can readily use the asymptotic results of Theorem \ref{thm:limits_YT} to (numerically) find a pair of parameters $(\b^*,\g^*)$ to meet the performance requirements.
For instance, given that we want the delay probability to be at most $\varepsilon$, we first solve the equation $g^b(\b^*,\g^*)=\varepsilon$ and then assign $s = \lceil R_1 + \b^*\sqrt{R_1}\rceil$ and $n = \lceil R_1/r+\g^*\sqrt{R_1/r}\rceil$.  Note that there could be multiple solutions to that problem, i.e.\ there could be multiple combinations of number of beds and number of nurses that can result in the same value of a single performance level.
The system manager can ultimately decide which of these feasible solutions fits the environment best, for instance taking into account space and cost constraints.

We illustrate the resource allocation decisions in an MU setting, using data originated from two articles by Lundgren \& Segesten~\cite{LS2001} and Green \& Yankovic~\cite{GY2011}. Green \& Yankovic describe an MU that has 42 beds, with average occupancy level 78\%, and Average Length of Stay (ALOS) 4.3 days. Lundgren \& Segesten studied nurses' service times in a medical-surgical ward. They found that the average service time in their unit was 15.3 minutes per service, and that the average demand rate for each patient is 0.38 requests per hour. Therefore, we take an average service time of 15 minutes and assume that there are 0.4 requests per hour from each patient. Fitting this data to our model results in the following parameters values: $\lambda =0.32, \mu =4, \delta =0.4$, $p=0.975$ and the fraction of needy time is then approximately $r=0.09$. 
This yields nominal offered load $R_1 = 3.2$ and $R_1/r = 34.4$.

\begin{figure}[htbp]
\centering
\begin{subfigure}{0.48\textwidth}
\centering
\begin{tikzpicture}[scale = 0.85]
\begin{axis}[
	xmin = -2,
	xmax = 2,
	ymin = 0,
	ymax = 1,
	grid = both, 
	axis line style={->},
	axis lines = left,
	xlabel = $\beta$,
	legend cell align=left,
	legend style = {at = {(axis cs: -1.9,0.05)},anchor = south west},
	xscale = 0.95,
	yscale = 0.9
]

\addplot[thick,col1,mark repeat = 3, mark = *] table[x=beta,y=delay_gmin1] {staffing_example_with_blocking1.txt};
\addplot[thick,col1,opacity = 0.75,mark repeat = 3, mark = x] table[x=beta,y=delay_g0] {staffing_example_with_blocking1.txt};
\addplot[thick,col1,opacity = 0.5,mark repeat = 3, mark = square] table[x=beta,y=delay_g1] {staffing_example_with_blocking1.txt};
\addplot[thick,col1,opacity = 0.25,mark repeat = 3, mark = triangle*] table[x=beta,y=delay_g2] {staffing_example_with_blocking1.txt};

\draw[->,col1,very thick,dashed] (axis cs: -2,0.5) -- (axis cs: -0.0552366,0.5) -- (axis cs: -0.0552366,0);
\draw[->,col1,opacity = 0.75,very thick,dashed] (axis cs: -2,0.5) -- (axis cs: 0.179728,0.5) -- (axis cs: 0.179728,0);
\draw[->,col1,opacity = 0.5,very thick,dashed] (axis cs: -2,0.5) -- (axis cs: 0.359034,0.5) -- (axis cs: 0.359034,0);
\draw[->,col1,opacity = 0.25,very thick,dashed] (axis cs: -2,0.5) -- (axis cs: 0.459825,0.5) -- (axis cs: 0.459825,0);

\legend{$\gamma = -1$, $\gamma =0$,$\gamma=1$, $\gamma=2$}
\end{axis}

\end{tikzpicture}
\caption{Delay probability}
\label{fig:ratio01_delay}
\end{subfigure}
\begin{subfigure}{0.48\textwidth}
\centering
\begin{tikzpicture}[scale = 0.85]
\begin{axis}[
	xmin = -2,
	xmax = 2,
	ymin = 0,
	ymax = 1,
	grid = both, 
	axis line style={->},
	axis lines = left,
	xlabel = $\beta$,
	legend cell align=left,
	legend style = {at = {(axis cs: 1.9,0.95)},anchor = north east},
	xscale = 0.95,
	yscale = 0.9
]

\addplot[thick,col1,mark repeat = 3, mark = *] table[x=beta,y=delay_gmin1] {staffing_example_with_blocking2.txt};
\addplot[thick,col1,opacity = 0.75,mark repeat = 3, mark = x] table[x=beta,y=delay_g0] {staffing_example_with_blocking2.txt};
\addplot[thick,col1,opacity = 0.5,mark repeat = 3, mark = square] table[x=beta,y=delay_g1] {staffing_example_with_blocking2.txt};
\addplot[thick,col1,opacity = 0.25,mark repeat = 3, mark = triangle*] table[x=beta,y=delay_g2] {staffing_example_with_blocking2.txt};

\draw[ultra thick, col1,dashed,->] (axis cs: -0.0552,0) -- (axis cs: -0.0552,0.292798) -- (axis cs: -2,0.292798);
\draw[ultra thick,col1,opacity = 0.75,dashed,->] (axis cs: 0.179728,0)  -- (axis cs: 0.179728,0.164903)   -- (axis cs: -2,0.164903);
\draw[ultra thick,col1,opacity = 0.5,dashed,->] (axis cs:  0.359034,0)  -- (axis cs:  0.359034,0.0705547)   -- (axis cs: -2,0.0705547);
\draw[ultra thick,col1,opacity = 0.25,dashed,->] (axis cs: 0.459825,0)  -- (axis cs: 0.459825,0.0207909)   -- (axis cs: -2,0.0207909);

\legend{$\gamma = -1$, $\gamma =0$,$\gamma= 1$, $\gamma=2$}
\end{axis}

\end{tikzpicture}
\caption{Blocking probability}
\label{fig:ratio01_block}
\end{subfigure}
\caption{Approximate performance of restricted Erlang-R with blocking for $r \approx 0.09$ and $R_1 = 3.2$, as functions of $\beta$.}
\label{fig:ratio01}
\end{figure}
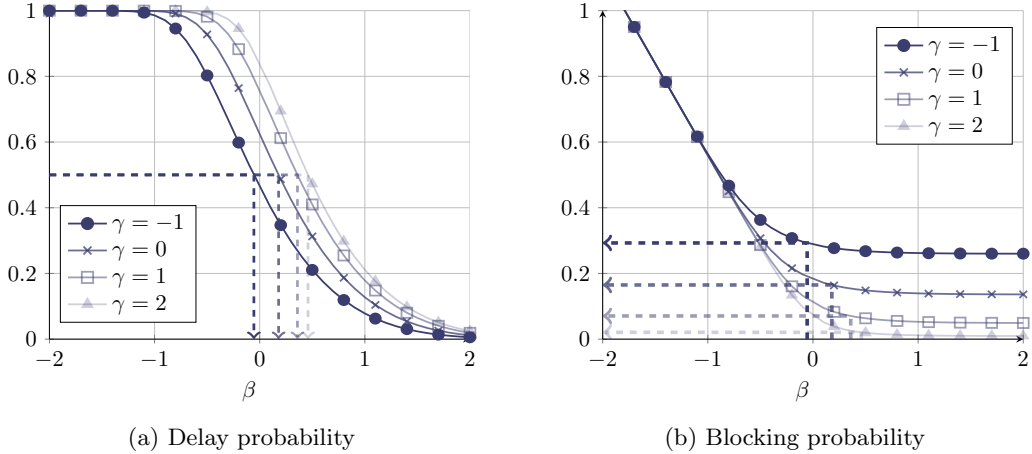

Figure \ref{fig:ratio01} visualizes the dimensioning procedure for this particular MU. 
The hospital management can find a pair of $n$ and $s$ to meet certain criteria, for example to achieve target delay probability $\varepsilon = 0.5$ with reasonable blocking probability.
Figure \ref{fig:ratio01}a indicates that this target can be achieved by a variety of pairs, for instance $(\beta_1,\gamma_1) = (-0.06,-1)$, $(\beta_2,\gamma_2) = (0.16,0)$, $(\beta_3,\gamma_3) = (0.36,1)$ or $(\beta_4,\gamma_4) = (0.46,2)$, among infinitely many others.
According to Figure \ref{fig:ratio01}b, the pairs named above lead to blocking probabilities 0.293, 0.165, 0.071 and 0.021, respectively. 
If the manager decides that probability of blocking of more than 10 percent is not acceptable, this leaves the choices $(\beta_3,\gamma_3) = (0.36,1)$ or $(\beta_4,\gamma_4) = (0.46,2)$ as candidate parameter pairs.
Using the two-fold square-root staffing rule $s_i = \lceil R_1 + \beta_i \sqrt{R_1}\rceil$ and $n_i = [R_1/r + \gamma\sqrt{R_1/r}]$, this yields feasible staffing levels $(s_3,n_3) = (4,40)$ and $(s_4,n_4)=(5,46)$. 
The ultimate decision to apply any of these solutions can be based on external factors, such as operational costs or space limitations of number of beds.


%

\subsection{Capacity setting for Erlang-R with holding}
For the holding model, we need a more sophisticated approach, exploiting the asymptotic approximation with the fixed-point equation in \eqref{eq:fixedpoint}.  We propose the following dimensioning procedure to achieve a preset target delay probability at the needy queue.

\begin{algorithm}\label{alg:stationarydimensioning}
\textsc{Stationary dimensioning algorithm}.\\
\textbf{Input}: Target delay probability $\varepsilon$. Parameters $\l,\m,\d$ and $p$.\\
\textbf{Output}: Staffing levels $s$ and $n$.
\begin{enumerate}[noitemsep]
\item Set $R_1:= \frac{\l}{(1-p)\m}$ and $r = \frac{\d}{\d+p\m}$. 
\item Determine parameters $(\b^*,\g^*)$ such that $g^b(\b^*,\g^*) = \varepsilon$. 
\item Set $\b = \b^* + f^b(\b^*,\g^*)$ and $\g = \g^* + f^b(\b^*,\g^*)/\sqrt{r}$.
\item Return $s = \left\lceil R_1 + \b\sqrt{R_1}\right\rceil$ and $n = \lfloor R_1/r + \g \sqrt{R_1/r}\rfloor$.
\end{enumerate}
\end{algorithm}

\begin{remark}\label{rem:upperboundHW}
In Step 2 of Algorithm \ref{alg:stationarydimensioning} infinitely many pairs $(\b^*,\g^*)$ satisfy the delay probability equation. 
For practical purposes, it is convenient to fix either $\b^*$ or $\g^*$ beforehand, and then solve $g^b(\b^*,\g^*) = \varepsilon$ for the remaining unknown. 
The preset value should however be chosen with care, since $g^b(\b^*,\g^*)$ is upper bounded by the Halfin-Whitt delay probability formula 
\[g_{\rm HW}(\b^*) = \left( 1 + \frac{\b^* \Phi(\b^*)}{\phi(\b^*)}\right)^{-1}.\]
Hence, if $\varepsilon > g_{\rm HW}(\b^*)$, then no feasible solution to $g^b(\b^*,\g^*)=\varepsilon$ exists.
This should be considered when choosing $\beta^*$. 
Furthermore, it is evident from Step 3 that the final values $(\b,\g)$ are always larger than $(\b^*,\g^*)$. 
\end{remark}


We now use the same example as in \S\ref{sec:dimensioning_block} to demonstrate capacity allocation decisions for the model with holding. This can be viewed as the additional capacity the MU needs in terms of nurses and beds, in order to account for the fact that patients are waiting in the ED to be admitted, into the preferred MU, instead of being blocked and transferred to a less preferred unit. 
Observe that the holding model leaves less flexibility for management in choosing system parameters due to stability constraints. For example, the policy with $n=30$ ($\gamma=-0.75$) is infeasible in the holding model.
For similar reasons, only nurse staffing levels with $\beta>0$, or $s > R_1=3.2$ are feasible.

\begin{figure}[htbp]
\centering
\begin{tikzpicture}[scale = 0.9]
\begin{axis}[
	xmin = 0,
	xmax = 2,
	ymin = 0,
	ymax = 1,
	grid = both, 
	axis line style={->},
	axis lines = left,
	xlabel = $\beta$,
	ylabel = {$g^h(\beta,\gamma)$},
	legend cell align=left,
	legend style = {at = {(axis cs: 1.9,0.95)},anchor = north east},
	xscale = 0.95,
	yscale = 0.9
]

\addplot[thick,col1,mark repeat = 3, mark = *] table[x=beta,y=delay_n35] {staffing_example_with_holding1.txt};
\addplot[thick,col1,opacity = 0.75,mark repeat = 3, mark = x] table[x=beta,y=delay_n40] {staffing_example_with_holding1.txt};
\addplot[thick,col1,opacity = 0.5,mark repeat = 3, mark = square] table[x=beta,y=delay_n45] {staffing_example_with_holding1.txt};

\draw[ultra thick, col1,opacity = 0.75,dashed,->] (axis cs: -2,0.5) -- (axis cs: 0.475,0.5) -- (axis cs: 0.475,0);

\legend{$\gamma = -0.75$, $\gamma =0.102$,$\gamma= 0.955$, $\gamma=1.807$};

\end{axis}

\end{tikzpicture}
%
%
%
%
%
\caption{Approximate delay probability of restricted Erlang-R system with holding for $r\approx 0.09$ and $R_1=3.2$ }
\label{fig:ratio01_hold}
\end{figure}
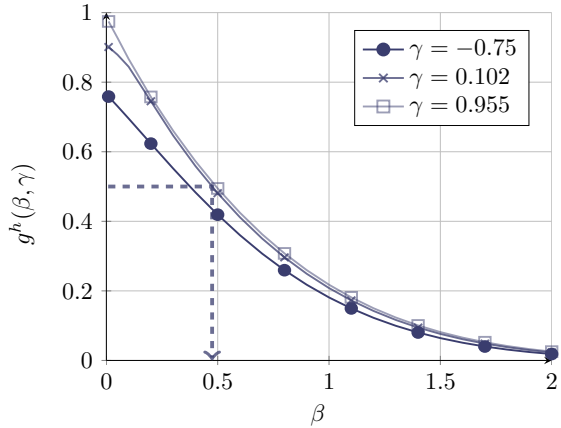

Targeting a delay probability of $0.5$ with $n=40$, Figure \ref{fig:ratio01_hold} shows that operating an MU with holding room requires $\beta = 0.475$ or $s=5$. 
Recall that under the blocking policy, only $s=4$ nurses were needed to achieve a delay probability of $0.5$. 
This example hence shows how the managerial decision to have a holding room, rather than deferring patients to less preferred medical units, requires additional workforce in that unit (as well as the ED). 
This example also shows that the facility with holding room is able to treat fewer patients simultaneously than under blocking constraints, in line with the bounds in \S\ref{sec:bounds} and Conjecture \ref{conj:stochorder}.

\section{Model analysis and managerial implications}
\label{sec:analysis}
In this section, we use the analysis and algorithms developed in earlier sections to gain insights into the importance of the capacity restrictions and customer returns in a restricted Erlang-R system by drawing a comparison to related models studied in the literature. 
\subsection{The influence of customer returns or the role of $r$}
Here we study how the parameter $r$ affects the service level in the restricted Erlang-R model with blocking, on the basis of the asymptotic expressions in Theorem \ref{thm:limits_YT}.


To better understand the connection with the single-station model and the importance of returns we examine the role of $r$. 
Recall the interpretation of $r$ as the fraction of time a patient is needy during his stay within the system in the idealized scenario with infinite capacity, i.e. for $r\in(0,1)$. 
The case $r=1$ corresponds to the setting in which patients are needy all the time, in this case customers get service in one time.
When $r=1$ the infinite-server queue, describing the number of content patients, disappears from the queueing system and we end up with a standard loss model---$M/M/s/n$ queue---in which capacity is scaled as
\[ s = R_1+\beta\sqrt{R_1}, \qquad n = R_1+\gamma\sqrt{R_1}. \]
This staffing rule only makes sense in case $\beta<\gamma$, since no delay is experienced if $n\leq s$. 
If indeed $\gamma>\beta$, then the asymptotic delay probability and scaled blocking probability are given by \cite{masseywallace},
\[
g_B(\beta,\gamma) = \frac{1-{\rm e}^{-\beta(\gamma-\beta)}}{1-{\rm e}^{-\beta(\gamma-\beta)}+\beta\Phi(\beta)/\phi(\beta)},
\qquad f_B(\beta,\gamma) = \frac{\beta{\rm e}^{-\beta(\gamma-\beta)}}{1-{\rm e}^{-\beta(\gamma-\beta)}+\beta\Phi(\beta)/\phi(\beta)}. 
\]

\begin{figure}
\centering
\begin{subfigure}{0.48\textwidth}
\centering
\begin{tikzpicture}[scale=0.9]
\begin{axis}[
	xmin = 0,
	xmax = 1,
	ymin = 0,
	ymax = 0.8,
	grid = both,
	xlabel = $r$,
	axis line style={->},
	legend cell align=left,
	legend style={at={(0.98,1.28)},anchor= north east},
	yscale = 0.75
]

\addplot[thick,col1,mark repeat = 8, mark = *] file {PdelayB_g1_b025.txt};
\addplot[thick,col1,opacity = 0.75,mark repeat = 8, mark = x] file {PdelayB_g1_b05.txt};
\addplot[thick,col1,opacity = 0.5,mark repeat = 8, mark = square] file {PdelayB_g1_b1.txt};
\addplot[thick,col1,opacity = 0.25,mark repeat = 8, mark = triangle*] file {PdelayB_g1_b2.txt};

\legend{$\beta = 0.25$,$\beta=0.5$,$\beta=1.0$,$\beta=2.0$}
\end{axis}
\end{tikzpicture}
\caption{Delay probability}
\end{subfigure}
\begin{subfigure}{0.48\textwidth}
\centering
\begin{tikzpicture}[scale=0.9]
\begin{axis}[
	xmin = 0,
	xmax = 1,
	ymin = 0,
	ymax = 0.4,
	grid = both,
	xlabel = $r$,
	axis line style={->},
	legend cell align=left,
	legend style={at={(0.98,0.05)},anchor= south east},
	yscale = 0.75
]

\addplot[thick,col1,mark repeat = 8, mark = *] file {PblockB_g1_b025.txt};
\addplot[thick,col1,opacity = 0.75,mark repeat = 8, mark = x] file {PblockB_g1_b05.txt};
\addplot[thick,col1,opacity = 0.5,mark repeat = 8, mark = square] file {PblockB_g1_b1.txt};
\addplot[thick,col1,opacity = 0.25,mark repeat = 8, mark = triangle*] file {PblockB_g1_b2.txt};
\addplot[thick,dashed] file {PblockB_g1_inf.txt};
\legend{$\beta = 0.25$,$\beta=0.5$,$\beta=1.0$,$\beta=2.0$}
\end{axis}
\end{tikzpicture}
\caption{Scaled blocking probability}
\label{fig:influence_of_r_b}
\end{subfigure}
\caption{Asymptotic performance measures as a function of $r$ in the restricted Erlang-R model with blocking for $\gamma=1$.}
\label{fig:influence_of_r}
\end{figure}
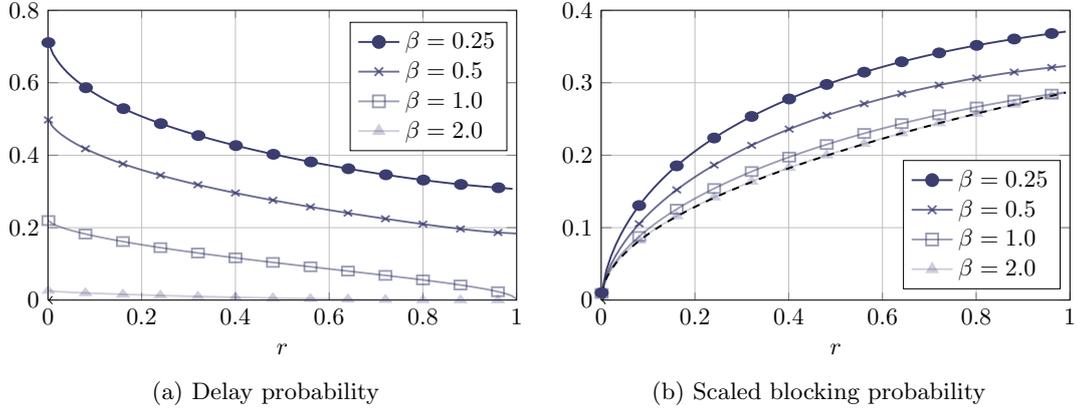

We can see that $f^b(\beta,\gamma)$ of increasing $\b$ approaches a lower bound that is a function of $r$.
To see this, observe that as $\beta$ grows, delays at the nurse queue vanish. 
Then the sojourn time of an admitted patient only consists of a geometric number of needy and content periods with mean $(1/\mu+p/\delta)/(1-p) = 1/r\mu(1-p)$. 
The blocking model can in this case be modeled as an $M/G/n/n$ queue, with offered load $\l/(r\mu(1-p)) =R_1/r$, in which the scaled blocking probability is known to be, see \cite{Janssen2008},
\[\sqrt{R_1} \, \P({\rm block}) = \sqrt{R_1} \, \frac{(R_1/r)^n/n!}{\sum_{k=0}^n (R_1/r)^k / k!} \to \sqrt{r} \, \frac{\phi(\gamma)}{\Phi(\gamma)},\]
as $R_1\to\infty$. 
This function of $r$ is plotted in Figure \ref{fig:influence_of_r_b} as the dashed line.

We observe that in general the probability of blocking increases with $r$, regardless of the capacity constraints on the needy station. 
We can explain this by observing that $r$ influences only $n$ in the QED staffing rule. When $n$ reduces, more patients are blocked. Therefore, if customers spend relatively more time in needy state, which usually indicates services that are less interrupted, blocking will increase. Delays, on the other hand, will decrease in such situations---the minimal delay possible can be achieved if service is given in one time ($r=1$). Returns or interruptions increase delays significantly under QED staffing.


\begin{insight}
\qquad \\
\vspace{-5mm}
\begin{enumerate}
\item Returning customers should be explicitly accounted for in determining staffing in a system with space constraints both in steady-state and transient conditions.
\item The above becomes more important as the proportion of time spend in needy state becomes small, since then the number of customer contributing to the space constrains increases.
\item As returns become more spread over the patient's length-of-stay ($r$ decreases), delay increases and blocking decreases.
\end{enumerate}
\end{insight}

\subsection{Comparing restricted and unrestricted Erlang-R models}

\begin{figure}
\centering
\begin{subfigure}{0.32\textwidth}
\begin{tikzpicture}[scale = 0.66]
\begin{axis}[
	xmin = 0,
	xmax = 3,
	ymin = 0,
	ymax = 1,
	axis line style={->},
	axis lines = left,
	xlabel = $\to \gamma$,
	xscale = 0.9,
	yscale = 0.8,
	legend style = {at = {(0.575,-0.3)},anchor = north},
	legend columns=3,
	legend entries={$\beta=0.1$\ ,$\beta=0.5$\ ,$\beta=1$\ }
]

\addlegendimage{no markers,ultra thick,col1};
\addlegendimage{no markers,ultra thick,col1,opacity=0.5};
\addlegendimage{no markers,ultra thick,col1,opacity=0.25};

\addplot[thick,col1,dashed] file {PdelayR_b01.txt};
\addplot[thick,col1,mark=*, mark repeat = 2] file {PdelayB_r01_b01.txt};
\addplot[thick,col1,mark=square, mark repeat = 2] file {PdelayH_r01_b01.txt};

\addplot[thick,col1,opacity = 0.5,dashed] file {PdelayR_b05.txt};
\addplot[thick,col1,opacity = 0.5,mark=*, mark repeat = 2] file {PdelayB_r01_b05.txt};
\addplot[thick,col1,opacity = 0.5,mark=square, mark repeat = 2] file {PdelayH_r01_b05.txt};

\addplot[thick,col1,opacity = 0.25,dashed] file {PdelayR_b1.txt};
\addplot[thick,col1,opacity = 0.25,mark=*, mark repeat = 2] file {PdelayB_r01_b1.txt};
\addplot[thick,col1,opacity = 0.25,mark=square, mark repeat = 2] file {PdelayH_r01_b1.txt};

\end{axis}
\end{tikzpicture}
\caption{$r=0.1$.}
\end{subfigure}
\begin{subfigure}{0.32\textwidth}
\begin{tikzpicture}[scale = 0.66]
\begin{axis}[
	xmin = 0,
	xmax = 3,
	ymin = 0,
	ymax = 1,
	axis line style={->},
	axis lines = left,
	xlabel = $\to \gamma$,
	xscale = 0.9,
	yscale = 0.8,
	legend style = {at = {(0.575,-0.3)},anchor = north},
	legend columns=3,
	legend entries={$\beta=0.1$\ ,$\beta=0.5$\ ,$\beta=1$\ }
]

\addlegendimage{no markers,ultra thick,col1};
\addlegendimage{no markers,ultra thick,col1,opacity=0.5};
\addlegendimage{no markers,ultra thick,col1,opacity=0.25};
 
\addplot[thick,col1,dashed] file {PdelayR_b01.txt};
\addplot[thick,col1,mark=*, mark repeat = 2] file {PdelayB_r025_b01.txt};
\addplot[thick,col1,mark=square, mark repeat = 2] file {PdelayH_r025_b01.txt};

\addplot[thick,col1,opacity = 0.5,dashed] file {PdelayR_b05.txt};
\addplot[thick,col1,opacity = 0.5,mark=*, mark repeat = 2] file {PdelayB_r025_b05.txt};
\addplot[thick,col1,opacity = 0.5,mark=square, mark repeat = 2] file {PdelayH_r025_b05.txt};

\addplot[thick,col1,opacity = 0.25,dashed] file {PdelayR_b1.txt};
\addplot[thick,col1,opacity = 0.25,mark=*, mark repeat = 2] file {PdelayB_r025_b1.txt};
\addplot[thick,col1,opacity = 0.25,mark=square, mark repeat = 2] file {PdelayH_r025_b1.txt};

\end{axis}
\end{tikzpicture}
\caption{$r=0.25$.}
\end{subfigure}
\begin{subfigure}{0.32\textwidth}
\begin{tikzpicture}[scale = 0.66]
\begin{axis}[
	xmin = 0,
	xmax = 3,
	ymin = 0,
	ymax = 1,
	axis line style={->},
	axis lines = left,
	xlabel = $\to \gamma$,
	xscale = 0.9,
	yscale = 0.8,
	legend style = {at = {(0.575,-0.3)},anchor = north},
	legend columns=3,
	legend entries={$\beta=0.1$\ ,$\beta=0.5$\ ,$\beta=1$\ }
]

\addlegendimage{no markers,ultra thick,col1};
\addlegendimage{no markers,ultra thick,col1,opacity=0.5};
\addlegendimage{no markers,ultra thick,col1,opacity=0.25};

\addplot[thick,col1,dashed] file {PdelayR_b01.txt};
\addplot[thick,col1,mark=*, mark repeat = 2] file {PdelayB_r05_b01.txt};
\addplot[thick,col1,mark=square, mark repeat = 2] file {PdelayH_r05_b01.txt};

\addplot[thick,col1,opacity = 0.5,dashed] file {PdelayR_b05.txt};
\addplot[thick,col1,opacity = 0.5,mark=*, mark repeat = 2] file {PdelayB_r05_b05.txt};
\addplot[thick,col1,opacity = 0.5,mark=square, mark repeat = 2] file {PdelayH_r05_b05.txt};

\addplot[thick,col1,opacity = 0.25,dashed] file {PdelayR_b1.txt};
\addplot[thick,col1,opacity = 0.25,mark=*, mark repeat = 2] file {PdelayB_r05_b1.txt};
\addplot[thick,col1,opacity = 0.25,mark=square, mark repeat = 2] file {PdelayH_r05_b1.txt};

\end{axis}
\end{tikzpicture}
\caption{$r=0.5$.}
\end{subfigure}
\caption{Asymptotic delay probability in open Erlang-R (dashed), restricted Erlang-R with blocking ($\bullet$) and restricted Erlang-R with holding ($\square$), as function of $\gamma$.}
\label{fig:comparison_delay}
\end{figure}
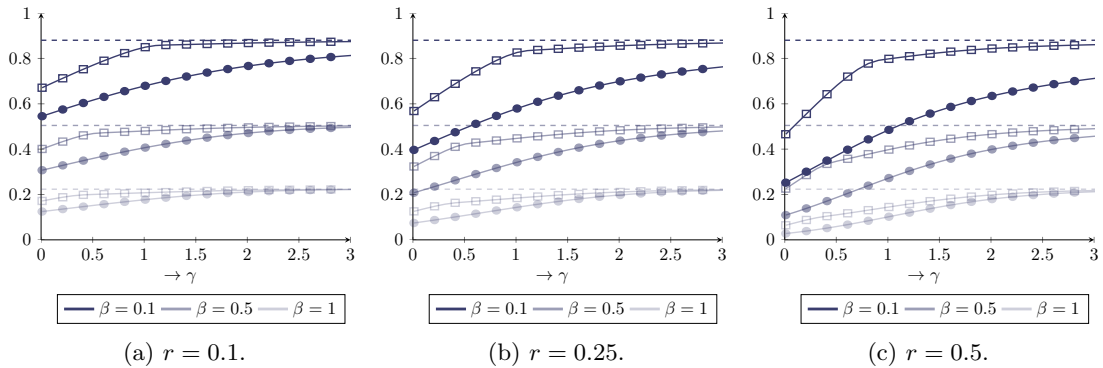

Given the expressions for the asymptotic delay probability in the open Erlang-R model, and its restricted versions with blocking and holding, we compare the three policies for various values of $\beta$, $\gamma$ and $r$. 
Figure \ref{fig:comparison_delay} plots the delay probability for blocking ($g^b(\b,\g)$), holding ($g^h(\b,\g)$) and Erlang-R ($g_{\rm HW}(\beta)$) models,  as functions of $\gamma$, while keeping $\beta$ fixed, for three values of $r$. 
We make a couple of observations.
Notice that
\[ g^b(\beta,\gamma) \leq g^h(\beta,\gamma) \leq g_{\rm HW}(\beta) \]
for all $\beta,\gamma>0$ and $r$. 
In that sense, the holding model is an interpolation between the blocking and the open model. 
As expected, the delay probabilities in the restricted models converge to those of the open Erlang-R model, because increasing $\gamma$ is tantamount to lifting the stringent constraints on the system size. Note that the rate of conversion is fast---one can provide probability of waiting close to that of the open model with small values of $\gamma$. Indeed, the fact that when using QED staffing not much of excessive delay results from the beds restriction is important by itself.
Also, we observe that the difference between delay probabilities increases with $r$.

\subsection{The impact of visit number}
\label{subsec:num_visit}
We next reflect on the impact of operational capacity decisions on different customer populations. We measure patient's complexity by the number of times she needs to see the nurse or the physician during her stay. In the ED context, simple-to-treat patients will need to see the physician once, while complex ones will need multiple visits. Hence, we divide the patients into complexity groups by the number of visits in the needy station. Since the number of visits is geometrically distributed, we have a higher proportion of simple patients than complex ones; that fits well the health care environment.

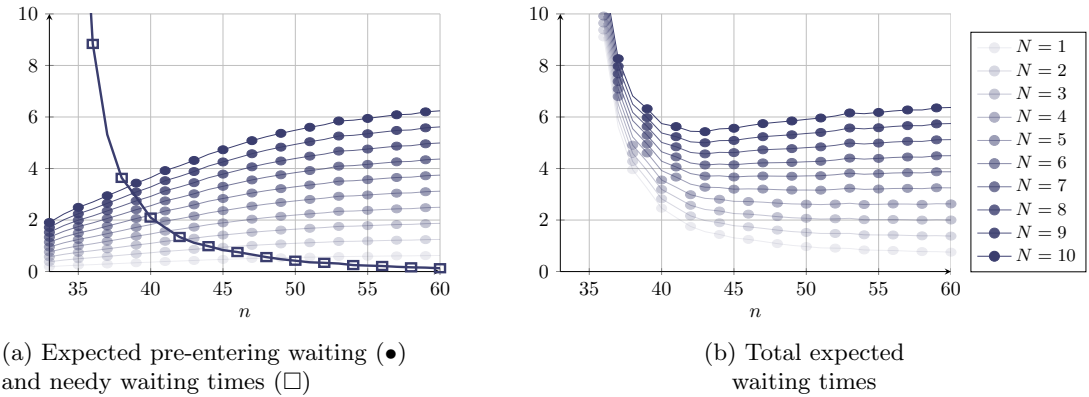
\begin{figure}[h]
\centering
\begin{subfigure}{0.41\textwidth}
\centering
\begin{tikzpicture}[scale=0.75]
\begin{axis}[
	xmin = 33,
	xmax = 60,
	ymin = 0,
	ymax = 10,
	grid = both, 
	axis line style={->},
	axis lines = left,
	xlabel = $n$,
	yscale = 0.8
]
	 	
\addplot[mark=*,col1,opacity=0.1,mark repeat = 2] table[x=n,y=in1] {inner_vs_outer_wait.txt};	
\addplot[mark=*,col1,opacity=0.2,mark repeat = 2] table[x=n,y=in2] {inner_vs_outer_wait.txt};
\addplot[mark=*,col1,opacity=0.3,mark repeat = 2] table[x=n,y=in3] {inner_vs_outer_wait.txt};	
\addplot[mark=*,col1,opacity=0.4,mark repeat = 2] table[x=n,y=in4] {inner_vs_outer_wait.txt};	
\addplot[mark=*,col1,opacity=0.5,mark repeat = 2] table[x=n,y=in5] {inner_vs_outer_wait.txt};	
\addplot[mark=*,col1,opacity=0.6,mark repeat = 2] table[x=n,y=in6] {inner_vs_outer_wait.txt};	
\addplot[mark=*,col1,opacity=0.7,mark repeat = 2] table[x=n,y=in7] {inner_vs_outer_wait.txt};
\addplot[mark=*,col1,opacity=0.8,mark repeat = 2] table[x=n,y=in8] {inner_vs_outer_wait.txt};	
\addplot[mark=*,col1,opacity=0.9,mark repeat = 2] table[x=n,y=in9] {inner_vs_outer_wait.txt};	
\addplot[mark=*,col1,opacity=1,mark repeat = 2] table[x=n,y=in10] {inner_vs_outer_wait.txt};	
 	
\addplot[very thick,mark=square,col1, mark repeat = 2] table[x=n,y=hold] {inner_vs_outer_wait.txt};	

\end{axis}
\end{tikzpicture}
\centering
\caption{Expected pre-entering waiting ($\bullet$) \\and needy waiting times ($\square$)}
\end{subfigure}
\begin{subfigure}{0.58\textwidth}
\centering
\begin{tikzpicture}[scale=0.75]
\begin{axis}[
	xmin = 33,
	xmax = 60,
	ymin = 0,
	ymax = 10,
	grid = both, 
	axis line style={->},
	axis lines = left,
	xlabel = $n$,
	yscale = 0.8,
	legend cell align=left,
	legend style = {at = {(1.05,0.58)}, anchor = west}
]

\addplot[mark=*,col1,opacity=0.1,mark repeat = 2] table[x=n,y=in1] {total_wait.txt};	
\addplot[mark=*,col1,opacity=0.2,mark repeat = 2] table[x=n,y=in2] {total_wait.txt};
\addplot[mark=*,col1,opacity=0.3,mark repeat = 2] table[x=n,y=in3] {total_wait.txt};	
\addplot[mark=*,col1,opacity=0.4,mark repeat = 2] table[x=n,y=in4] {total_wait.txt};	
\addplot[mark=*,col1,opacity=0.5,mark repeat = 2] table[x=n,y=in5] {total_wait.txt};	
\addplot[mark=*,col1,opacity=0.6,mark repeat = 2] table[x=n,y=in6] {total_wait.txt};	
\addplot[mark=*,col1,opacity=0.7,mark repeat = 2] table[x=n,y=in7] {total_wait.txt};
\addplot[mark=*,col1,opacity=0.8,mark repeat = 2] table[x=n,y=in8] {total_wait.txt};	
\addplot[mark=*,col1,opacity=0.9,mark repeat = 2] table[x=n,y=in9] {total_wait.txt};	
\addplot[mark=*,col1,opacity=1,mark repeat = 2] table[x=n,y=in10] {total_wait.txt};	

\legend{{\small $N=1$},{\small $N=2$},{\small $N=3$},{\small $N=4$},{\small $N=5$},{\small $N=6$},{\small $N=7$},{\small $N=8$},{\small $N=9$},{\small $N=10$}};
\end{axis}
\end{tikzpicture}
\caption{Total expected \\ waiting times}
\label{fig:wait_by_visit_b}
\end{subfigure}
\caption{Expected waiting times as a function of $n$ given the number of visits $N$ in the Erlang-R model with holding with $\lambda=2$ $\mu=1$, $\delta=0.25$, $p=0.75$ and $s=9$.}
\label{fig:wait_by_visit}
\end{figure}

Figure \ref{fig:wait_by_visit} shows the waiting time in the needy and pre-entring queues, and the total waiting time, as a function $n$ (number of beds), for each complexity group. 
Obviously, the expected waiting time in the pre-entring queue decreases with $n$, while the needy waiting time increases.
For patients who require a relative large number of visits of the physician, in this case more than 6, the total needy wait is the dominant part of the total waiting time. Therefore, as $n$ grows, the total waiting time first decreases and then increases.
In fact, Figure \ref{fig:wait_by_visit_b} suggests that there is an optimal number of beds $n$ that minimizes the total wait for each complexity type. 
Thus, size restrictions reduce the length-of-stay of patients with complex health conditions (given that the constraint is not too tight).
On the other hand, this figure also shows that no such $n$ exists for patients who only require little assistance.
Hence, there is no $n$ that improves the sojourn time of all patients in the ED simultaneously. 
This leaves the decision to the hospital manager to weigh the importance of patients of different complexity levels.

\begin{remark}
From a different perspective, note that in communication queueing systems, the partition of a job to sizable quantities and scheduling those jobs in a similar dynamic to the Erlang-R model became a popular way for increasing throughput. This is because this effectively schedule jobs by their size even though the total job requirements are uncertain. This in fact creates a shortest-job-first policy without prior knowledge of job size \cite{Comte2016}. Considering that perspective we note that the Erlang-R model actually prioritize simple jobs over complex ones. But without restrictions, when load is too high, such procedures may lead to very long LOS of long jobs. The capacity restriction we analyze in this paper, in both of its versions, limits such delays. Hence, even in cases in which the returns themselves are created by a managerial decision, imposing the additional managerial restriction on entering the system has benefits. 
\end{remark}

\subsection{Case study: comparison of operational decisions}
\label{sec:case_study}
\begin{figure} 
\begin{subfigure}{0.48\textwidth}
\centering
\begin{tikzpicture}[scale=0.75]
\begin{axis}[
	xmin = 0,
	xmax = 24,
	ymin = 0.0,
	ymax = 45,
	ytick = {0,10,20,30,40},	
	grid = both, 
	axis line style={->},
	axis lines = left,
	xlabel = $\to t$,
	xscale=1.25,
	yscale=1,
	legend cell align=left,
	legend style = {at = {(0.01,0.95)},anchor = north west}
]

\addplot[very thick,black] file {lambdaFunction.txt};
\addplot[very thick,col1,mark = *,mark repeat = 2] file {R1.txt};
\addplot[very thick,col1,opacity=0.4,mark = x,mark repeat = 2,mark size = 2] file {R1R2.txt};

\legend{ $\lambda(t)$, $R_1(t)$, $R_1(t)+R_2(t)$};
\end{axis}
\end{tikzpicture}
\caption{Dynamic arrival rate and load functions}
\label{fig:Case_study_arrival_pattern_a}
\end{subfigure}
\begin{subfigure}{0.48\textwidth} 
\centering
\begin{tikzpicture}[scale=0.75]

\begin{axis}[
	xmin = 0,
	xmax = 24,
	ymin = 0.0,
	ymax = 45,
	ytick = {0,10,20,30,40},	
	grid = both, 
	axis line style={->},
	axis lines = left,
	xlabel = $\to t$,
	xscale=1.25,
	yscale=1,
	legend cell align=left,
	legend pos = north west
]

\addplot[very thick,col1,mark = *,mark repeat = 4] file {sFunction.txt};
\addplot[very thick,col1,opacity = 0.4,mark = x,mark repeat = 4,mark size = 2] file {nFunction.txt};

\legend{ $s(t)$, $n(t)$};
\end{axis}
\end{tikzpicture}
\caption{Capacity functions for $\beta=\gamma=0.5$}
\label{fig:Case_study_arrival_pattern_b}
\end{subfigure}
\caption{Empirical arrival rate and offered load functions $R_1(t)$ and $R_1(t)+R_2(t)$ in Israeli ED and corresponding capacity functions}
\label{fig:Case_study_arrival_pattern}
\end{figure}
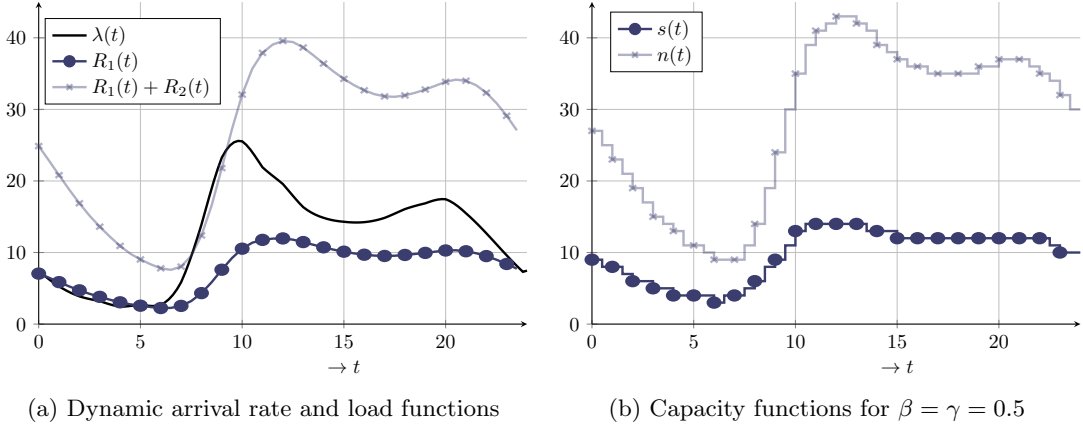

We now illustrate how the managerial decision to operate under a specific operational regime affects ED performance in terms of efficiency and quality-of-care, through a case study. 
The practical environment we investigate is the ED of a moderately-sized hospital, which faces the arrival pattern $\lambda(t)$ plotted in Figure \ref{fig:Case_study_arrival_pattern_a} on a typical workday. 
Other parameters of the model are estimated to be $\mu = 6.67,\ \delta = 2.18$ and $p = 0.76$, so that $r = 0.301$. (Parameters were taken from \cite{YomTov2014}). In order to set time-varying staffing levels $s(t)$ and $n(t)$, we adopt the MOL approximation of the demand process of~\cite{Jennings1996}.
This approach initially presumes infinite capacity to obtain the number of customers $R(t)$ in the queueing system as a function of time. 
This offered load function then replaces (constant) value of $R$ in the stationary dimensioning scheme under consideration, to determine the adequate number of servers at each point in time. 
Following this idea in our two-dimensional queueing system, we find the offered load function for the nurses $R_1(t)$ and the offered load function for the beds $R_1(t)+R_2(t)$ as the solution of the system of ODEs,
\begin{align} \label{eq:offeredloadODE}
\frac{d}{dt} R_1(t) &= \l(t) + \d R_2(t) - \m R_1(t),\\
\frac{d}{dt} R_2(t) &= p\m R_1(t) - \d R_2(t),
\end{align}
see \cite[Thm.~2]{YomTov2014} for details.
For this case study's parameters, these offered load functions are also plotted in Figure \ref{fig:Case_study_arrival_pattern_a}.
While the number of nurses can be adjusted in a relatively flexible manner, the value of $n$, which echoes a hard restriction on the ED capacity, is naturally less amenable to fluctuations. The reason is that the maximum ED capacity is to a large extent determined by its hardware, such as beds and medical equipment.
However, the ED manager might deliberately consider reducing $n$ during more quiet periods of the day, e.g.\ during the night, by imposing bed-to-physician constraints. This is done, for example, when setting a case management constraint \cite{EDexperiment,Campello2016}. 
Therefore, we consider the scenario in which both $s$ and $n$ are time-dependent but we do not force a constant case management quantity, rather let our new methodology to recommend an appropriate one. 

Extrapolating Algorithm \ref{alg:stationarydimensioning} to the time-varying case, Step 4 is replaced by
\begin{align*}
s(t) &= R_1(t) + \b\sqrt{R_1(t)},\\
n(t) &= R_1(t)+R_2(t) + \g\sqrt{R_1(t)+R_2(t)},
\end{align*}
for some $\b,\g>0$. 
Since $R_1(t)$ and $R_2(t)$ are given, the QED staffing problem again reduces to finding the pair $(\beta,\gamma)$.

Figure \ref{fig:Case_study_arrival_pattern_b} plots the capacity functions for $\beta = 0.5$ and $\gamma=0.5$, assuming capacity can only be adjusted every 30 minutes. 
In this case study, we consider three pairs of parameters $(\beta,\gamma)$.
For each of these we investigate, using simulation, the differences in the time-varying performance indicators between the policy with blocking and holding.

\begin{figure}[h]
\centering

\begin{subfigure}{0.33\textwidth}
\centering
\begin{tikzpicture}[scale=0.54]

\begin{axis}[
	xmin = 0,
	xmax = 24,
	ymin = 0.0,
	ymax = 1.0,
	ytick = {0,0.2,0.4,0.6,0.8,1.0},	
	xtick = {0,3,6,9,12,15,18,21,24},
	grid = both, 
	axis line style={->},
	axis lines = left,
	xlabel = $\to t$,
	xscale=1.1,
	yscale=1
]

\addplot[very thick,col1] table[x=t,y=delay_b01g2] {case_study_block.txt};
\addplot[very thick,col1,dashed] table[x=t,y=delay_b01g2] {case_study_hold.txt};

\addplot[very thick,col5] table[x=t,y=delay_b1g15] {case_study_block.txt};
\addplot[very thick,col5,dashed] table[x=t,y=delay_b1g15] {case_study_hold.txt};

\addplot[very thick,col3] table[x=t,y=delay_b2g1] {case_study_block.txt};
\addplot[very thick,col3,dashed] table[x=t,y=delay_b2g1] {case_study_hold.txt};
\end{axis}
\end{tikzpicture}
\caption{$\P({\rm delay})$}
\label{fig:simulation_results_a}
\end{subfigure}
\begin{subfigure}{0.33\textwidth}
\centering
\begin{tikzpicture}[scale=0.54]

\begin{axis}[
	xmin = 0,
	xmax = 24,
	ymin = 0.0,
	ymax = 0.52,
	ytick = {0,0.1,0.2,0.3,0.4,0.5},
	xtick = {0,3,6,9,12,15,18,21,24},	
	grid = both, 
	axis line style={->},
	axis lines = left,
	xlabel = $\to t$,
	xscale=1.1,
	yscale=1,
	legend cell align=left,
	legend style = {at = {(0.9,0.95)}, anchor = north east}
]
\addplot[very thick,col1] table[x=t,y=block_b01g2] {case_study_block.txt};
\addplot[very thick,col5] table[x=t,y=block_b1g15] {case_study_block.txt};
\addplot[very thick,col3] table[x=t,y=block_b2g1] {case_study_block.txt};

\addplot[very thick,col1,dashed] table[x=t,y=block_b01g2] {case_study_hold.txt};
\addplot[very thick,col5,dashed] table[x=t,y=block_b1g15] {case_study_hold.txt};
\addplot[very thick,col3,dashed] table[x=t,y=block_b2g1] {case_study_hold.txt};

\legend{{$(\beta,\gamma)=(0.1,2)$},{$(\beta,\gamma)=(1,1.5)$},{$(\beta,\gamma)=(2,1)$}};
\end{axis}
\end{tikzpicture}
\caption{$\P({\rm block})$ or $\P({\rm hold})$}
\label{fig:simulation_results_b}
\end{subfigure}
\begin{subfigure}{0.32\textwidth}
\centering
\begin{tikzpicture}[scale=0.54]

\begin{axis}[
	xmin = 0,
	xmax = 24,
	ymin = 0.0,
	ymax = 5.5,
	xtick = {0,3,6,9,12,15,18,21,24},	
	grid = both, 
	axis line style={->},
	axis lines = left,
	xlabel = $\to t$,
	xscale=1.1,
	yscale=1
]

\addplot[very thick,col1] table[x=t,y=ratio_b01g2] {case_study_block.txt};
\addplot[very thick,col1,dashed] table[x=t,y=ratio_b01g2] {case_study_hold.txt};

\addplot[very thick,col5] table[x=t,y=ratio_b1g15] {case_study_block.txt};
\addplot[very thick,col5,dashed] table[x=t,y=ratio_b1g15] {case_study_hold.txt};

\addplot[very thick,col3] table[x=t,y=ratio_b2g1] {case_study_block.txt};
\addplot[very thick,col3,dashed] table[x=t,y=ratio_b2g1] {case_study_hold.txt};
\end{axis}
\end{tikzpicture}
\caption{Nurse-to-patient ratio.}
\label{fig:simulation_results_c}
\end{subfigure}

\caption{Simulation results for case study. Solid and dashed lines represent time-varying performance in the blocking and holding model, respectively.}
\label{fig:simulation_results}
\end{figure}
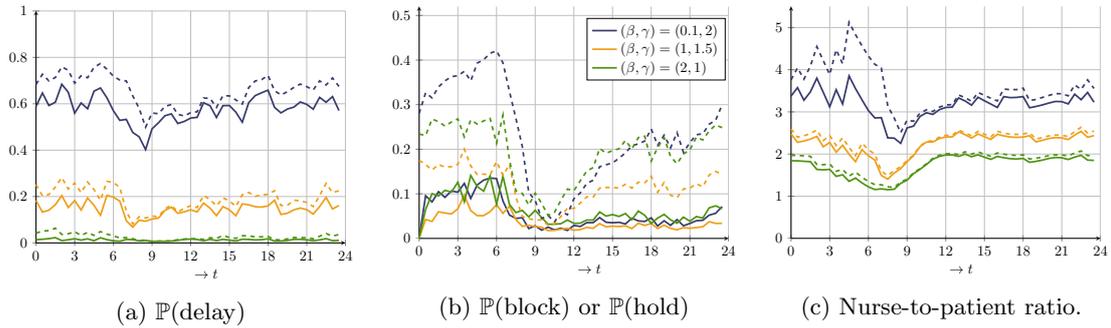

The simulation results are presented in Figure \ref{fig:simulation_results}. 
Figure \ref{fig:simulation_results_a} shows that the MOL approach for capacity allocation roughly stabilizes the delay probability. 
Figure \ref{fig:simulation_results_b} shows that the fraction of patients not entering the ED on arrival in the blocking model is reasonable for all parameter pairs considered and are ordered according to $\gamma$.
We also see a significant difference with holding.
Observe also that the holding probability drops in the period 8--13, which is exactly the period when the system is experiencing peak offered load. 
Hence, this temporary reduction is in line with our asymptotic findings that the probability of blocking/holding is $O(1/\sqrt{R_1})$.

Finally note that the three parameter settings lead to different nurse-to-patient ratios. 
Clearly, larger $\beta$ leads to small nurse-to-patient ratios (due do larger staffing). 
Figure \ref{fig:simulation_results_c} demonstrates that for $(\beta,\gamma) = (1,1.5)$ and $(\beta,\gamma) = (2,1)$ the difference between the holding policy and the blocking policy is small. However, for $(\beta,\gamma) = (0.1,2)$ we see a significant increase in the ratio during night hours. 
This may be due to the tight nurse schedule, that causes the holding queue to build up just before midnight. 
This queue then empties latter on, causing an increase in workload per nurse in the period 24--7.

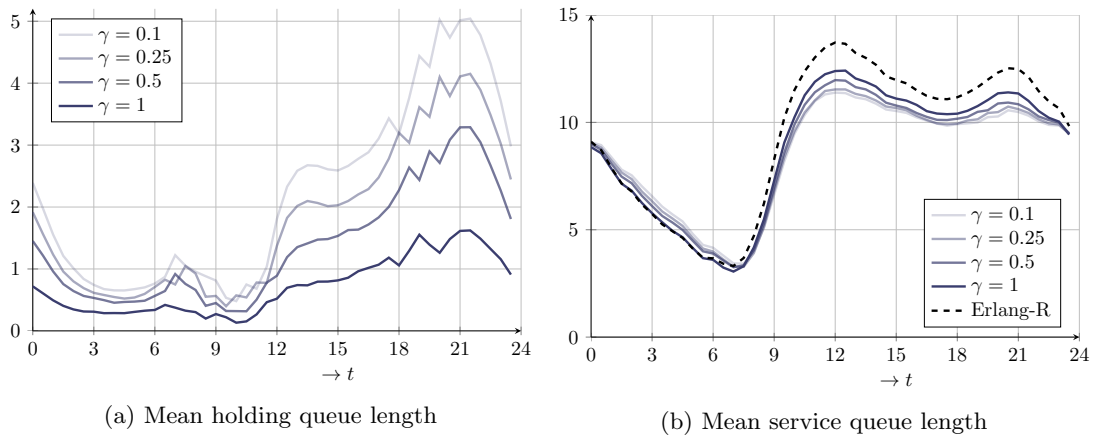
\begin{figure}
\centering
\begin{subfigure}{0.48\textwidth}
\centering
\begin{tikzpicture}[scale=0.75]
\begin{axis}[
	xmin = 0,
	xmax = 24,
	ymin = 0.0,
	ymax = 5.2,
	xtick = {0,3,6,9,12,15,18,21,24},	
	grid = both, 
	axis line style={->},
	axis lines = left,
	xlabel = $\to t$,
	xscale=1.25,
	yscale=1,
	legend cell align=left,
	legend pos = north west
]

\addplot[very thick,col1,opacity=0.2] file {holdingQueue_g01.txt};
\addplot[very thick,col1,opacity=0.45] file {holdingQueue_g025.txt};
\addplot[very thick,col1,opacity=0.7] file {holdingQueue_g05.txt};
\addplot[very thick,col1] file {holdingQueue_g1.txt};

\legend{$\gamma=0.1$, $\gamma=0.25$,$\gamma=0.5$,$\gamma=1$}
\end{axis}
\end{tikzpicture}
\caption{Mean holding queue length}
\end{subfigure}
\begin{subfigure}{0.48\textwidth}
\centering
\begin{tikzpicture}[scale=0.75]
\begin{axis}[
	xmin = 0,
	xmax = 24,
	ymin = 0.0,
	ymax = 15,
	xtick = {0,3,6,9,12,15,18,21,24},	
	grid = both, 
	axis line style={->},
	axis lines = left,
	xlabel = $\to t$,
	xscale=1.25,
	yscale=1,
	legend cell align=left,
	legend style = {at = {(0.77,0.03)},anchor = south east}
]

\addplot[very thick,col1,opacity=0.2] file {serviceQueue_g01.txt};
\addplot[very thick,col1,opacity=0.45] file {serviceQueue_g025.txt};
\addplot[very thick,col1,opacity=0.7] file {serviceQueue_g05.txt};
\addplot[very thick,col1] file {serviceQueue_g1.txt};
\addplot[very thick,dashed] file {serviceQueue_R.txt};
\legend{$\gamma=0.1$, $\gamma=0.25$,$\gamma=0.5$,$\gamma=1$,Erlang-R}
\end{axis}
\end{tikzpicture}
\caption{Mean service queue length}
\end{subfigure}
\caption{Simulated queue length of holding model with different values of $\gamma$.}
\label{fig:simulation_queuelengths}
\end{figure}

To see the direct effect of the size restriction on the queue lengths, we plotted the mean holding and service queue lengths in the holding model as a function of the parameter $\gamma$ in Figure \ref{fig:simulation_queuelengths}.
Note that for all $\gamma$ considered, the holding queue length are, as expected, of a smaller order than the number of patients admitted. 
Also, the holding queue length decreases as we increase $\gamma$. 
The service queue lengths naturally approach the expected queue lengths in the Erlang-R model as $\gamma$ is increased.

\section{Conclusions and future research}
In this research we developed and analyzed a queueing network with repeated services that has capacity restrictions both on service and accessibility. We not only developed approximations for the performance of the system with a blocking policy, but also constructed a new innovative way to use those results to approximate performance of the intractable model with holding. Those approximations were based on a fixed-point analysis, and enabled us to take the first step towards characterizing the pre-entering queue behavior in the QED regime. 
In addition, the unified analysis enabled us to connect the two policies and to deepen the understanding of the practical application of the two. We showed in two realistic case studies the application of those two models to medical unit resource allocation. We discussed the implication of holding patients both in terms of the range of possible combinations that results from stability constraints, and in terms of its impact on waiting in the system. Finally, we showed that the approximations are accurate for a wide range of parameters.


\bibliographystyle{plain}
\bibliography{references_rer}

\begin{thebibliography}{10}

\bibitem{AHA2007}
{American Hospital Association}.
\newblock Survey of hospital leaders, 2007.

\bibitem{Armony2015}
M.~Armony, S.~Israelit, A.~Mandelbaum, Y.N. Marmor, Y.~Tseytlin, and G.B.
  Yom-Tov.
\newblock On patient flow in hospitals: A data-based queueing-science
  perspective.
\newblock {\em Stochastic Systems}, 5(1):146--194, 2015.

\bibitem{Bekker2009a}
R.~{Bekker} and A.M. de~{Bruin}.
\newblock Time-dependent analysis for refused admissions in clinical wards.
\newblock {\em Annals of Operations Research}, 178(1):45--65, 2009.

\bibitem{israelit}
R.~Bennidor and S.H. Israelit.
\newblock Emergency department intermediate stay unit - a failed model.
\newblock Technion working paper, 2015.

\bibitem{Comte2016}
T.~Bonald and C.~Comt\'e.
\newblock Networks of multi-server queues with parallel processing.
\newblock Unpublished manuscript, 2016.

\bibitem{Borst2004}
S.C. Borst, A.~Mandelbaum, and M.I. Reiman.
\newblock Dimensioning large call centers.
\newblock {\em Operations Research}, 52(1):17--34, January-February 2004.

\bibitem{Campello2016}
F.~Campello, A.~Ingolfsson, and R.A. Shumsky.
\newblock Queueing models of case managers.
\newblock {\em Management Science}, 2016.
\newblock Articles in advance.

\bibitem{Carmen2016}
R.~Carmen and I.~van Nieuwenhuyse.
\newblock How inpatient boarding impacts {ED} performance: A queueing analysis.
\newblock Working paper, KU Leuven, 2016.

\bibitem{ChanYao}
H.~Chen and D.D. Yao.
\newblock {\em Fundamentals of Queueing Networks: Performance, Asymptotics, and
  Optimization}.
\newblock Springer, 2001.

\bibitem{Bekker2009b}
A.M. de~{Bruin}, R.~{Bekker}, L.~van {Zanten}, and G.M. {Koole}.
\newblock Dimensioning hospital wards using the {E}rlang loss model.
\newblock {\em Annals of Operations Research}, 178(1):23--43, 2009.

\bibitem{Denton2013}
B.T. {Denton}, editor.
\newblock {\em Handbook of healthcare operations management: Methods and
  applications}.
\newblock North Holland, New York, 2nd edition, 2013.

\bibitem{diwas}
KC~S. Diwas.
\newblock Does multitasking improve performance? evidence from the emergency
  department.
\newblock {\em Manufacturing \& Service Operations Management}, 16(2):168--183,
  2014.

\bibitem{Gamarnik2013}
D.~Gamarnik and D.A. Goldberg.
\newblock On the rate of convergence to stationarity of the ${M}/{M}/{N}$ queue
  in the {H}alfin-{W}hitt regime.
\newblock {\em The Annals of Applied Probability}, 23(5):1879--1912, 2013.

\bibitem{GY2011}
L.~Green and N.~Yankovic.
\newblock Identifying good nursing levels: A queuing approach.
\newblock {\em Operations Research}, 59(4):942--955, 2011.

\bibitem{Green2006}
L.V. {Green}.
\newblock Using queueing theory to increase the effectiveness of physician
  staffing in the emergency department.
\newblock {\em Academic Emergency Medicine}, 13:61--68, 2006.

\bibitem{Green2008}
L.V. Green and S.~Savin.
\newblock Reducing delays for medical appointments: A queueing approach.
\newblock {\em Operations Research}, 56(6):1526--1538, 2008.

\bibitem{HalfinWhitt1981}
S.~{Halfin} and W.~{Whitt}.
\newblock Heavy-traffic limits for queues with many exponential servers.
\newblock {\em Operations Research}, 29(3):567--588, 1981.

\bibitem{Hall2006}
R.W. {Hall}, editor.
\newblock {\em Patient Flow: Reducing Delay in Healthcare Delivery}.
\newblock Springer, 2006.

\bibitem{Hall2012}
R.W. {Hall}, editor.
\newblock {\em Handbook of Healthcare System Scheduling}.
\newblock Springer, 2012.

\bibitem{Junfei2015}
Junfei Huang, Boaz Carmeli, and Avishai Mandelbaum.
\newblock Control of patient flow in emergency departments, or multiclass
  queues with deadlines and feedback.
\newblock {\em Operations Research}, 63(4):892--908, 2015.

\bibitem{Jackson1963}
J.R. Jackson.
\newblock Jobshop-like queueing systems.
\newblock {\em Management Science}, 10(1):131--142, 1963.

\bibitem{Janssen2008}
A.J.E.M. {Janssen}, J.S.H. van {Leeuwaarden}, and B.~{Zwart}.
\newblock Gaussian expansions and bounds for the {P}oisson distribution applied
  to the {E}rlang-{B} formula.
\newblock {\em Advances in Applied Probability}, 40(1):122--143, 2008.

\bibitem{Janssen2011}
A.J.E.M. {Janssen}, J.S.H. van {Leeuwaarden}, and B.~{Zwart}.
\newblock Refining square-root safety staffing by expanding {E}rlang {C}.
\newblock {\em Operations Research}, 59(6):1512--1522, 2011.

\bibitem{Jennings2008}
O.B. {Jennings} and F.~de~{V\'ericourt}.
\newblock Dimensioning large-scale memberschip services.
\newblock {\em Operations Research}, 55(1):173--187, 2008.

\bibitem{Jennings2011}
O.B. {Jennings} and F.~de~{V\'ericourt}.
\newblock Nurse staffing in medical units: A queueing perspective.
\newblock {\em Operations Research}, 59(6):1320--1331, 2011.

\bibitem{Jennings1996}
O.B. Jennings, A.~Mandelbaum, W.A. Massey, and W.~Whitt.
\newblock Server staffing to meet time-varying demand.
\newblock {\em Management Science}, 42(10):1383--1394, 1996.

\bibitem{Kazahaya2005}
G.~Kazahaya.
\newblock Harnessing technology to redesign labor cost management reports.
\newblock {\em Healthcare Financial Management}, 59(4):94--100, 2005.

\bibitem{Khudyakov2010}
P.~Khudyakov, P.D. Feigin, and A.~Mandelbaum.
\newblock Designing a call center with an {I}{V}{R} ({I}nteractive {V}oice
  {R}esponse).
\newblock {\em Queueing Systems}, 66(3):215--237, 2010.

\bibitem{LS2001}
S.~Lundgren and K.~Segesten.
\newblock Nurses use of time in a medical-surgical ward with all-rn staffing.
\newblock {\em Journal of Nursing Management}, 9(1):13--20, 2001.

\bibitem{masseywallace}
W.A. Massey and R.B. Wallace.
\newblock An asymptotically optimal design of the {${M}/{M}/c/k$} queue.
\newblock Unpublished report, 2004.

\bibitem{Neuts1981}
M.F. {Neuts}.
\newblock {\em Matrix-Geometric Solutions in Stochastic Models}.
\newblock The John Hopkins University Press, Baltimore, 1981.

\bibitem{Sanders2016}
J.~Sanders, S.C. Borst, A.J.E.M. Janssen, and J.S.H. vvan Leeuwaarden.
\newblock Optimality gaps in asymptotic dimensioning of many-server systems.
\newblock {\em Operations Research Letters}, 44(3):369--365, 2016.

\bibitem{EDexperiment}
H.~Song, A.L. Tucker, and K.L. Murrell.
\newblock The diseconomies of queue pooling: An empirical investigation of
  emergency department length of stay.
\newblock {\em Management Science}, 61(12):3032--3053, 2015.

\bibitem{Leeuwaarden2011}
J.S.H. van Leeuwaarden and C.~Knessl.
\newblock Transient behavior of the halfin–whitt diffusion.
\newblock {\em Stochastic Processes and their Applications}, pages 1524--1545,
  2011.

\bibitem{Leeuwaarden2012}
J.S.H. van Leeuwaarden and C.~Knessl.
\newblock Spectral gap of the {E}rlang-{A} model in the {H}alfin-{W}hitt
  regime.
\newblock {\em Stochastic Systems}, 2(1):149--207, 2012.

\bibitem{Leeuwaarden2015}
J.S.H. van Leeuwaarden, B.W.J. Mathijsen, and F.~Sloothaak.
\newblock Delayed workload shifting in many-server systems.
\newblock {\em ACM SIGMETRICS Performance Evaluation Review}, 43(2):10--12,
  2015.

\bibitem{Leeuwaarden2016}
J.S.H. van Leeuwaarden, B.W.J. Mathijsen, and F.~Sloothaak.
\newblock Cloud provisioning in the {Q}{E}{D} regime.
\newblock {\em Proceedings of the 9th EAI International Conference on
  Performance Evaluation Methodologies and Tools}, pages 180--187, 2016.

\bibitem{YomTov2014}
G.B {Yom-Tov} and A.~{Mandelbaum}.
\newblock Erlang-{R}: A time-varying queue with reentrant customers, in support
  of healthcare staffing.
\newblock {\em Manufacturing \& Service Operations Management}, 16(2):283--299,
  2014.

\bibitem{erlanga}
B.~Zhang, J.S.H. van Leeuwaarden, and B.~Zwart.
\newblock Refining square-root staffing for call centers with impatient
  customers.
\newblock {\em Operations Research}, 60:461--474, 2012.

\bibitem{zychlinski2016bedblocking}
N.~Zychlinski, A.~Mandelbaum, P.~Mom{\v{c}}ilovi{\'c}, and I.~Cohen.
\newblock Bed blocking in hospitals due to scarce capacity in geriatric
  institutions -- cost minimization via fluid models.
\newblock Under review, 2016.

\end{thebibliography}

\appendix

\section{Description of the QBD process}
\label{app:QBDdescription}

\subsection{The QBD-process}
\label{app:theQBDprocess}
We consider the QBD-process $X=\{N,Q_1\}$ in stationarity. Let $\nu(i)=\min\{i,s\}\mu$. To determine the (outgoing) transition rates of the process $X$ we distinguish between the following cases:
\begin{itemize}
\item \emph{Transitions from $(0,0)$:} There are no patients in the Emergency Department and thus the only possible occurrence is when a new patient arrives. This results in a transition to $(1,1)$ and occurs with rate $\lambda$.
\item \emph{Transitions from $(i,0), 1 \leq i < n$:} There are exactly $i$ patients assigned to a bed of which none are seen by a nurse. Then either one of those patients becomes needy, or a new patient arrives at the Emergency Department that can immediately be seen by a nurse. The first results in a transition to $(i,1)$ and occurs at rate $i \delta$, and the second results in a transition to $(i+1,1)$ and occurs with rate $\lambda$.
\item \emph{Transitions from $(i,0), i \geq n$:} Again, the only possible transitions arises from either a newly arrived patient or a patient assigned to a bed becoming needy. However, a newly arrived patient finds all beds occupied and needs to wait. Thus, with rate $\lambda$ we have a transition to $(i+1,0)$ and with rate $n \delta$ a transition to $(i,1)$.
\item \emph{Transitions from $(i,i), i < n$:} In this case all patients assigned to a bed are in need of service. With rate $\lambda$ a new patient arrives at the Emergency Department. She joins the (possible) queue to be seen by a nurse immediately, so this results in a transition to $(i+1,i+1)$. Moreover, since there are only $s < n$ nurses, a service completion occurs with rate $\nu(i)$. With probability $p$ the patient turns to the holding phase, so in total we still have $i$ patients with one patient less in queue for a nurse. With probability $1-p$ the patient leaves the Emergency Department, decreasing both $N$ and $Q_1$ by one. In other words, with rate $p \nu(i)$ we have a transition to $(i,i-1)$ and with rate $(1-p)\nu(i)$ we have a transition to $(i-1,i-1)$.
\item \emph{Transitions from $(n,n)$:} Similar to the previous case, we have a transition to $(n,n-1)$ with rate $p s \mu$ and with rate $(1-p)s \mu$ we have a transition to $(n-1,n-1)$. In this case however, a newly arrived patient finds all beds occupied, resulting in a transition to $(n+1,n)$ with rate $\lambda$.
\item \emph{Transitions from $(i,n), i > n$:} We have a transition to $(i+1,n)$ with rate $\lambda$ and a transition to $(i,n-1)$ with rate $p s \mu$. In case that a patient leaves the Emergency Department there are $i-n>0$ patients in the holding room waiting for an available bed. Thus, one of the $i-n$ patients in the holding room is assigned to the available bed in need of service. That is, with rate $(1-p) s \mu$ we have a transition to $(i-1,n)$.
\item \emph{Transitions from $(i,j), 1 \leq j < i  < n$:} There are four possible transitions. First, with rate $\lambda$ there is a new arrival which results in a transition to $(i+1,j+1)$. Second, with rate $(i-j) \delta$ a patient in one of the beds becomes needy, which results in a transition to $(i,j+1)$. Third, with rate $p \nu(j)$ a patient turns to the content state after service completion, which results in a transition to $(i,j-1)$. Last, with rate $(1-p) \nu(j)$ a patient leaves the Emergency Department after service completion, which results in a transition to $(i-1,j-1)$.
\item \emph{Transitions from $(n,j), 1 \leq j  < n$:} This case is similar to the previous one. The only difference arises when a new patient arrives, since all $n$ beds are already occupied. Thus, with rate $\lambda$ we have a transition to $(n+1,j)$.
\item \emph{Transitions from $(i,j), i > n, 1 \leq j  \leq n$:} This case is the previous one, except when a patient leaves the Emergency Department after service completion. Then one of the $(i-n)$ patients in the holding room will be assigned to a bed in need of service. This results in a transition to $(i-1,j)$ with rate $(1-p) \nu(j)$.
\end{itemize}

\noindent
The state space and transition rates of the Erlang-R model with holding are illustrated in Figure~\ref{fig:QBDIllustration}.

\begin{figure}
 \centering
\begin{tikzpicture}[scale=0.8]
\draw[step=1cm,gray!50!,very thin] (0,0) grid (15.5,8.5);
\draw[thick,->] (0,0) -- (15.5,0);
\draw[thick,->] (0,0) -- (0,8.5);
\draw[thick] (0,0) -- (8,8);
\draw[thick,dashed,black!50!] (8,0) -- (8,8);
\draw[thick] (8,8) -- (15.5,8);
\foreach \x in {0,1,2,3,4,5,6,7,8,9,10,11,12,13,14,15}
	\foreach \y in {0,1,2,3,4,5,6,7,8}
    \draw[fill] (\x,\y) circle [radius=0.025];

\node [below left] at (0,0) {$0$};
\node [left] at (0,2) {$j$};
\node [left] at (0,5) {$k$};
\node [left] at (0,8) {$n$};
\node [below] at (6,0) {$i$};
\node [below] at (8,0) {$n$};
\node [above left] at (0,8.5) {$Q_1$}; 
\node [below right] at (15.5,0) {$N$}; 

\path [->,thick,-latex] (0,0) edge [bend right] (1,1);

\path [->,thick,-latex] (3,0) edge (3,1);
\path [->,thick,-latex] (3,0) edge (4,1);

\path [->,thick,-latex] (4,4) edge [bend right] (5,5);
\path [->,thick,-latex] (4,4) edge (4,3);
\path [->,thick,-latex] (4,4) edge [bend right] (3,3);

\path [->,thick,-latex] (6,2) edge (6,3);
\path [->,thick,-latex] (6,2) edge (7,3);
\path [->,thick,-latex] (6,2) edge (6,1);
\path [->,thick,-latex] (6,2) edge (5,1);

\path [->,thick,-latex] (8,5) edge (8,6);
\path [->,thick,-latex] (8,5) edge (8,4);
\path [->,thick,-latex] (8,5) edge (9,5);
\path [->,thick,-latex] (8,5) edge (7,4);

\path [->,thick,-latex] (8,8) edge (8,7);
\path [->,thick,-latex] (8,8) edge [bend right] (7,7);
\path [->,thick,-latex] (8,8) edge [bend left] (9,8);

\path [->,thick,-latex] (11,8) edge (11,7);
\path [->,thick,-latex] (11,8) edge [bend left] (12,8);
\path [->,thick,-latex] (11,8) edge [bend left] (10,8);

\path [->,thick,-latex] (11,0) edge (11,1);
\path [->,thick,-latex] (11,0) edge [bend left] (12,0);

\path [->,thick,-latex] (12,5) edge (12,6);
\path [->,thick,-latex] (12,5) edge (12,4);
\path [->,thick,-latex] (12,5) edge (13,5);
\path [->,thick,-latex] (12,5) edge (11,5);

\node [above] at (12.75,5) {\scriptsize $\lambda$};
\node [above] at (11.25,5) {\scriptsize $(1-p)\nu(k)$};
\node [right] at (12,5.75) {\scriptsize $(n-k)\delta$};
\node [right] at (12,4.25) {\scriptsize $p \nu(k)$};

\node [above] at (6.75,2.25) {\scriptsize $\lambda$};
\node [below] at (5,2) {\scriptsize $(1-p)\nu(j)$};
\node [above] at (6,3) {\scriptsize $(i-j)\delta$};
\node [right] at (6,1.25) {\scriptsize $p \nu(j)$};

\end{tikzpicture}
\caption{Illustration of state space and the transitions for the Erlang-R model with holding.}
\label{fig:QBDIllustration}
\end{figure}

The state space can be partitioned according to its levels, where level $i$ corresponds to a total queue length $N=i$ patients. This results in an infinite-sized matrix consisting of blocks, where each block corresponds to the transition flow from one level to another. Since the only transitions allowed are within the same level or between two adjacent levels in a QBD-process, we obtain a tridiagonal block structure. Each block consists of elements representing the transition rate of one state to another, and therefore each block is a matrix of size at most $(n+1) \times (n+1)$.

For the Erlang-R model with holding this gives the following result. Let $P$ denote the transition matrix of the process $\{N(t),Q_1(t)\}$. We have the boundary levels $\{1,2,...,n\}$ and $P$ is of the form

\[
P = \left( \begin{array}{cccccccccc}
B_{00} & B_{01} & & & & & & & & \\
B_{10} & B_{11} & B_{12} & & & & & & & \\
 & B_{21} & B_{22} & B_{23} & & & & & & \\
 & & \ddots & \ddots &\ddots & & & & & \\
 & & & & & B_{n \, n-1} & & & & \\
 & & & & B_{n-1 \, n} & B_{nn} & A_0 & & & \\
 & & & & & A_2 & A_1 & A_0 & &  \\
 & & & & & & A_2 & A_1 & A_0 & \\
 & & & & & & & \ddots & \ddots & \ddots \\
\end{array} \right),
\]
where $B_{ii} \in \mathbb{R}^{(i+1) \times (i+1)}$, $B_{i \, i-1} \in \mathbb{R}^{(i+1) \times i}$, $B_{i-1 \, i} \in \mathbb{R}^{i \times (i+1)}$, and $A_0,A_1,A_2 \in \mathbb{R}^{(n+1)\times(n+1)}$. The matrices of transition rates for the boundary states are given by
\[
\begin{array}{cc}
B_{00}=(-\lambda), & \qquad	B_{i-1 \, i} = \left( \begin{array}{ccccc}
0 & \lambda & & &  \\
& \ddots & \lambda & & \\
 & & \ddots & \ddots &\\
 & & & 0 & \lambda  \\
\end{array} \right),
\end{array}
\]
\[
 B_{i \, i-1} = \left( \begin{array}{cccc}
0 & & &   \\
 (1-p)\mu & 0 & &  \\
 & (1-p)\nu(2)& \ddots &  \\
 & & \ddots & 0 \\
 & & & (1-p)\nu(i)\\
\end{array} \right)
\]
and
\[
B_{ii} =
{\small \left(
\begin{array}{ccccccccc}
-(\lambda+i \delta) & i \delta & & & &\\
p \mu & -(\lambda+\mu+(i-1)\delta) & (i-1)\delta & & &\\
& \ddots & \ddots & \ddots & & \\
& & p \nu(i-1) & -(\lambda+\nu(i-1)+\delta) & \delta \\
& & & & p \nu(i) & -(\lambda+\nu(i)) \\
\end{array} \right).
}
\]
Moreover, the transition rates are given by
\[
A_0 = \left( \begin{array}{ccccc}
\lambda & & & & \\
& \lambda & & & \\
 & & & \ddots & \\
 & & & & \lambda  \\
\end{array} \right),
\]
\[ A_2 = \left( \begin{array}{ccccccc}
0 & & & & & & \\
& (1-p)\mu & & & & & \\
& & 2(1-p)\mu & & & & \\
& & & \ddots & & & \\
& & & & s(1-p)\mu & & \\
& & & & & \ddots & \\
& & & & & & s(1-p)\mu \\
\end{array} \right)
\]
and
\[
{\small
A_1 =  \left( \arraycolsep=0.55pt
\begin{array}{cccccccc}
-(\lambda+n \delta) & n \delta & & & & & & \\
p \mu & -(\lambda+\mu+(n-1)\delta) & (n-1)\delta & & & & & \\
& \ddots & \ddots & \ddots & & & & \\
& & s p\mu & -(\lambda+s\mu+(n-s)\delta) & (n-s)\delta & & & \\
& & & \ddots & \ddots & \ddots & & \\
& & & & s p\mu & -(\lambda+s\mu+\delta) & \delta & \\
& & & & & s p\mu & -(\lambda+s\mu)\\
\end{array} \right) 
\hspace{-2mm} .
}
\]

\subsection{Stability condition}
\label{app:stability}
From the general theory of QBD processes \cite{Neuts1981} follows that the Markov process $\{N(t),Q_1(t)\}$ is ergodic (stable) if and only if
\begin{equation}
\pi A_0 e < \pi A_2 e,
\label{eq:QBDstableCondition}
\end{equation}
where $e$ is the all one column vector and $\pi=(\pi_0,...,\pi_n)$ is the equilibrium distribution of the Markov process with generator $A^* := A_0+A_1+A_2$. In other words, $\pi$ is such that
\begin{equation}
\begin{array}{ll}
\pi A^* =0, & \pi e =1,
\end{array}
\label{eq:QBDstableProbabilityVector}
\end{equation}
and
\begin{align*}
A^*
= 
{\small
\left(
\begin{array}{cccccccc}
-n \delta & n \delta & & & & & & \\
p \mu & -(p\mu+(n-1)\delta) & (n-1)\delta & & & & & \\
& \ddots & \ddots & \ddots & & & & \\
& & s p\mu & -(ps\mu+(n-s)\delta) & (n-s)\delta & & & \\
& & & \ddots & \ddots & \ddots & & \\
& & & & p s \mu & -(ps\mu+\delta) & \delta & \\
& & & & & p s \mu & -ps\mu\\
\end{array} \right).
}
\end{align*}
Then $\pi$ must satisfy the balance equations
\begin{align*}
- n \delta \pi_0 + p \mu \pi_1 &= 0, \\
(n-j+1)\delta \pi_{j-1} - (p\nu(j) +(n-j)\delta) \pi_j + p \nu(j+1) \pi_{j+1} &= 0, \\
\delta \pi_{n-1} - p s \mu \pi_n &= 0,
\end{align*}
with $\nu(j)=\min\{j,s\}\mu$, and the normalization condition
\[
\sum_{i=0}^n \pi_i=1.
\]
It is readily verified that
\begin{equation}
\pi_i =
\left\{\begin{array}{ll}
\pi_0 \binom{n}{i} \left(\frac{\delta}{p \mu}\right)^i & \textrm{\normalfont for } 0 \leq i \leq s, \\
\pi_0 \binom{n}{i} \frac{i!}{s!} s^{s-i} \left(\frac{\delta}{p \mu}\right)^i & \textrm{\normalfont for } s+1 \leq i \leq n \\
\end{array} \right.
\label{eq:eqdistr}
\end{equation}
with
\begin{align*}
\pi_0= \left(\sum_{i=0}^{s} \binom{n}{i} \left(\frac{\delta}{p \mu}\right)^i + \sum_{i=s+1}^{n} \binom{n}{i} \frac{i!}{s!} s^{s-i} \left(\frac{\delta}{p \mu}\right)^i \right)^{-1}.
\end{align*}
satisfies the balance equations and the normalization condition.
\endproof

\begin{proposition}
The distribution of the closed two-node Jackson network illustrated in Figure~\ref{fig:Jennings} is given by
\begin{equation}
\hat{\pi_i} =
\left\{\begin{array}{ll}
\hat{\pi}_0 \binom{n}{i} \left(\frac{\delta}{p \mu}\right)^i & \textrm{\normalfont for } 0 \leq i \leq s, \\
\hat{\pi_0} \binom{n}{i} \frac{i!}{s!} s^{s-i} \left(\frac{\delta}{p \mu}\right)^i & \textrm{\normalfont for } s+1 \leq i \leq n \\
\end{array} \right.
\label{eq:eqdistrTwistedJennings}
\end{equation}
with
\begin{align*}
\hat{\pi}_0= \left[\sum_{i=0}^{s} \binom{n}{i} \left(\frac{\delta}{p \mu}\right)^i + \sum_{i=s+1}^{n} \binom{n}{i} \frac{i!}{s!} s^{s-i} \left(\frac{\delta}{p \mu}\right)^i \right]^{-1}.
\end{align*}
\label{prop:CriticalTilburgdistr}
\end{proposition}

\proof
We have a two-node closed Jackson network, with probability transition matrix
\[
P = \left(
\begin{array}{cc}
1-p & p \\
1 & 0
\end{array} \right).
\]
Let $r_i(m)$ denote the rate of service when there are $m$ patient at queue $i$, so $r_1(m)=\min\{m,s\}$ and $r_2(m)=m$. The throughput vector $\gamma = (\gamma_1,\gamma_2) \in \mathbb{R_1}^2$ must satisfy $\gamma = \gamma P$ and we find that $\gamma=(p,1)$ suffices. From the general theory of Jackson networks, see \cite{Jackson1963}, it follows that the stationary distribution is given by
\begin{align*}
\pi_i = G^{-1} g_1(i) g_2(n-i)
\end{align*}
with
\begin{align*}
\begin{array}{ll}
g_1(i)= \frac{(\gamma_1/\mu)^i}{\prod_{m=1}^i r_1(m)}, & g_2(n-i)= \frac{(\gamma_2/\delta)^{n-i}}{\prod_{m=1}^{n-i} r_2(m)},
\end{array}
\end{align*}
and normalization constant $G= \sum_{i=0}^n g_1(i) g_2(n-i)$. Then,
\begin{align*}
g_1(i) &=  \left\{\begin{array}{ll}
\frac{1}{i! \mu^i} & \textrm{\normalfont for } 0 \leq i \leq s, \\
\frac{1}{s! s^{i-s} \mu^i} & \textrm{\normalfont for } s+1 \leq i \leq n, \\
\end{array} \right.\\
g_2(n-i) &=\frac{1}{(n-i)!} \left(\frac{p}{\delta}\right)^n \left(\frac{\delta}{p}\right)^i,
\end{align*}
and rewriting the expressions yields~\eqref{eq:eqdistrTwistedJennings}.
\endproof

\subsection{Stationary distribution}
\label{app:StationaryDistributrion}
Assuming that the stability condition is satisfied, we can determine the unique stationary distribution of the Markov process $\{N(t),Q_1(t)\}$. The vector $\pi_i$ can be written as $\pi_{n+i}= \pi_n G^{i}$ for $i=0,1,...$, where $G$ is the minimal nonnegative solution of the non-linear matrix equation
\begin{equation}
A_0+G A_1 + G^2 A_2=0.
\label{eq:MG-G}
\end{equation}
The balance equations can be written as
\[
\begin{array}{ll}
\pi_{i-1} A_0+ \pi_i A_1 + \pi_{i+1} A_2=0, & i=n+1,n+2,...
\end{array}
\]
and using $\pi_{n+i}= \pi_n G^{i-n}$ for $i=0,1,...$, this find
\[
\begin{array}{ll}
\pi_n G^{i-n-1} \left(A_0+ G A_1 + G A_2\right)=0, & i=n+1,n+2,....
\end{array}
\]

\noindent
Moreover, we have the boundary equations
\begin{align*}
\pi_0 B_{00} + \pi_1 B_{10} &= 0 \\
\pi_0 B_{01} + \pi_1 B_{11} + \pi_2 B_{21} &= 0 \\
\pi_1 B_{12} + \pi_1 B_{22} + \pi_2 B_{32} &= 0 \\
 &\vdots&  \\
\pi_{n-2} B_{n-2 \, n-1} + \pi_{n-1} B_{n-1 \, n-1} + \pi_{n} B_{n \, n-1} &= 0 \\
\pi_{n-1} B_{n-1 \, n} + \pi_{n} B_{nn} + \pi_{n+1} A_2 &= 0, 
\end{align*}
along with the normalization equation
\[
1 = \sum_{i=0}^{\infty} \pi_i e = \sum_{i=0}^{n-1} \pi_i e + \pi_n(I-G)^{-1}e,
\]
where we slightly abuse notation by using $e$ as the all ones vector of appropriate size. We note that the matrix $G$ has a spectral radius less than one and therefore $(I-G)$ is invertible.

These equations provide the tools for finding the equilibrium probabilities. Although it is hard to solve $G$ analytically from Equation~\eqref{eq:MG-G}, it is easy to solve numerically by using the following algorithm (matrix-geometric method). Rewriting~\eqref{eq:MG-G} gives
\[
G=-(A_0+G^2 A_2) A_1^{-1},
\]
where $A_1$ is invertible, since it is a transient generator matrix. Let
\[
G_{k+1}=-(A_0+G_k^2 A_2) A_1^{-1},
\]
starting with $G_k=0$. We note that $G_k \uparrow G$ as $k$ grows to infinity \cite{Neuts1981}. Once $||G_{k+1}-G_{k}||_2$ is below a certain preset threshold, we approximate $G$ by $G_{k+1}$.

\section{Proof of Proposition \ref{thm:stochasticordering}}\label{app:stochastic_ordering}

First, note that by definition of the Erlang-R model with holding, in which no more that $n$ patients can be admitted in the ED simultaneously, that $Q_1^h(t)+Q_2^h(t) \leq n = Q_1^J(t) + Q_2^J(t)$ follows directly. 
Therefore, we only consider the relation between the states in the blocking and holding variants Erlang-R model.

As noted Section \ref{sec:Markov_process}, the model with holding can be characterized as a three-dimensional Markov chain $X^h(t) = (H(t),Q^h_1(t),Q^h_2(t))$ in which the components denote the number of holding, needy and content patients respectively. The Erlang-R model with blocking similarly admits a Markov process description, but with two dimensions, namely $X^b(t) = (Q^b_1(t),Q^b_1(t))$.

We prove the result by constructing a coupling between the Markov processes $X^h$ and $X^b$. Let $Z(t) := \left(\hat{X}^h(t),\hat X^b(t)\right) = \left(\hat{H}(t),\hat{Q}_1^h(t),\hat{Q}_2^h(t),\hat{Q}^b_1(t),\hat{Q}^b_2(t)\right)$.

We first define the transition rates of this five-dimensional Markov process, which naturally only depend on the current state of the system. 
After that we show that the transition rates relevant to $\hat{X}^h(t)$ and $X^{h}(t)$ coincide with those of either $X^h(t)$ or $\hat{X}^b(t)$, respectively. The latter implies that the marginal transitions of $\hat{X}^h(t)$ and $X^b(t)$ (and $\hat{X}^b(t)$ and $X^h(t)$) are equal, and hence so are their probability distribution of the Markov processes. 

Let $Z(t) = (h,\qh,\ch,\qy,\cy)$. While defining the reachable states from this state and associated transition rates, we distinguish four transition types, and further differentiate the transition rates depending on the current state.\\
\\*
\textbf{Arrival.}
Arrivals to occur in both models simultaneously, but are handled differently according to the current queue lengths.
\begin{enumerate}
\item If $\qh+\ch < n$ and $\qy+\cy < n$,
\begin{equation}
\label{eq:arr1}
(h,\qh+1,\ch,\qy+1,\cy) \qquad \text{with rate }\l,
\end{equation}
\item if $\qh+\ch = n$ and $\qy+\cy < n$,
\begin{equation}
\label{eq:arr2}
(h+1,\qh,\ch,\qy+1,\cy) \qquad \text{with rate }\l,
\end{equation}
\item if $\qh+\ch < n$ and $\qy+\cy = n$,
\begin{equation}
\label{eq:arr3}
(h,\qh+1,\ch,\qy,\cy) \qquad \text{with rate }\l,
\end{equation}
\item if $\qh+\ch = n$ and $\qy+\cy = n$,
\begin{equation}
\label{eq:arr4}
(h+1,\qh+1,\ch,\qy,\cy) \qquad \text{with rate }\l,
\end{equation}
\end{enumerate}
 
\noindent \textbf{Departure.}
Basically, we align service completions in the two models, but allow a completion occurring solely in either of one of the two models, only if the queue length in this model is strictly larger than in the other one.
\begin{enumerate}
\item If $\qh \geq \qy$ and $h > 0$
\begin{equation}
\label{eq:dep1}
\left\{
\begin{array}{ll}
(h-1,\qh,\ch,\qy-1,\cy) & \text{with rate }(\qy \wedge s)(1-p)\mu,\\
(h-1,\qh,\ch,\qy,\cy) & \text{with rate }[(\qh \wedge s)-(\qy \wedge s)](1-p)\mu.
\end{array}
\right.
\end{equation}

\item If $\qh < \qy$ and $h > 0$
\begin{equation}
\label{eq:dep2}
\left\{
\begin{array}{ll}
(h-1,\qh,\ch,\qy-1,\cy) & \text{with rate }(\qh \wedge s)(1-p)\mu,\\
(h,\qh,\ch,\qy-1,\cy) & \text{with rate }[(\qy \wedge s)-(\qh \wedge s)](1-p)\mu.
\end{array}
\right.\end{equation}

\item If $\qh \geq \qy$ and $h = 0$
\begin{equation}
\label{eq:dep3}
\left\{
\begin{array}{ll}
(0,\qh-1,\ch,\qy-1,\cy) & \text{with rate }(\qy \wedge s)(1-p)\mu,\\
(0,\qh-1,\ch,\qy,\cy) & \text{with rate }[(\qh \wedge s)-(\qy \wedge s)](1-p)\mu.
\end{array}
\right.\end{equation}

\item If $\qh < \qy$ and $h = 0$
\begin{equation}
\label{eq:dep4}
\left\{
\begin{array}{ll}
(0,\qh-1,\ch,\qy-1,\cy) & \text{with rate }(\qh \wedge s)(1-p)\mu,\\
(0,\qh,\ch,\qy-1,\cy) & \text{with rate }[(\qy \wedge s)-(\qh \wedge s)](1-p)\mu.
\end{array}
\right.\end{equation}
\end{enumerate}

\noindent\textbf{Become content} 
The differentiation between transitions is similar to those in the \textit{departure} transition type.
\begin{enumerate}
\item If $\qh \geq \qy$,
\begin{equation}
\label{eq:con1}
\left\{
\begin{array}{ll}
(h,\qh-1,\ch+1,\qy-1,\cy+1) & \text{with rate }(\qy \wedge s)p\mu,\\
(h,\qh-1,\ch+1,\qy,\cy) & \text{with rate }[(\qh \wedge s)-(\qy \wedge s)]p\mu.
\end{array}
\right.\end{equation}

\item If $\qh < \qy$,
\begin{equation}
\label{eq:con2}
\left\{
\begin{array}{ll}
(h,\qh-1,\ch+1,\qy-1,\cy+1) & \text{with rate }(\qh \wedge s)p\mu,\\
(h,\qh,\ch,\qy-1,\cy+1) & \text{with rate }[(\qy \wedge s)-(\qh \wedge s)]p\mu.
\end{array}
\right.\end{equation}
\end{enumerate}
\noindent
\textbf{Become needy}
 \begin{enumerate}
\item If $\ch \geq \cy$,
\begin{equation}
\label{eq:ne1}
\left\{
\begin{array}{ll}
(h,\qh+1,\ch-1,\qy+1,\cy-1) & \text{with rate } \cy\delta,\\
(h,\qh+1,\ch-1,\qy,\cy) & \text{with rate }(\ch-\cy)\delta,\\
\end{array}
\right.\end{equation}

\item If $\ch < \cy$,
\begin{equation}
\label{eq:ne2}
\left\{
\begin{array}{ll}
(h,\qh+1,\ch-1,\qy+1,\cy-1) & \text{with rate } \ch\delta,\\
(h,\qh,\ch,\qy+1,\cy-1) & \text{with rate }(\cy-\ch)\delta,\\
\end{array}
\right.\end{equation}
\end{enumerate}

This set of transitions defines the dynamics of the Markov process $Z(t) = (\hat{X}^h(t),\hat{X}^b(t))$.
Let us now restrict our attention to the transitions in which (at least one of the) first three coordinates of $Z(t)$ changes, that is, the marginal transitions of the process $\hat{X}^h$. 
Let $\hat{X}^h(t) = (h,\qh,\ch)$, then according to the transition scheme above, $\hat{X}^h$ moves to state

\begin{enumerate}
\item If $\qh+\ch < n$ (and hence necessarily $h=0$),
\[
\left\{
\begin{array}{ll}
(0,\qh+1,\ch) & \text{with rate } \l,\\
(0,\qh-1,\ch) & \text{with rate }(\qh\wedge s)(1-p)\mu,\\
(0,\qh-1,\ch+1) & \text{with rate }(\qh\wedge s)p\mu,\\
(0,\qh+1,\ch-1) & \text{with rate }\ch \delta.
\end{array}
\right.\]

\item if $\qh+\ch = n$ and $h=0$,
\[
\left\{
\begin{array}{ll}
(1,\qh,\ch) & \text{with rate } \l,\\
(0,\qh,\ch) & \text{with rate }(\qh\wedge s)(1-p)\mu,\\
(0,\qh-1,\ch+1) & \text{with rate }(\qh\wedge s)p\mu,\\
(0,\qh+1,\ch-1) & \text{with rate }\ch \delta.
\end{array}
\right.\]

\item if $h>0$ (and hence necessarily $\qh+\ch = n$),
\[
\left\{
\begin{array}{ll}
(h+1,\qh,\ch) & \text{with rate } \l,\\
(h-1,\qh,\ch) & \text{with rate }(\qh\wedge s)(1-p)\mu,\\
(h,\qh-1,\ch+1) & \text{with rate }(\qh\wedge s)p\mu,\\
(h,\qh+1,\ch-1) & \text{with rate }\ch \delta.
\end{array}
\right.\]
\end{enumerate}
One can check that these transitions indeed coincide with the transitions in the original holding model, hence $\hat{X}^h(t) \ed X^h(t)$.

Similarly, when the focusing on transitions of $Z(t)$ that are relevant for $\hat{X}^b(t)$, we deduce the following transition scheme. If $\hat{X}^b(t) = (\qy,\cy)$, then the next move according to the transitions of $Z(t)$ is
\[
\left\{
\begin{array}{ll}
(\qy+1_{\{\qy + \cy < n\}},\cy) & \text{with rate } \l,\\
(\qy-1,\cy) & \text{with rate }(\qy\wedge s)(1-p)\mu,\\
(\qy-1,\cy+1) & \text{with rate }(\qy\wedge s)p\mu,\\
(\qy+1,\cy-1) & \text{with rate }\cy \delta.
\end{array}
\right.\]

These transition rates clearly coincide with the original Erlang-R model with blocking, and also hence $\hat{X}^b(t) \ed X^h(t)$.

Next, we show that under this coupling scheme we have that if $\hat{H}(0) = 0$,  $\hat{Q}_1^h(0)=\hat{Q}_1^b(0)$ and $\hat{Q}_1^h(0)=\hat{Q}^b(0)$ then for all $t\geq 0$, $Z(t)$ satisfies the hypothesis:
\begin{itemize}
\item[(i)] $\hat{Q}_1^b(t) + \hat{Q}_2^b(t) \leq \hat{Q}_1^h(t) + \hat{Q}_2^h(t)$,
\item[(ii)] $\hat{Q}_2^b(t) \leq \hat{Q}_2^h(t)$,
\item[(iii)] $\hat{Q}_1^b(t) \leq \hat{Q}_1^h(t) + H(t)$.
\end{itemize}
We do so by induction on the next state reached after a transition of the joint Markov process $Z=(\hat{X}^h,\hat{X}^b)$.
First of all, $Z(0)$ clearly satisfies (i)-(iii). 
Next, assume $Z(t^-) = (h,\qh,\ch,\qy,\cy)$ satisfies the hypothesis and a transition occurs at $t$. 
We show that under the specified coupling scheme, the state reached after the next transition, $Z(t)$ must satisfy (i)-(iii) as well. To do so, we differentiate between the four types of transitions that could occur: arrival, departure, become content and become needy.

\noindent\textbf{Arrival.}\\
Recall that under our coupling scheme an arrival always occurs in both the holding and blocking model simultaneously, see \eqref{eq:arr1}--\eqref{eq:arr4}. Furthermore, $\ch$ and $\cy$ are unchanged during this transition, rendering (ii) trivial.

By hypothesis $\qy + \cy \leq \qh+\cy$, hence the event $\qh+\ch < n$ and $\qh+\cy =n$, with resulting state $(0,\qh+1,\ch,\qy,\cy)$,  can be excluded from our analysis
We check the conditions for the remaining three cases.

\begin{enumerate}[noitemsep]
\item If $Z(t)= (0,\qh+1,\ch,\qy+1,\cy)$, then $\qy + \cy<n$ and $\qh + \ch <n$.
\begin{itemize}[noitemsep]
\item[(i)] $\QY +\CY = \qy+\cy+1 \less[i] \qh+\ch+1 =\QH+\CH$.
\item[(iii)] $\QY  = \qy+1 \less[iii] \qh+1 = \QH = \QH+\hat{H}(t)$.
\end{itemize}

\item If $Z(t)= (h+1,\qh,\ch,\qy+1,\cy)$, then $\qy + \cy<n$ and $\qh + \ch =n$.
\begin{itemize}[noitemsep]
\item[(i)] $\QY +\CY = \qy+\cy+1 \leq n = \qh+\ch =\QH+\CH$.
\item[(iii)] $\QY = \qy+1 \less[iii] \qh +1= \QH +\hat{H}(t)$.
\end{itemize}

\item If $Z(t)= (h+1,\qh,\ch,\qy,\cy)$, then $\qy + \cy = \qh+\ch=n$.
\begin{itemize}[noitemsep]
\item[(i)] $\QY +\CY = \qy+\cy \less[i] \qh+\ch =\QH+\CH$.
\item[(iii)] $\QY = \qy \less[iii] \qh+h < \qh+h+1 = \hat{H}(t)$.
\end{itemize}
\end{enumerate}

\noindent\textbf{Departure.}
By carefully examining the possible state transitions of $Z(t)$ following a departure, we list six reachable states. However, by (iii), we have that if $h=0$, then $\qy \leq \qh$, which excludes the state $(0,\qh,\ch,\qy,\cy)$ in \eqref{eq:dep4} from the reachability graph. 
We check the remaining states for conditions (i)--(iii). Again, during a departure, $\cy$ and $\ch$ are unchanged, so (ii) is automatically satisfied by the induction hypothesis.
\begin{enumerate}[noitemsep]
\item If $Z(t) = (h-1,\qh,\ch,\qy-1,\cy)$, then $h>0$. 
\begin{itemize}
\item[(i)] $\QY+\CY = \qy+\cy-1 \less[i] \qh+\ch-1 < \qh+\ch = \QH+\CH$. 
\item[(iii)] $\QY = \qy-1 \less[iii] \qh + h-1 = \QH +\hat{H}(t)$.
\end{itemize}

\item If $Z(t) = (h-1,\qh,\ch,\qy,\cy)$, then $h>0$ and $\qh \geq \qy$ (*). 
\begin{itemize}
\item[(i)] $\QY+\CY = \qy+\cy \less[i] \qh+\ch = \QH+\CH$.
\item[(iii)] $\QY = \qy \less[*] \qh-1 \leq \qh+h-1 = \QH+\hat{H}(t)$.
\end{itemize}

\item If $Z(t) = (h,\qh,\ch,\qy-1,\cy)$, then $h>0$ and $\qh < \qy$ (*). 
\begin{itemize}
\item[(i)] $\QY+\CY = \qy+\cy-1 < \qy+\cy \less[i] \qh+\ch = \QH+\CH$.
\item[(iii)] $\QY = \qy-1 < \qy \less[*] \qh + h = \QH+\hat{H}(t)$.
\end{itemize}

\item If $Z(t) = (h,\qh-1,\ch,\qy-1,\cy)$, then $h=0$. 
\begin{itemize}
\item[(i)] $\QY+\CY = (\qy-1)+\cy-1 <  \less[i] (\qh-1)+\ch = \QH+\CH$.
\item[(iii)] $\QY  = \qy-1 \less[iii] \qh-1 = \QH + \hat{H}(t)$.
\end{itemize}

\item If $Z(t) = (0,\qh-1,\ch,\qy,\cy)$, then $h=0$ and $\qh>\qy$ (*). 
\begin{itemize}
\item[(i)] $\QY+\CY = \qy+\cy  \less[i] (\qh-1)+\cy \less[ii] (\qh-1)+\ch = \QH+\CH$.
\item[(iii)] $\QY = qy \less[*] \qh-1 =\QH+ \hat{H}(t)$.
\end{itemize}

\end{enumerate}

\noindent\textbf{Content start.}
On the event of a patient becoming content, it is clear that the sums $\QH+\CH$ and $\QY+\CY$ and $H(t)$ are unaffected. This means that (i) is directly satisfied by the induction hypothesis.
According to \eqref{eq:con1}--\eqref{eq:con2}, three states can be reached.  
\begin{enumerate}[noitemsep]
\item If $Z(t) = (h,\qh-1,\ch+1,\qy-1,\cy+1)$,
\begin{itemize}[noitemsep]
\item[(ii)] $\CY = \cy+1 \less[ii] \ch+1 = \CH$. 
\item[(iii)] $\QY = \qy-1 \less[iii] \qh+h-1 = \QH+\hat{H}(t)$.
\end{itemize}
\item If $Z(t) = (h,\qh-1,\ch+1,\qy,\cy)$, then $\qh > \qy$ (*),
\begin{itemize}[noitemsep]
\item[(ii)] $\CY = \cy \less[ii] \ch < \ch+1 = \CH$. 
\item[(iii)] $\QY = \qy \less[iii] \qh+h < \qh+1+h = \QH+ \hat{H}(t)$.
\end{itemize}
\item If $Z(t) = (h,\qh,\ch,\qy-1,\cy+1)$, then $\qy > \qh$ and hence by (iii) $h > 0$. The latter is only possible if $\qh+\ch=n$ (*),
\begin{itemize}[noitemsep]
\item[(ii)] $\CY = \cy+1 \leq n-\qy+1 = (\qh+\ch)-\qy+1 \less[*] \ch = \CH$. 
\item[(iii)] $\QY = \qy-1 < \qh+h-1 \less[*] \qh+h = \QH+\hat{H}(t)$.
\end{itemize}
\end{enumerate}

\noindent \textbf{Become needy.}\\
Just as in the event of content start, the sums $\QH+\CH$ and $\QY+\CY$ and $H(t)$ are unaffected, whereby  (i) is directly satisfied by the induction hypothesis.
By (ii), we have $\ch \geq \cy$. This excludes the state $(h,\qh,\ch,\qy+1,\cy-1)$ from being reached, see \eqref{eq:ne2}. 
We check the remaining to possibilities.

\begin{enumerate}
\item If $Z(t) = (h,\qh+1,\ch-1,\qy+1,\cy-1)$. 
\begin{itemize}[noitemsep]
\item[(ii)] $\CY = \cy-1 \less[ii] \ch-1 = \CH$. 
\item[(iii)] $\QY = \qy+1 \less[iii] \qh+h+1 = \QH+\hat H(t)$.
\end{itemize}
\item If $Z(t) = (h,\qh+1,\ch-1,\qy,\cy)$, then $\ch > \cy$ (*). 
\begin{itemize}
\item[(ii)] $\CY = \cy \less[*] \ch-1 = \CH$.
\item[(iii)] $\QY = \qy \less[iii] \qh+h < \qh+1+h =\QH + \hat H(t)$.
\end{itemize}
\end{enumerate}

Hence, the state reached after any feasible transition under the coupling scheme satisfies the conditions (i)--(iii).
Thus we conclude that the joint process $(\hat{H}(t),\QH,\CH,\QY,\CY)$ adheres to (i)--(iii) for all $t$. Consequently, we have that (i) implies
\begin{align*}
\P\left(Q_1^b(t) + Q_2^b(t) \geq k\right) &= \P\left(Q_1^b(t) + Q_2^b(t) \geq k\right)\\
&=\sum_{j=0}^n \P\left( \QY + \CY \geq k , \QH+\CH = j \right) \\
&=\sum_{j=k}^n \P\left( \QY + \CY \geq k , \QH+\CH = j \right) \\
&\leq \sum_{j=h}^n \P\left( \QH+\CH = j \right)\\
&= \P\left( Q_1^h(t) + Q_2^h(t) \geq k\right) = \P\left(Q_1^h(t) + Q_2^h(t) \geq k\right).
\end{align*}
The other two orderings follow similarly.

\begin{remark}
Note that under this coupling scheme we cannot get the ordering $\QH(t) \geq \QY(t)$ for all $t\geq 0$. A minimal counter example occurs for $s=n=1$. Let $Z(0) = ((0,0,0),(0,0))$. First, two arrivals occur, such that state $((1,1,0),(1,0))$ is reached, followed by a departure transition, yielding $((0,1,0),(0,0))$. Next, the one patient left in the model with holding system becomes content, so that we obtain $((0,0,1),(0,0))$. 
At this stage, if an arrival occurs, the arriving patient will be put in the holding queue in the model with holding, and admitted to nurse queue in the model with blocking. Hence we end up in state $((1,0,1),(1,0))$, in which $\QH < \QY$. 
\end{remark}

\section{Proof of Theorem \ref{thm:limits_YT} - Performance measures of Erlang-R with blocking}\label{app:proof_block}
For convenience we state the theorem here again.
\begin{theorem}
Let $s$ and $n$ scale as in \eqref{eq:twofoldscaling} with $\b,\g>0$ as $\l\to\infty$. Then, if $\b \neq 0$
\begin{equation}
g^b(\b,\g)
:= \lim_{\l\to\iy} \P^b({\rm delay}) = 
\left(1 + 
\frac{ \b \, \int_{-\iy}^\b \Phi\left(\frac{\g-t\sqrt{r}}{\sqrt{1-r}}\right)\, d\Phi(t) }
{\phi(\b)\Phi(\eta) -  \phi(\sqrt{\b^2+\eta^2}){\rm e}^{\tfrac{1}{2} \omega^2} \Phi(\omega)}
\right)^{-1},
\end{equation}
\begin{equation}
f^b(\b,\g) 
:= \lim_{\l\to\iy} \sqrt{R_1}\cdot\P^b({\rm block}) = 
\frac{
\sqrt{r}\phi(\g)\Phi(-\omega\sqrt{r}) + \phi(\sqrt{\b^2+\eta^2})\,{\rm e}^{\frac{1}{2} \omega^2} \Phi(\omega)
}{
\int_{-\iy}^\b \Phi\left(\frac{\g-t\sqrt{r}}{\sqrt{1-r}}\right)\, d\Phi(t) +
\frac{\phi(\b)\Phi(\eta)}{\b} -  \frac{\phi(\sqrt{\b^2+\eta^2})}{\b}{\rm e}^{\tfrac{1}{2} \omega^2} \Phi(\omega)
 },
\end{equation}
\begin{equation}
h^b(\b,\g) := \lim_{\l\to\iy} \sqrt{R_1}\cdot\E[W]
=
\frac{
\frac{\phi(\b)\Phi(\eta)}{\b^2} +
 \left(\frac{\b}{r}-\frac{\g}{\sqrt{r}}-\frac{1}{\b}\right)\,\frac{\phi(\sqrt{\eta^2+\b^2})}{\b}\, {\rm e}^{\tfrac{1}{2}\omega^2}\, \Phi(\omega) 
 - \sqrt{\frac{1-r}{r}}\,\frac{\phi(\b)\phi(\eta)}{\b}
}{
\int_{-\iy}^\b \Phi\left(\frac{\g-t\sqrt{r}}{\sqrt{1-r}}\right)\, d\Phi(t) +
\frac{\phi(\b)\Phi(\eta)}{\b} -  \frac{\phi(\sqrt{\b^2+\eta^2})}{\b}{\rm e}^{\tfrac{1}{2} \omega^2} \Phi(\omega)
 },
\end{equation}
%
and if $\b=0$,
\begin{equation}
g^b_0(\g) 
:= \lim_{\l\to\iy} \P^b({\rm delay})
\left(1+
\frac{
\int_{-\iy}^0 \Phi\left(\frac{\g-t\sqrt{r}}{\sqrt{1-r}}\right)\, d\Phi(t)
}{
\sqrt{\frac{1-r}{r}} \frac{1}{\sqrt{2\pi}}\,\left(\eta \,\Phi(\eta) + \phi(\eta) \right)
}
\right)^{-1}
\end{equation}
\begin{equation}
f^b_0(\g) 
:= \lim_{\l\to\iy} \sqrt{R_1}\cdot\P^b({\rm block}) = 
\frac{
\sqrt{r}\,\phi(\g)\Phi(-\omega\sqrt{r}) + \frac{1}{\sqrt{2\pi}} \Phi(\eta)
}{
\int_{-\iy}^\b \Phi\left(\frac{\g-t\sqrt{r}}{\sqrt{1-r}}\right)\, d\Phi(t) +
\sqrt{\frac{1-r}{r}} \frac{1}{\sqrt{2\pi}}\,\left(\eta \,\Phi(\eta) + \phi(\eta) \right)
 },
\end{equation}

\begin{equation}
h_0^b(\g) := \lim_{\l\to\iy} \sqrt{R_1}\cdot\E[W] 
= \frac{1}{2\mu}\, \frac{ \left( \gamma^2/r+1\right) \Phi(\eta) + \eta \phi(\eta) }
{ \frac{r}{1-r} \sqrt{2\pi} \int_{-\infty}^0 \Phi\left(\frac{\g-t\sqrt{r}}{\sqrt{1-r}}\right)\, d\Phi(t) + \sqrt{\frac{r}{1-r}} \left(\eta \Phi(\eta)+\phi(\eta)\right)},
\end{equation}
where $\eta = \frac{\g - \b\sqrt{r}}{\sqrt{1-r}}$ and $\omega := \frac{\g - \b/\sqrt{r}}{\sqrt{1-r}}$.
\end{theorem}

In this appendix we prove the heavy-traffic approximations of the system-measures introduced in Theorem \ref{thm:limits_YT}. As a first stage we present and prove four lemmas.

\begin{lemma} \label{lem1}
Let the variables $\lambda$, $s$ and $n$ tend to $\infty$ simultaneously and satisfy the QED scaling conditions in \eqref{eq:twofoldscaling} with $\beta \neq 0$. 
Define $B_1$ as the expression 
\begin{equation*}
B_1=\frac{e^{-R}}{s!} R_1^{s} \frac{1}{1-\rho} \sum_{l=0}^{n-s-1}  \frac{1}{l!} R_2^{l} e^{-R_2}. 
\end{equation*}
Then
\[ \lim_{\lambda \rightarrow \infty} B_1 = \frac{\phi(\beta)\Phi(\eta)}{\beta } .\]
\end{lemma}

\begin{lemma} \label{lem2}
Let the variables $\lambda$, $s$ and $n$ tend to $\infty$ simultaneously and satisfy the QED scaling conditions in \eqref{eq:twofoldscaling} with $\beta \neq 0$. 
Define $B_2$ as the expression 
\begin{equation*}
B_2 = \frac{e^{-\left(R_1 + R_2 \right) }}{s!} R_1^{s} \frac{\rho^{n-s}}{1-\rho} \sum_{l=0}^{n-s-1} \frac{1}{l!} \left(\frac{R_2}{\rho} \right)^{l}. 
\end{equation*}
Then
\[ \lim_{\lambda \rightarrow \infty} B_2 = \frac{\phi(\sqrt{\eta^2 + \beta^2})}{\beta} e^{\frac{1}{2} \omega ^2} \Phi(\omega) .\]
\end{lemma}

\begin{lemma} \label{lem3}
Let the variables $\lambda$, $s$ and $n$ tend to $\infty$ simultaneously and satisfy the QED scaling conditions in \eqref{eq:twofoldscaling}. 
Define $A$ as the expression 
\begin{equation}
\label{eq:A}
	A = \sum_{\substack{i,j|i \leq s,\\ i+j \leq n-1}} \frac{1}{i!j!} R_1^i  R_2^j  e^{-\left(R_1 + R_2 \right)}.
\end{equation}
Then
\[ \lim_{\lambda \rightarrow \infty} A =  \int_{-\infty}^{\beta} \Phi\left( \eta + (\beta-t) \sqrt{\frac{\delta }{\mu p}} \right) d\Phi(t) .\]
\end{lemma}

\begin{lemma} \label{lem4}
Let the variables $\lambda$, $s$ and $n$ tend to $\infty$ simultaneously and satisfy the QED scaling conditions in \eqref{eq:twofoldscaling} with $\beta = 0$. Define $B$ as the expression 
\begin{equation*}
	B = e^{-(R_1 + R_2)} \frac{1}{s!}  {R_1}^{s} \sum_{j=0}^{n-s-1} \frac{1}{j!} {R_2}^{j} \sum_{i=0}^{n-s-j-1} \rho^{i} .
\end{equation*}
Then
\[ \lim_{\lambda \rightarrow \infty} B = \sqrt{\frac{\mu p}{\delta}} \frac{1}{\sqrt{2\pi}} \left(\eta \Phi(\eta) + \phi(\eta) \right) = .\]
\end{lemma}

\subsection{Proof of Lemma \ref{lem1}}
\proof
By using Stirling's formula $\left(s! \approx \sqrt{2\pi s} \left(\frac{s}{e}\right)^s \right)$, and QED assumption that $\sqrt{s}(1-R_1/s) \rightarrow \beta$ as $\lambda \rightarrow \infty$, one obtains for $B_1$:
\begin{equation*}
\begin{split}
B_1 &\approx \frac{e^{s-R}}{\sqrt{2\pi s} } \rho ^{s} \frac{\sqrt{s}}{\beta} \sum_{l=0}^{n-s-1}  \frac{1}{l!} \left(R_2  \right)^{l} e^{- R_2 } = \frac{e^{s-R} \rho ^{s}}{\sqrt{2\pi}\beta }    \sum_{l=0}^{n-s-1}  \frac{1}{l!} R_2^{l} e^{-R_2}  \\
&= \frac{e^{s(1-\rho)}}{\sqrt{2\pi}\beta }  \rho^{s}  P(X_{\lambda} \leq n-s-1)
\end{split}
\end{equation*}
where $\rho=\frac{\lambda}{(1-p)s\mu}$, and $X_{\lambda}$ is a random variable with the Poisson distribution with parameter $R_2$. When $\lambda \rightarrow \infty$, $R_2 \rightarrow \infty$ too, since $p$, and $\delta$  are fixed.  Note that
\begin{equation*}
	\P (X_{\lambda} \leq n-s-1) = \P\left(\frac{X_{\lambda}-R_2}{\sqrt{R_2}} \leq \frac{n-s-1-R_2}{\sqrt{R_2}}\right)
\end{equation*}
Now we need to find the limit for the following fraction
\[\frac{n-s-R_2}{\sqrt{R_2}}\]
as $\lambda \rightarrow \infty$ using assumption that $\sqrt{s}(1-R_1/s) \rightarrow \beta$ as $\lambda \rightarrow \infty$.
\begin{equation}
\begin{split}
	&\lim_{\lambda \rightarrow \infty} \frac{n-s-R_2}{\sqrt{R_2}} = \lim_{\lambda \rightarrow \infty} \frac{n- \frac{R_1}{r} -s + R_1 }{\sqrt{\frac{R_1}{r}-R_1}} =  \frac{\gamma \sqrt{\frac{R_1}{r}}- \beta \sqrt{R_1}}{\sqrt{\frac{R_1}{r}-R_1}}  = \frac{\gamma - \beta \sqrt{r}}{\sqrt{1-r}}. 
\end{split}
\end{equation}
Hence, define $\eta=\frac{\gamma - \beta \sqrt{r}}{\sqrt{1-r}}$.
Thus, when $\lambda \rightarrow \infty$, by the Central Limit Theorem (Normal approximation to
Poisson) we have
\begin{equation*}
	\left(\frac{X_{\lambda}-R_2}{\sqrt{R_2}}\right) \Rightarrow N(0,1)
\end{equation*}
and due to assumption QED (i) of the lemma we get 
\begin{equation}
\label{lm1_eq1}
	\P(X_{\lambda} \leq n-s-1) \rightarrow \P(N(0,1) \leq \eta) = \Phi(\eta), ~~ as ~\lambda \rightarrow \infty 
\end{equation}
where $N(0,1)$ is a standard normal random variable with distribution function $\Phi(\cdot)$. 
It follows thus that
\begin{equation*}
B_1 \approx \frac{e^{s(1-\rho)}}{\sqrt{2\pi}\beta } \rho^{s} \Phi(\eta) =\frac{e^{s(1-\rho+ \ln \rho)}}{\sqrt{2\pi}\beta } \Phi(\eta).
\end{equation*}
Making use of the expansion
\[\ln \rho = \ln(1-(1-\rho)) = -(1-\rho)-\frac{(1-\rho)^2}{2}+ o(1 -\rho)^2, ~~ (\rho \rightarrow 1)\]
one obtains
\begin{equation*}
B_1 \approx \frac{e^{s(1-\rho -(1-\rho)-\frac{(1-\rho)^2}{2})}}{\sqrt{2\pi}\beta } \Phi(\eta) = \frac{e^{-\frac{s(1-\rho)^2}{2}}}{\sqrt{2\pi}\beta } \Phi(\eta)
\end{equation*}
by QED assumption that $\sqrt{s}(1-R_1/s) \rightarrow \beta$, it is clear that $s(1-\rho)^2 \rightarrow \beta^2$, when $\lambda \rightarrow \infty$. This implies
\begin{equation*}
\lim_{\lambda \rightarrow \infty} B_1 = \frac{\phi(\beta)\Phi(\eta)}{\beta } 
\end{equation*}
where $\phi(\cdot)$ is the standard normal density function, and $\Phi(\cdot)$ is the standard
normal distribution function. This proves Lemma \ref{lem1}.
\endproof

\subsection{Proof of Lemma \ref{lem2}}
Again according to Stirling's formula, and the QED assumption that $\frac{n-s-R_1-R_2}{\sqrt{R_1+R_2}} \rightarrow \eta$ as $\lambda \rightarrow \infty$, one obtains for $B_2$:
\begin{equation*}
\begin{split}
B_2 &\approx \frac{e^{s-R_1-R_2 }}{\sqrt{2\pi s}}  \frac{\rho^{n}}{1-\rho} \sum_{l=0}^{n-s-1} \frac{1}{l!} \left(\frac{R_2}{ \rho} \right)^{l} = \frac{e^{s(1-\rho)-R_2}}{\sqrt{2\pi s}} \frac{\sqrt{s}\rho^{n}}{\beta} e^{\frac{R_2}{ \rho} } \sum_{l=0}^{n-s-1} \frac{1}{l!} \left(\frac{R_2}{ \rho} \right)^{l} e^{ - \frac{R_2}{ \rho} } \\
&= \frac{e^{s(1-\rho)+ R_2 \left(\frac{1-\rho}{\rho}\right)} }{\sqrt{2\pi } \beta} \rho^{n} \P(Y_{\lambda}\leq n-s-1) ,
\end{split}
\end{equation*}
where $\rho=\frac{\lambda}{(1-p)s\mu}$,   and $Y_{\lambda}$ is a random variable with the Poisson distribution with parameter $\frac{R_2}{\rho}$. 
Note that
\begin{equation*}
	\P(Y_{\lambda} \leq n-s-1) = \P\left(\frac{Y_{\lambda}-\frac{R_2}{\rho}}{\sqrt{\frac{R_2}{\rho}}} \leq \frac{n-s-1-\frac{R_2}{\rho}}{\sqrt{\frac{R_2}{\rho}}}\right).
\end{equation*}
Now we need to find the limit for the following fraction
$\frac{n-s-\frac{R_2}{\rho}}{\sqrt{\frac{R_2}{\rho}}}$
as $\lambda \rightarrow \infty$ using assumption QED (i).
\begin{equation}
\label{eq:lem2_app_n-s}
\begin{split}
	&\lim_{\lambda \rightarrow \infty} \frac{n-s-\frac{R_2}{\rho}}{\sqrt{\frac{R_2}{\rho}}} = \lim_{\lambda \rightarrow \infty} \frac{\eta \sqrt{R_2} +R_2 -\frac{R_2}{\rho}} {\sqrt{\frac{R_2}{\rho}}} = \lim_{\lambda \rightarrow \infty} \eta \sqrt{\rho}+ \frac{\sqrt{R_2}(\rho-1)}{\sqrt{\rho}} \\
	&= \eta - \lim_{\lambda \rightarrow \infty} \sqrt{\frac{s p \mu}{\delta}}(1-\rho) = \eta - \sqrt{\frac{ p\mu}{\delta }} \beta. 
\end{split}
\end{equation}
Denote $\omega = \eta - \beta \sqrt{\frac{ p\mu}{\delta }}  $. Thus, when $\lambda \rightarrow \infty$, by the Central Limit Theorem (Normal approximation to
Poisson) we have
\begin{equation*}
\left(\frac{Y_{\lambda}-\frac{R_2}{\rho}}{\sqrt{\frac{R_2}{\rho}}}\right) \Rightarrow N(0,1)
\end{equation*}
and
\begin{equation*}
	\P(Y_{\lambda} \leq n-s-1) \rightarrow \P(N(0,1) \leq \omega) = \Phi(\omega), \textnormal{ as } \lambda \rightarrow \infty ,
\end{equation*}
where $N(0,1)$ is a standard normal random variable with distribution function $\Phi$. It follows thus that
\begin{equation*}
B_2 \approx \frac{e^{s(1-\rho)+R_2 \left(\frac{1-\rho}{\rho}\right)}}{\sqrt{2\pi } \beta} \rho^{n} \Phi(\omega) = \frac{e^{s(1-\rho)+R_2\left(\frac{1-\rho}{\rho}\right)+n \ln \rho}}{\sqrt{2\pi } \beta} \Phi(\omega).
\end{equation*}
Making use of the expansion
\[\ln \rho = \ln(1-(1-\rho)) = -(1-\rho)-\frac{(1-\rho)^2}{2}+ o(1 -\rho)^2, ~~ (\rho \rightarrow 1)\]

 and using our assumptions that as $\lambda \rightarrow \infty$: $\rho \rightarrow 1$,  
$s \approx R_1 + \beta \sqrt{R_1},$ $ n \approx \frac{R_1}{r}+\gamma \sqrt{\frac{R_1}{r}}$, and $ n-s \approx R_2+\eta \sqrt{R_2}$, one obtains 
\begin{equation*}
\begin{split}
&s(1-\rho)+R_2\left(\frac{1-\rho}{\rho}\right)+n \ln \rho  = s(1-\rho)+R_2 \left(\frac{1-\rho}{\rho}\right) - n \left(1-\rho+\frac{(1-\rho)^2}{2}\right) \\
&= -\left(n-s-\frac{R_2}{\rho}\right) (1-\rho)-\frac{n(1-\rho)^2}{2} \approx -\left(\eta \sqrt{R_2}+R_2-\frac{R_2}{\rho} \right) (1-\rho)-\frac{n(1-\rho)^2}{2} \\
&= \left(\frac{R_2}{\rho}-\frac{n}{2}\right)(1-\rho)^2 -\eta \sqrt{R_2}(1-\rho) \approx \left(\frac{R_2}{\rho}-\frac{1}{2}\left(\frac{R_1}{r} +\gamma  \sqrt{\frac{R_1}{r}}\right) \right) (1-\rho)^2 -\eta \sqrt{R_2}(1-\rho) \\
&= \left(\frac{R_2}{\rho}-\frac{R_1} {2r} \right) (1-\rho)^2 - \frac{1}{2}\gamma \sqrt{\frac{R_1}{r}}(1-\rho)^2 -\eta \sqrt{R_2}(1-\rho)\\
&= \left(\frac{R_2}{\rho}-\frac{R_1} {2r} \right) \frac{\beta^2 }{R_!} - \frac{1}{2}\gamma \sqrt{\frac{R_1}{r}}\frac{\beta^2 }{R_1} - \eta \sqrt{R_2}\beta \sqrt{\frac{1}{R_1}} \\
&\approx \frac{p \mu}{\delta \rho}\beta^2  - \frac{1}{2}\left(\frac{p \mu}{\delta }  + 1\right)\beta^2  - \eta \beta \sqrt{\frac{p\mu}{\delta } } = -\frac{1}{2}(\eta^2 + \beta^2) + \frac{1}{2} \omega ^2. \\
\end{split}
\end{equation*}
Therefore,
\begin{equation*}
\begin{split}
&\lim_{\lambda \rightarrow \infty} B_2 \approx \frac{e^{s(1-\rho)+ R_2 \left(\frac{1-\rho}{\rho}\right)+n \ln \rho}}{\sqrt{2\pi } \beta} \Phi(\omega) \approx \frac{e^{-\frac{1}{2}(\eta^2 + \beta^2) + \frac{1}{2} \omega ^2}}{\sqrt{2\pi } \beta} \Phi(\omega) = \frac{\phi(\sqrt{\eta^2 + \beta^2})}{\beta} e^{\frac{1}{2} \omega ^2} \Phi(\omega).
\end{split}
\end{equation*}
This proves Lemma \ref{lem2}.

\subsection{Proof of Lemma \ref{lem3}}

We will find the asymptotic behavior of $A$ by finding its lower and upper bounds. Let us consider a partition $\{s_h\}_{h=0}^l$ of the interval $[0,s]$.
\[ s_h = s - h \tau ,~~ h=0,1,...,\ell;~~ s_{\ell+1} = 0 \]
where $ \tau = \left[ \epsilon \sqrt{R_1} \right]$, $\epsilon$ is an arbitrary non-negative real and $\ell$ is a positive integer.

If $\lambda$ and $s$ tend to infinity and satisfy the QED assumption that $\frac{n-s-R_1-R_2}{\sqrt{R_1+R_2}} \rightarrow \eta$ as $\lambda \rightarrow \infty$, then $\ell<\frac{s}{\tau}$ for $\lambda$ big enough and all the $s_h$ belong to $[0, s];~~ h=0,1,...,\ell$. 
Emphasize that the length $\tau$ of every interval $[s_{h-1}, s_h ]$ depends on $\lambda$. 
The variable $A$ is given by the formula \eqref{eq:A}. Let us consider a lower estimate for $A$ given by the following sum:
\begin{equation}
\label{eq_lem3_xi_1}
\begin{split}
	A \geq A_1 &= \sum_{h=0}^{\ell}  \sum_{i=s_{h+1}}^{s_h} \frac{1}{i!} R_1^i e^{-R_1} \cdot  \sum_{j=0}^{n-s_h-1} \frac{1}{j!} R_2^j e^{-R_2} \\
	&= \sum_{h=0}^{\ell}  \sum_{i=s_{h+1}}^{s_h} \frac{1}{i!} R_1^i e^{-R} \P(Y_n \leq n-s_h-1)  \\
	&= \sum_{h=0}^{\ell}  \P(s_{h+1} \leq X_n \leq s_h) \P(Y_n \leq n-s_h-1) 
\end{split}
\end{equation}
where $X_n$ and $Y_n$ are independent Poisson random variables with parameters $R_1$ and $R_2$, respectively.

If $\lambda \rightarrow \infty$ then $R_1 \rightarrow \infty $, since $p$ and $\mu$  are fixed. Note that
\[\P(s_{h+1} \leq X_n \leq s_h)= \P\left( \frac{s_{h+1}-R_1}{\sqrt{R_1}} \leq \frac{X_n-R_1}{\sqrt{R_1}} \leq \frac{s_h-R_1}{\sqrt{R_1}}\right).
\]
Thus, when $\lambda \rightarrow \infty$, by the Central Limit Theorem (Normal approximation to
Poisson) we have 
\[ \frac{X_n-R}{\sqrt{R_1}} \Rightarrow N(0,1).
\]
Since
\begin{equation*}
\begin{split}
&\lim_{\lambda \rightarrow \infty} \frac{s_h-R_1}{\sqrt{R_1}} = \lim_{\lambda \rightarrow \infty} \frac{s-h \epsilon \sqrt{R_1} - R_1} {\sqrt{R_1}} = \lim_{\lambda \rightarrow \infty} \frac{ R_1 - \beta \sqrt{R_1} -h \epsilon \sqrt{R_1} - R_1} {\sqrt{R_1}} = \beta -h \epsilon
\end{split}
\end{equation*}
we obtain:
\begin{equation}
\label{eq_lem3_x}
\begin{split}
&\P(s_{h+1} \leq X_n \leq s_h) = \Phi(\beta -h \epsilon)-\Phi(\beta - (h+1) \epsilon), ~ h=0,..,\ell-1\\
&\P(0 \leq X_n \leq s_{\ell}) = \Phi(\beta -\ell \epsilon).\\
\end{split}
\end{equation}

Similarly, if $\lambda \rightarrow \infty$ then $R_2 \rightarrow \infty $, since $p$,$\delta$ and $\gamma$ are fixed. Note that
\[\P(Y_n \leq n-s_h)= \P\left( \frac{Y_n-R_2} {\sqrt{R_2}} \leq \frac{n-s_h-R_2} {\sqrt{R_2}}\right).
\]
Thus, when $\lambda \rightarrow \infty$, by the Central Limit Theorem (Normal approximation to
Poisson) we have 
\[ \frac{Y_n - R_2} {\sqrt{R_2}} \Rightarrow N(0,1).
\]
 Since
\begin{equation*}
\begin{split}
&\lim_{\lambda \rightarrow \infty} \frac{n-s_h-R_2} {\sqrt{R_2}} 
= \lim_{\lambda \rightarrow \infty} \frac{n-s-h \epsilon \sqrt{R_1} - R_2} {\sqrt{R_2}} 
= \lim_{\lambda \rightarrow \infty} \frac{R_2 + \eta \sqrt{R_2} -h \epsilon \sqrt{R_1} - R_2} {\sqrt{R_2}} \\
&= \eta -h \epsilon \frac{\sqrt{R_1}}{\sqrt{R_2}} = \eta -h \epsilon \sqrt{\frac{\delta}{p\mu}} , 
\end{split}
\end{equation*}
we obtain:
\begin{equation}
\label{eq_lem3_y}
\P(Y_n \leq n-s_h) = \Phi \left(\eta -h \epsilon \sqrt{\frac{\delta}{p\mu}}\right), ~ h=0,..,\ell
\end{equation}

It follows from \eqref{eq_lem3_xi_1}, \eqref{eq_lem3_x}, and \eqref{eq_lem3_y} that
\begin{multline}
\label{eq_lm3_upperbound}
	\lim_{\lambda \rightarrow \infty} A \geq \sum_{h=0}^{\ell-1} ( \Phi(\beta -h \epsilon)-\Phi(\beta - (h+1) \epsilon)) \Phi \left(\eta -h \epsilon \sqrt{\frac{\delta}{p\mu}}\right)
	+ \Phi(\beta -\ell \epsilon) \Phi \left(\eta -\ell \epsilon \sqrt{\frac{\delta}{p\mu}}\right)
\end{multline}
which is the lower Riemann-Stieltjes sum of the integral
\begin{equation}
\label{eq_lm3_RS}
-\int_0^{\infty} \Phi\left( \eta + x \sqrt{\frac{\delta}{p\mu}} \right) d\Phi(\beta-x)
= \int_{-\infty}^{\beta} \Phi\left( \eta + (\beta-t) \sqrt{\frac{\delta}{p\mu}}\right) d\Phi(t)
\end{equation}
corresponding to the partition $\{\beta - h \epsilon\}_{h=0}^{\ell}$ of the semi axis $(-\infty,\beta)$. 
Similarly, let us take the upper estimate for $A$ as the following sum:
\begin{equation}
\label{eq_lem3_xi_2}
\begin{split}
	A \leq A_2 &= \sum_{h=0}^{\ell}  \sum_{i=s_{h+1}}^{s_h} \frac{1}{i!} R_1^i e^{-R_1} \sum_{j=0}^{n-s_{h+1}-1} \frac{1}{j!} R_2^j e^{-R_2}   \\
	&= \sum_{h=0}^{\ell}  \sum_{i=s_{h+1}}^{s_h} \frac{1}{i!} R_1^i e^{-R} \P(Y_n \leq n-s_{h+1}-1)  \\
	&= \sum_{h=0}^{\ell}  \P(s_{h+1} \leq X_n \leq s_h) \P(Y_n \leq n-s_{h+1}-1) 
\end{split}
\end{equation}
where $X_n$ and $Y_n$ are the same random variable as before.
Using the same calculation that were computed for the upper boundary we obtain

\begin{equation}
\label{eq_lm3_lowerbound}
\begin{split}
	\lim_{\lambda \rightarrow \infty} A \leq& \sum_{h=0}^{\ell-1} \left(\Phi(\beta -h \epsilon)-\Phi(\beta - (h+1) \epsilon)\right) \Phi \left(\eta - (h+1) \epsilon \sqrt{\frac{\delta }{p\mu }}\right)+ \Phi(\beta -\ell \epsilon)
\end{split}
\end{equation}
which is the upper Riemann-Stieltjes sum for the integral \eqref{eq_lm3_RS}.
When $\epsilon \rightarrow 0$ the boundaries \eqref{eq_lm3_upperbound} and \eqref{eq_lm3_lowerbound} lead to the following equality

\begin{equation*}
\lim_{\lambda \rightarrow \infty} A = \int_{-\infty}^{\beta} \Phi\left( \eta + (\beta-t) \sqrt{\frac{\delta }{p\mu}} \right) d\Phi(t) = \int_{-\infty}^{\beta} \Phi\left( \frac{ \gamma - t \sqrt{r}}{\sqrt{1-r}} \right) d\Phi(t)
\end{equation*}

This proves Lemma \ref{lem3}.
\endproof

\subsection{Proof of Lemma \ref{lem4}}

\proof
First, we will rewrite Equation \eqref{eq:A}:
\begin{equation*}
	B= e^{-(R_1 + R_2)} \frac{1}{s!}  {R_1}^{s} \sum_{j=0}^{n-s-1} \frac{1}{j!} {R_2}^{j}   \sum_{i=0}^{n-s-j-1} \rho^{i} = e^{-(R_1 + R_2)} \frac{1}{s!}  {R_1}^{s}  \sum_{j=0}^{n-s-1} \frac{1}{j!} {R_2}^{j}   \frac{1-\rho^{n-s-j}}{1-\rho}.
\end{equation*}
When $\beta = 0$, as $\lambda \rightarrow \infty$, by the QED assumption that $\sqrt{s}(1-\rho) \rightarrow \beta$, $\rho = 1$. Therefore,
\[\sum_{i=0}^{n-s-j-1} \rho^{i} = n-s-j .\]
When $\beta \rightarrow 0$, $\rho \rightarrow 1$ but still $\rho \neq 1$, the expression $\frac{1-\rho^{n-s-j}}{1-\rho}$ can be approximated by $\frac{1-\rho^{i}}{1-\rho} \approx i$. Thus, \[\lim_{\rho \rightarrow 1} \sum_{i=0}^{n-s-j-1} \rho^{i} = n-s-j, \] which is the same phrase as when $\rho = 1$.
Hence,  
\begin{equation*}
\begin{split}
	B &\approx e^{-(R_1 + R_2)} \frac{1}{s!}  {R_1}^{s}  \sum_{j=0}^{n-s-1} \frac{1}{j!} {R_2}^{j}  (n-s-j)\\
	&= e^{-(R_1 + R_2)} \frac{1}{s!}  {R_1}^{s} \left((n-s)  \sum_{j=0}^{n-s-1} \frac{1}{j!}  {R_2}^{j}  -  \sum_{j=0}^{n-s-1} \frac{1}{j!} {R_2}^{j}  -  \sum_{j=0}^{n-s-1} \frac{j}{j!} {R_2}^{j} \right) \\
	&= e^{-(R_1 + R_2)} \frac{1}{s!}  {R_1}^{s} \left( (n-s) \sum_{l=0}^{n-s-1} \frac{1}{l!} R_2^{l} - \sum_{j=0}^{n-s-2} \frac{1}{j!} {R_2}^{j} - R_2  \sum_{j=0}^{n-s-2} \frac{1}{j!} {R_2}^{j} \right) \\
	&= e^{-(R_1 + R_2)} \frac{1}{s!}  {R_1}^{s} \left((n-s) \sum_{l=0}^{n-s-1} \frac{1}{l!} R_2^{l}  -  \sum_{l=0}^{n-s-2} \frac{1}{l!} R_2^{l} - R_2 \sum_{l=0}^{n-s-2} \frac{1}{l!} R_2^{l} \right)\\
	&= e^{-(R_1 + R_2)} \frac{1}{s!}  {R_1}^{s} \left((n-s-R_2) \sum_{l=0}^{n-s-1} \frac{1}{l!} R_2^{l} + \frac{R_2^{n-s}}{(n-s-1)!} \right) \\
	&\approx e^{s-R} \frac{1}{\sqrt{2\pi s}}  {\rho}^{s} \left((n-s-R_2) \sum_{l=0}^{n-s-1} \frac{1}{l!} R_2^{l} e^{-R_2} + \frac{R_2^{n-s} e^{-R_2}}{(n-s-1)!} \right) \\
	&= \frac{1}{\sqrt{2\pi s}}  \left((n-s-R_2) \sum_{l=0}^{n-s-1} \frac{1}{l!} R_2^{l} e^{-R_2}	+ \frac{R_2^{n-s} e^{-R_2}}{(n-s-1)!} \right) .
\end{split}
\end{equation*}
As seen in Equation \eqref{lm1_eq1}
\begin{equation*}
(n-s-R_2) \sum_{l=0}^{n-s-1} \frac{1}{l!} R_2^{l} e^{-R_2} \approx \eta \sqrt{R_2} \Phi(\eta).
\end{equation*}
By using Stirling's formula:
\begin{equation*}
\begin{split}
& \frac{R_2^{n-s} e^{-R_2}}{(n-s-1)!} = \frac{(n-s)R_2^{n-s} e^{-R_2}}{(n-s)!} \approx  \frac{(n-s)e^{n-s-R_2}}{\sqrt{2 \pi (n-s)}} \left(\frac{R_2}{n-s}\right)^{n-s} \\
&= \frac{(n-s)e^{n-s-R_2+(n-s)\ln{\left(\frac{R_2}{n-s}\right)}}} {\sqrt{2 \pi (n-s)}} = \sqrt{\frac{n-s}{2 \pi}} e^{(n-s)\left(1-\frac{R_2}{n-s}+ \ln{\left(\frac{R_2}{n-s}\right)}\right)}.
\end{split}
\end{equation*}
By assuming that $\frac{n-s-R_1-R_2}{\sqrt{R_1+R_2}} \rightarrow \eta$  when $\lambda \rightarrow \infty$ 
\begin{equation*}
\begin{split}
&(n-s)\left(1-\frac{R_2}{n-s}+\ln{\left(\frac{R_2}{n-s}\right)}\right) = (n-s)\left(1-\frac{R_2}{n-s}-\left(1-\frac{R_2}{n-s}\right)-\frac{1}{2}\left(1-\frac{R_2}{n-s}\right)^2\right)\\
&=-\frac{n-s}{2}\left(1-\frac{R_2}{n-s}\right)^2 = -\frac{1}{2}\frac{\left(n-s-R_2\right)^2}{n-s} \approx -\frac{1}{2} \frac{\left(\eta \sqrt{R_2}\right)^2}{\eta \sqrt{R_2}+R_2} \approx -\frac{1}{2} \frac{\left(\eta \sqrt{R_2}\right)^2}{\eta \sqrt{R_2}+R_2} \approx -\frac{\eta^2}{2} .
\end{split}
\end{equation*}
Therefore, by  the QED assumption that $\frac{n-s-R_1-R_2}{\sqrt{R_1+R_2}} \rightarrow \eta$ . 
\begin{equation*}
\begin{split}
&\sqrt{\frac{n-s}{2 \pi}} e^{(n-s)\left(1-\frac{R_2}{n-s}+ \ln{\left(\frac{R_2}{n-s}\right)}\right)} \approx \sqrt{\frac{\eta \sqrt{R_2}+R_2}{2 \pi}} e^{-\frac{\eta^2}{2}} = \sqrt{\eta \sqrt{R_2}+R_2} \phi(\eta) \approx \sqrt{R_2} \phi(\eta).
\end{split}
\end{equation*}
Combining the above approximations and the assumption that $\beta=0$ and therefore $s=R_1$ yields
\begin{equation*}
\begin{split}
	B &\approx \frac{1}{\sqrt{2\pi s}}  \left((n-s-R_2) \sum_{l=0}^{n-s-1} \frac{1}{l!} R_2^{l} e^{-R_2} + \frac{R_2^{n-s} e^{-R_2} }{(n-s-1)!} \right)\\
	& \approx \frac{1}{\sqrt{2\pi s}}  \left( \eta \sqrt{R_2} \Phi(\eta) + \sqrt{R_2} \phi(\eta) \right) = \frac{\sqrt{R_2}}{\sqrt{2\pi s}}  \left( \eta  \Phi(\eta) + \phi(\eta) \right) \\
	&= \frac{\sqrt{R_2}}{\sqrt{2\pi R}}  \left( \eta  \Phi(\eta) + \phi(\eta) \right)  = \sqrt{\frac{p \mu}{\delta}} \frac{1}{\sqrt{2 \pi}} \left(\eta \Phi(\eta) +  \phi(\eta) \right) = \sqrt{\frac{1-r}{r}} \frac{1}{\sqrt{2 \pi}} \left(\eta \Phi(\eta) +  \phi(\eta) \right).
	\end{split}
\end{equation*}
This proves Lemma \ref{lem4}.
\endproof

\subsection{Approximation of the Probability of Delay}
The first approximation will be for the measure: the probability of waiting or the probability of delay. 

\begin{theorem} \label{thm1}
Let the variables $\lambda$, $s$ and $n$ tend to $\infty$ simultaneously and satisfy the QED scaling conditions in \eqref{eq:twofoldscaling}. Then if $\beta \neq 0$
\begin{equation*}
\lim_{\lambda \rightarrow \infty} \P(W>0)= \left( 1+  \frac{\beta \int_{-\infty}^{\beta} \Phi\left( \eta + (\beta-t) \sqrt{\xi} \right) d\Phi(t)}{{\phi(\beta)\Phi(\eta)}-{\phi(\sqrt{\eta^2 + \beta^2})} e^{\frac{1}{2} \omega ^2} \Phi(\omega)} \right)^{-1},
\end{equation*}
and if $\beta=0$,
\begin{equation*}
\lim_{\lambda \rightarrow \infty} \P(W>0)= \left( 1+\frac{\int_{-\infty}^{0} \Phi\left( \frac{\gamma -t \sqrt{r}}{\sqrt{1-r}} \right) d\Phi(t)}{\sqrt{\frac{1-r}{r}} \frac{1}{\sqrt{2\pi}} \left(\eta \Phi(\eta) + \phi(\eta) \right)} \right)^{-1},
\end{equation*}
where $\xi=\frac{R_1}{R_2}=\frac{\delta }{p\mu}$, $\omega = \eta - \beta \sqrt{\xi^{-1}}$.
\end{theorem}

\proof{Proof: }
By the Arrival theorem for closed networks,
\begin{equation*}
\begin{split}
\P_n(W > 0) &= \P_{n-1}(Q_1(\infty) \geq s) = \sum_{m=s}^{n-1} \sum_{i=s}^{m} \pi^b_{n-1}(i,m-i)\\
&= \pi_0^{n-1} \sum_{m=s}^{n-1} \sum_{i=s}^{m}  \frac{1}{s!s^{i-s}}  R_1^i  \frac{1}{(m-i)!}  R_2 ^{m-i} , 
\end{split}
\end{equation*}
where
\begin{equation*}
\pi_0^{n-1} =\left( \sum_{l=0}^{n-1}  \frac{1}{l!} \left( R_1 + R_2  \right)^{l}  + \sum_{m=s}^{n-1} \sum_{i={s}}^m \left( \frac{1}{s!s^{i-s}} - \frac{1}{i!} \right) \frac{1}{(m-i)!}  R_1^i R_2^{m-i}  \right)^{-1} .
\end{equation*}
Thus,
\[\P_n(W > 0) =\left( 1+\frac{A}{B} \right)^{-1} ,\]
where
\begin{equation*}
\begin{split}
A&= \sum_{l=0}^{n-1}  \frac{1}{l!} \left( R_1 + R_2  \right)^{l} e^{-\left(R_1 + R_2  \right)}- \sum_{m=s}^{n-1} \sum_{i=s}^m  \frac{1}{i!(m-i)!} R_1^i  R_2^{m-i} e^{-\left(R_1 + R_2  \right)} \\
&= \sum_{\substack{i,j| i \leq s, \\i+j \leq n-1 }} \frac{1}{i!j!} R_1^i  R_2^j e^{-\left(R_1 + R_2  \right)} ,
\end{split}
\end{equation*}
\begin{equation}
\label{eq:B}
\begin{split}
	B&= \sum_{m=s}^{n-1} \sum_{i=s}^m  \frac{1}{s!s^{i-s}}  R_1^i \frac{1}{(m-i)!} R_2^{m-i}  e^{-\left(R_1 + R_2  \right)} = \sum_{j=0}^{n-s-1}  \sum_{i=s}^{n-j-1}  \frac{1}{s!s^{i-s}}   R_1^{i} \frac{1}{j!} R_2^{j} e^{-\left(R_1 + R_2  \right)} \\
	&= \sum_{j=0}^{n-s-1} \sum_{i=0}^{n-s-j-1}   \frac{1}{s!s^{i}}  \frac{1}{j!} R_1^{i+s} R_2^{j} e^{-\left(R_1 + R_2  \right)} = \frac{1}{s!} R_1^{s} e^{-\left(R_1 + R_2  \right)} \sum_{j=0}^{n-s-1} \frac{1}{j!} R_2^{j}  \sum_{i=0}^{n-s-j-1}  \rho^{i} ,
\end{split}
\end{equation}	
and $\rho=\frac{\lambda}{(1-p)s\mu} = \frac{R_1}{s}$.

Then under the QED assumption that
$\sqrt{s} (1-\rho) \rightarrow \beta, ~ -\infty < \beta < \infty$, if $ \beta \neq 0$ as $\lambda \rightarrow \infty$, we can rewrite the right-hand side in the following way:
\begin{equation}
\label{eq:B1B2}
\begin{split}
	B&=\frac{1}{s!} R_1^{s} e^{-\left(R_1 + R_2  \right)} \sum_{j=0}^{n-s-1} \frac{1}{j!} R_2^{j}  \frac{1-\rho^{n-s-j}}{1-\rho} \\
	&=\frac{1}{s!} R_1^{s} e^{-\left(R_1 + R_2  \right)}\frac{1}{1-\rho} \sum_{j=0}^{n-s-1} \frac{1}{j!} R_2^{j} - \frac{1}{s!} R_1^{s} e^{-\left(R_1 + R_2  \right)}\frac{\rho^{n-s}}{1-\rho} \sum_{j=0}^{n-s-1} \frac{1}{j!} \left(\frac{R_2}{ \rho} \right)^{j}  = B_1 -B_2. 
\end{split}
\end{equation}

Applying Lemmas \ref{lem1},\ref{lem2}, and \ref{lem3} if $\beta \neq 0$ we get 
\begin{equation*}
\lim_{\lambda \rightarrow \infty} \P(W>0)= \left( 1+ \frac{\beta \int_{-\infty}^{\beta} \Phi\left( \eta + (\beta-t) \sqrt{\frac{\delta }{p\mu}} \right) d\Phi(t)}{\phi(\beta)\Phi(\eta)-\phi(\sqrt{\eta^2 + \beta^2}) e^{\frac{1}{2} \omega ^2} \Phi(\omega)} \right)^{-1}.
\end{equation*}
%
%
%
%
Applying Lemma \ref{lem3} and \ref{lem4} when $\beta = 0$ we get
\begin{equation*} 
\lim_{\lambda \rightarrow \infty} \P(W>0)= \left( 1+\frac{\int_{-\infty}^{0} \Phi\left( \frac{\gamma -t \sqrt{r}}{\sqrt{1-r}} \right) d\Phi(t)}{\sqrt{\frac{1-r}{r}} \frac{1}{\sqrt{2\pi}} \left(\eta \Phi(\eta) + \phi(\eta) \right)} \right)^{-1} .
\end{equation*}
This proves Theorem \ref{thm1}.
\endproof

\subsection{Approximation of the Expected Waiting Time}

In this appendix we will prove the approximation for the expected waiting time, stated in Section \ref{sec:QED_limit_block}. The exact measure was defined in Section \ref{sec:performance_metrics}. For convenience we state the theorem here again. 
The first theorem gives the approximation for the case where $\beta \neq 0$.
\begin{theorem} \label{thm3}
Let the variables $\lambda$, $s$ and $n$ tend to $\infty$ simultaneously and satisfy the QED scaling conditions in \eqref{eq:twofoldscaling} with $\beta \neq 0$.
Then 
\begin{equation*}
\lim_{\lambda \rightarrow \infty} \sqrt{s} \E[W]= \frac{1}{\mu} \, \frac{\left(\frac{1}{\xi}-\frac{\eta}{\beta\sqrt{\xi}}-\frac{1}{\beta^2} \right) \phi(\sqrt{\beta^2+\eta^2})\, {\rm e}^{\tfrac{1}{2}\omega^2}\,\Phi(\omega)+ \frac{\phi(\beta)\Phi(\eta)}{\beta^2} - \frac{\phi(\beta)\phi(\eta)}{\sqrt{\xi}\b}}{
\int_{-\infty}^\beta \Phi(\eta+(\beta-t)\sqrt{\xi})\, d\Phi(t) + \frac{\phi(\beta)\Phi(\eta)}{\beta} -\frac{\phi(\sqrt{\beta^2+\eta^2})}{\beta}\, {\rm e}^{\tfrac{1}{2}\omega^2}\,\Phi(\omega)},
\end{equation*}
where $\xi=\frac{R_1}{R_2}=\frac{\delta }{p\mu}$, $\omega = \eta - \beta \sqrt{\xi^{-1}}$.
\end{theorem}

\proof{Proof:}
It follows from \eqref{eq:EW_exact} that the expectation of the waiting time is given by
\begin{equation}
\begin{split}
	\E[W]&= \int_0^{\infty} p_n(s; t) dt = \frac{1}{\mu s} \sum_{m=s}^{n-1}\sum^m_{i=s} \pi_{n-1}(i,m-i) (i-s+1)\\
	  &= \frac{1}{\mu s} \sum_{m=s}^{n-1} \sum^m_{i=s} \pi_{n-1}(i,m-i) (i-s) + \frac{1}{\mu s} \sum_{m=s}^{n-1} \sum^m_{i=s}  \pi_{n-1}(i,m-i)\\
	  &= \frac{1}{\mu s} \sum_{m=s}^{n-1} \sum^m_{i=s} \pi_{n-1}(i,m-i) (i-s) + \frac{1}{\mu s} \P(W>0) =C+D,
\end{split}
\end{equation}
where $D$ is given by
\begin{equation*}
D = \frac{1}{\mu s} \P(W>0) = \frac{1}{\mu s} \frac{B}{A+B}
\end{equation*}
and $A$ and $B$ were defined in \eqref{eq:A} and \eqref{eq:B}, respectively, 
and $C$ is given by,
\begin{equation*}
\begin{split}
	C &= \frac{1}{\mu s} \sum_{m=s}^{n-1}\sum^m_{i=s} \pi_{n-1}(i,m-i) \,(i-s)\\
	  &= \frac{1}{\mu s}\, \pi_0\, \sum_{m=s}^{n-1}\sum^m_{i=s} \frac{R_1^i}{s!s^{i-s}}\, \frac{R_2^{m-i}}{(m-i)!} \,(i-s) = \frac{1}{\mu s}\,\frac{G e^{-(R_1+R_2)}}{A+B}.
\end{split}
\end{equation*}

We will rewrite $G$ in the following way:
\begin{align*}
	G &= \sum_{m=s}^{n-1}\sum^m_{i=s} \frac{R_1^i}{s!s^{i-s}}\, \frac{R_2^{m-i}}{(m-i)!} \,(i-s) = \sum_{j=0}^{n-s-1} \sum_{i=s}^{n-j-1}\frac{R_1^i}{s!s^{i-s}}\, \frac{R_2^j}{j!} \,(i-s)\\
	&= \sum_{j=0}^{n-s-1} \sum_{i=0}^{n-s-j-1}\frac{R_1^{i+s}}{s!s^{i}}\, \frac{R_2^j}{j!} \,i
= \frac{R_1^s}{s!} \sum_{j=0}^{n-s-1} \sum_{i=0}^{n-s-j-1} i\,\rho^i\,\frac{R_2^j}{j!}
\end{align*}
Using the formula
\begin{equation}
\sum_{l=0}^M l\,\rho^l = \rho\Big( \sum_{l=0}^M \rho^l \Big)' = \rho\,\Big( \frac{1-\rho^{M+1}}{1-\rho}\Big)' = {-}M\frac{\rho^{M+1}}{1-\rho} + \rho\,\frac{1-\rho^M}{(1-\rho)^2},
\end{equation}
we can rewrite $G$ as a sum $G = G_1+G_2$, where
\begin{equation*}
G_1 = {-} \frac{R_1^s}{s!}\, \sum_{j=0}^{n-s-1} (n-s-j-1)\,\frac{\rho^{n-s-j}}{1-\rho}\, \frac{R_2^j}{j!},
\end{equation*}
and 
\begin{equation*}
G_2 = \frac{R_1^s}{s!}\sum_{j=0}^{n-s-1} \rho \, \frac{1-\rho^{n-s-j-1}}{(1-\rho)^2}\, \frac{R_2^j}{j!}.
\end{equation*}
Therefore,
\begin{align*}
G_1 &= {-} \frac{R_1^s}{s!}\, \sum_{j=0}^{n-s-1} (n-s-j-1)\,\frac{\rho^{n-s-j}}{1-\rho}\, \frac{R_2^j}{j!} \\ 
&= {-}(n-s-1)\,\frac{R_1^s}{s!} \sum_{j=0}^{n-s-1}\frac{\rho^{n-s-j}}{1-\rho}\, \frac{R_2^j}{j!} 
+ \frac{R_1^s}{s!} \sum_{j=0}^{n-s-1}j\,\frac{\rho^{n-s-j}}{1-\rho}\, \frac{R_2^j}{j!} \\
&= {-}(n-s-1)\,\frac{\rho^{n-s}}{1-\rho}\frac{R_1^s}{s!} \sum_{j=0}^{n-s-1}\, \frac{(R_2/\rho)^j}{j!} 
+ \frac{R_1^s}{s!} \,\frac{\rho^{n-s}}{1-\rho} \sum_{j=0}^{n-s-1} j\,\frac{(R_2/\rho)^j}{j!} \\
&= {-}(n-s-1)\,\frac{\rho^{n-s}}{1-\rho}\frac{R_1^s}{s!} \sum_{j=0}^{n-s-1}\, \frac{(R_2/\rho)^j}{j!} 
+ \frac{R_1^s}{s!} \, \frac{R_2}{\rho} \,\frac{\rho^{n-s}}{1-\rho} \sum_{j=0}^{n-s-2} \frac{(R_2/\rho)^j}{j!} \\
&= {-}\left(n-s-R_2/\rho - 1\right)\,\frac{\rho^{n-s}}{1-\rho}\frac{R_1^s}{s!} \sum_{j=0}^{n-s-1}\, \frac{(R_2/\rho)^j}{j!} 
- \frac{R_1^s}{s!} \,\frac{1}{1-\rho}\frac{R_2^{n-s}}{(n-s-1)!} \\
&= {-}\left(n-s-R_2/\rho - 1\right)\, {\rm e}^{R_1+R_2}\,B_2 
- \frac{1}{1-\rho}\,\frac{R_1^s}{s!}\,\frac{R_2^{n-s}}{(n-s-1)!},
\end{align*}
where $B_2$ was defined in Lemma \ref{lem2}. 
For $G_2$ we have, 
\begin{align*}
G_2 &= \frac{R_1^s}{s!}\,\frac{\rho}{(1-\rho)^2}\sum_{j=0}^{n-s-1} \left( 1- \rho^{n-s-j-1}\right) \frac{R_2^j}{j!}\\
&=  \frac{R_1^s}{s!}\,\frac{\rho}{(1-\rho)^2}\sum_{j=0}^{n-s-1} \frac{R_2^j}{j!}
 - \frac{R_1^s}{s!}\,\frac{\rho^{n-s}}{(1-\rho)^2} \sum_{j=0}^{n-s-1} \frac{(R_2/\rho)^j}{j!}\\
&= {\rm e}^{R_1+R_2}\,\frac{1}{1-\rho}\,(\rho B_1 - B_2 ),
\end{align*}
with $B_1$ as in Lemma \ref{lem1}. 
Note that
\begin{align*}
&n-s-\frac{R_2}{\rho}-1 = R_2 - \eta \sqrt{R_2} - \frac{R_2}{\rho} - 1 = R_2\left(1-1/\rho\right) + \eta\sqrt{R_2} - 1\\
&= R_2\left( \frac{R_1-s}{\sqrt{R_1}}\right) + \eta\sqrt{R_2} - 1 = R_2 \cdot {-} \frac{\beta\sqrt{R_1}}{R_1} + \eta\sqrt{R_2} -1 \approx \sqrt{R_1}\left( \eta\sqrt{\xi} - \beta\xi\right) \\
&\approx \sqrt{s}\left( \eta/\sqrt{\xi} - \beta/\xi\right),
\end{align*}
for $s$ large, where $\xi = R_1 / R_2 = \delta/p\mu$.
Furthermore, 
\begin{align*}
\frac{1}{1-\rho}\,\frac{R_1^s}{s!}\,\frac{R_2^{n-s}}{(n-s-1)!}
&= \frac{n-s}{1-\rho} \, \frac{R_1^s}{s!} \, \frac{R_2^{n-s}}{(n-s)!}\\
&= \frac{n-s}{1-\rho} \, {\rm e}^{R_1+R_2} \, \P\left({\rm Pois}(R_1) = s \right)\,\P\left({\rm Pois}(R_2) = n-s\right) \\
& = \frac{R_2 + \eta\sqrt{R_2}}{\beta/\sqrt{R_1}} {\rm e}^{R_1+R_2} \P\left({\rm Pois}(R_1) = s \right)\,\P\left({\rm Pois}(R_2) = n-s\right) \\
&= \sqrt{R_2} \, {\rm e}^{R_1+R_2}\, \frac{1}{\beta}\, \sqrt{R_1} \P\left({\rm Pois}(R_1) = s \right)\cdot
\sqrt{R_2} \P\left({\rm Pois}(R_2) = n-s\right))\\
& \approx \sqrt{R_1}\,{\rm e}^{R_1+R_2} \frac{1}{\sqrt{\xi}\beta}\,\phi(\beta)\phi(\eta),
\end{align*}
where we used that $\sqrt{R}\, \P({\rm Pois}(R) = R+x\sqrt{R}) \to \phi(x)$ as $R\to\infty$.
Hence, we have for $G_1$
\begin{align*}
\frac{1}{\sqrt{s}}{\rm e}^{-(R_1+R_2)} \, G_1 \to  - \left( \eta/\sqrt{\xi} - \beta/\xi\right)\cdot \frac{\phi(\sqrt{\beta^2+\eta^2})}{\beta}\, {\rm e}^{\tfrac{1}{2}\omega^2}\,\Phi(\omega)-\frac{1}{\sqrt{\xi}\beta}\,\phi(\beta)\phi(\eta)
\end{align*}
and for $G_2$
\begin{align*}
\frac{1}{\sqrt{s}}{\rm e}^{-(R_1+R_2)} \, G_2 &= \frac{1}{\sqrt{s}} \frac{1}{1-\rho}\left(\rho B_1-B_2\right)\approx \frac{1}{\beta}\,\left(\rho B_1 - B_2 \right) \\
&\to \frac{\phi(\beta)\Phi(\eta)}{\beta^2} -\frac{\phi(\sqrt{\beta^2+\eta^2})}{\beta^2}\, {\rm e}^{\tfrac{1}{2}\omega^2}\,\Phi(\omega),
\end{align*}
as $s\to\infty$.
In total, this gives
\begin{equation}
\frac{1}{\sqrt{s}}{\rm e}^{-(R_1+R_2)} \,G \to
\left(\frac{1}{\xi}-\frac{\eta}{\beta\sqrt{\xi}}-\frac{1}{\beta^2} \right) \phi(\sqrt{\beta^2+\eta^2})\, {\rm e}^{\tfrac{1}{2}\omega^2}\,\Phi(\omega)+ \frac{\phi(\beta)\Phi(\eta)}{\beta^2} - \frac{\phi(\beta)\phi(\eta)}{\sqrt{\xi}\b}.
\end{equation}
We can conclude, 
\begin{align*}
\sqrt{R_1}\,\E[W] &= \frac{\sqrt{R_1}}{\mu s} \frac{G e^{-(R_1+R_2)}}{A+B} + \frac{\sqrt{s}}{\mu s}\,\P( W>0)  = \frac{\sqrt{s}}{\mu s}\, \frac{G e^{-(R_1+R_2)}}{A+B} + \frac{1}{\mu \sqrt{s}}\,\P(W>0) \\
&\to \frac{1}{\mu} \, \frac{\left(\frac{1}{\xi}-\frac{\eta}{\beta\sqrt{\xi}}-\frac{1}{\beta^2} \right) \phi(\sqrt{\beta^2+\eta^2})\, {\rm e}^{\tfrac{1}{2}\omega^2}\,\Phi(\omega)+ \frac{\phi(\beta)\Phi(\eta)}{\beta^2} - \frac{\phi(\beta)\phi(\eta)}{\sqrt{\xi}\b}}{
\int_{-\infty}^\beta \Phi(\eta+(\beta-t)\sqrt{\xi})\, d\Phi(t) + \frac{\phi(\beta)\Phi(\eta)}{\beta} -\frac{\phi(\sqrt{\beta^2+\eta^2})}{\beta}\, {\rm e}^{\tfrac{1}{2}\omega^2}\,\Phi(\omega)},
\end{align*}
as $s\to\infty$.
\endproof

The second theorem gives the approximation for the case where $\beta=0$.
\begin{theorem} \label{thm4}
Let the variables $\lambda$, $s$ and $n$ tend to $\infty$ simultaneously and satisfy the QED scaling conditions in \eqref{eq:twofoldscaling} with $\beta = 0$. Then
\begin{equation*}
\lim_{\lambda \rightarrow \infty} \sqrt{s}\E[W]= \frac{1}{2\mu}  \frac{\xi^{-1} \left((\eta^2+1) \Phi(\eta) + \eta \phi(\eta) \right) }{\sqrt{2\pi} \int_{-\infty}^{0} \Phi\left( \eta -t \sqrt{\xi} \right) d\Phi(t) + \sqrt{\xi^{-1}}  \left(\eta \Phi(\eta) + \phi(\eta) \right)},
\end{equation*}
where $\xi=\frac{R_1}{R_2}=\frac{\delta }{p\mu }$, $\omega = \eta - \beta \sqrt{\xi^{-1}}$.
\end{theorem}

\proof{Proof:} As before
\begin{equation}
\E[W] = \frac{1}{\mu s} \frac{G e^{-(R_1+R_2)} +B}{A+B}.
\end{equation}
Since $\beta = 0$, we have $\rho = 1$ so that
\begin{align*}
G &= \frac{R_1^s}{s!} \, \sum_{j=0}^{n-s-1} \sum_{i=0}^{n-s-j-1} i\, \frac{R_2^j}{j!} = \frac{1}{2}\frac{R_1^s}{s!} \, \sum_{j=0}^{n-s-1} (n-s-j)(n-s-j-1)\, \frac{R_2^j}{j!} \\
&= \frac{1}{2}\,\frac{R_1^s}{s!}\,(n-s-1)\sum_{j=0}^{n-s-1} \frac{R_2^j}{j!}\,(n-s-j) - \frac{1}{2}\,\frac{R_1^s}{s!}\sum_{j=0}^{n-s-1} j\,\frac{R_2^j}{j!}\,(n-s-j)\\
&= \frac{1}{2}\,\frac{R_1^s}{s!}\,(n-s-1)\sum_{j=0}^{n-s-1} \frac{R_2^j}{j!}\,(n-s-j) - \frac{1}{2}\,\frac{R_1^s}{s!}\,R_2\,\sum_{j=0}^{n-s-2}\frac{R_2^j}{j!}\,(n-s-j-1)\\
&= \frac{1}{2}\,\frac{R_1^s}{s!}\,(n-s-1)\sum_{j=0}^{n-s-1} \frac{R_2^j}{j!}\,(n-s-j) - \frac{1}{2}\,\frac{R_1^s}{s!}\,R_2\,\sum_{j=0}^{n-s-1}\frac{R_2^j}{j!}\,(n-s-j-1)\\
&= \frac{1}{2}\,\frac{R_1^s}{s!}\,(n-s-1)\sum_{j=0}^{n-s-1} \frac{R_2^j}{j!}\,(n-s-j) - \frac{1}{2}\,\frac{R_1^s}{s!}\,R_2\,\sum_{j=0}^{n-s-1}\frac{R_2^j}{j!}\,(n-s-j)+\frac{1}{2}\frac{R_1^s}{s!}\,R_2\,\sum_{j=0}^{n-s-1} \frac{R_2^j}{j!}\\
&= \frac{1}{2}\,\frac{R_1^s}{s!}\,(n-s-R_2-1)\sum_{j=0}^{n-s-1} \frac{R_2^j}{j!}\,(n-s-j)+\frac{1}{2}\frac{R_1^s}{s!}\,R_2\,\sum_{j=0}^{n-s-1} \frac{R_2^j}{j!}.
\end{align*}
Here, 
\begin{align}
{\rm e}^{-R_2}\sum_{j=0}^{n-s-1} \frac{R_2^j}{j!}\,(n-s-j) 
&= {\rm e}^{-R_2}\,(n-s)\sum_{j=0}^{n-s-1}\frac{R_2^j}{j!} - {\rm e}^{-R_2}R_2\,\sum_{j=0}^{n-s-2}\frac{R_2^j}{j!}\nonumber\\
&= {\rm e}^{-R_2}\,(n-s-R_2)\sum_{j=0}^{n-s-1}\frac{R_2^j}{j!} + {\rm e}^{-R_2}(n-s)\,\frac{R_2^{n-s}}{(n-s)!}\nonumber\\
&= \eta\sqrt{R_2}\,\P(\Poi(R_2)\leq n-s-1) + (n-s-1)\P(\Poi(R_2) = n-s)\nonumber\\
&= \eta\sqrt{R_2}\,\P(\Poi(R_2)\leq n-s) + (R_2-1)\P(\Poi(R_2) = n-s)\nonumber\\
&\approx \sqrt{R_2}\left(\eta\,\Phi(\eta) + \phi(\eta)\right) = \sqrt{\xi^{-1}}\sqrt{R_1}\left(\eta\,\Phi(\eta) + \phi(\eta)\right).
\label{eq:x1}
\end{align}
Furthermore 
\begin{equation}
n-s-1-R_2 = \eta\sqrt{R_2}-1 \approx \eta\sqrt{\xi^{-1}}\sqrt{R_1}
\label{eq:x2}
\end{equation}
and 
\begin{align}
{\rm e}^{-R_1-R_2} \,\frac{R_1^s}{s!}\,R_2\,\sum_{j=0}^{n-s-1} \frac{R_2^j}{j!} 
&= \xi^{-1}\,\sqrt{R_1}\, \sqrt{R_1}\P(\Poi(R_1)=s)\,\P(\Poi(R_2\leq n-s-1) \nonumber\\
&\approx \sqrt{R_1}\xi^{-1}\, \phi(0) \Phi(\eta).
\label{eq:x3}
\end{align}
Combining \eqref{eq:x1}-\eqref{eq:x3}, we find
\begin{align*}
{\rm e}^{-R_1-R_2} G 
&\approx
\frac{1}{2}\,\xi^{-1}\,\eta\sqrt{R_1}\,\P(\Poi(R_1)=s)\,\left(\eta\Phi(\eta) +\phi(\eta)\right) + \frac{1}{2}\sqrt{R_1}\phi(0)\Phi(\eta) \\
&\approx \frac{1}{2}\,\xi^{-1}\eta\sqrt{R_1}\left(\eta\,\Phi(\eta) + \phi(\eta)\right) + \frac{1}{2}\,\xi^{-1}\sqrt{R_1}\,\phi(0)\Phi(\eta) \\
&= \frac{1}{2} \, \sqrt{R_1}\,\xi^{-1}\phi(0) \left[ (\eta^2+1)\Phi(\eta) + \eta\,\phi(\eta)\right],
\end{align*}
for $R_1$ large.
Hence, we can conclude 
\begin{align*}
\sqrt{R_1}\,\E[W] &= \frac{\sqrt{R_1}}{\mu s}\, \frac{G  e^{-R_1-R_2} }{A+B} + \frac{\sqrt{R_1}}{\mu s} \P(W>0) \approx \frac{1}{\sqrt{R_1}\mu} \,\frac{{\rm e}^{-R_1-R_2}G}{A+B} + \frac{1}{\sqrt{R_1}\mu}\, \P(W>0)\\
&\to  \frac{1}{2\mu} \frac{ \xi^{-1} \left[ (\eta^2+1)\Phi(\eta) + \eta\,\phi(\eta)\right]}
{\sqrt{2 \pi}\int_{-\infty}^0 \Phi(\eta-t\sqrt{\xi}) d\Phi(t) + \sqrt{\xi^{-1}}(\eta \Phi(\eta)+\phi(\eta))},
\end{align*}
where we used that $\phi(0) = 1/\sqrt{2\pi}$, and $\eta = \frac{\g - \b\sqrt{r}}{\sqrt{1-r}}$.

\subsection{Approximation of the Blocking Probability}
\label{Pb_theorem}

We prove here the approximations given in Theorem \ref{thm:limits_YT} of the probability of blocking. For convenience we repeat here the relevant theorem.
\begin{theorem} \label{thm5}
Let the variables $\lambda$, $s$ and $n$ tend to $\infty$ simultaneously and satisfy the QED scaling conditions in \eqref{eq:twofoldscaling} with $\beta \neq 0$. 
Then
\begin{equation}
\lim_{\lambda \rightarrow \infty} \sqrt{s} \P(block)= \frac{\sqrt{r}\phi(\gamma)\Phi(-\omega \sqrt{r})+\phi(\sqrt{\eta^2+\beta^2})e^{\frac{\omega^2}{2}}\Phi(\omega)}{\int_{-\infty}^{\beta} \Phi\left( \eta + (\beta-t) \sqrt{\xi} \right) d\Phi(t) + \frac{\phi(\beta)\Phi(\eta)}{\beta }-\frac{\phi(\sqrt{\eta^2 + \beta^2})}{\beta} e^{\frac{1}{2} \omega ^2} \Phi(\omega)}
\end{equation}
where $\eta = \frac{\gamma-\beta \sqrt{r}}{\sqrt{1-r}}$, and $\omega = \eta - \frac{\beta \sqrt{1-r}}{\sqrt{r}} = \frac{\gamma - \beta / \sqrt{r}}{\sqrt{1-r}} $.
\end{theorem}

\proof{Proof:}
It follows from \eqref{eq:delay_probability} that the probability of blocking is given by

\begin{equation*}
\begin{split}
\P_n &= \pi_0 \left(\sum_{i=0}^s \frac{1}{i!} R_1^i \frac{1}{(n-i)!} R_2^{n-i} + \sum_{i=s+1}^n \frac{1}{s!s^{i-s}} R_1^i \frac{1}{(n-i)!} R_2^{n-i} \right) 
	%
	= \frac{\delta_1+\delta_2}{\tilde{\xi}+\tilde{\zeta_1}-\tilde{\zeta_2}} ,
\end{split}
\end{equation*}
where
\begin{equation*}
\begin{split}
	{\delta_1}&=  \sum_{i=0}^s  \frac{1}{i!} R_1^i \frac{1}{(n-i)!} R_2^{n-i}  e^{-\left(R_1 + R_2 \right)}\\
	{\delta_2}&=  \sum_{i=s+1}^n \frac{1}{s!s^{i-s}} R_1^i \frac{1}{(n-i)!} R_2^{n-i} e^{-\left(R_1 + R_2  \right)}\\
	\tilde{A} &=  \sum_{\substack{i,j| i \leq s, \\i+j \leq n }} \frac{1}{i!j!} R_1^i  R_2^j   e^{-\left(R_1 + R_2  \right)}\\
	\tilde{B_1}&=  \frac{1}{s!} R_1^{s} \frac{1}{1-\rho} \sum_{l=0}^{n-s} \frac{1}{l!} R_2^{l} e^{-\left(R_1 + R_2  \right)}\\
	\tilde{B_2}&=  \frac{1}{s!} R_1^{s} \frac{\rho^{n-s+1}}{1-\rho} \sum_{l=0}^{n-s} \frac{1}{l!} \left(\frac{R_2}{\rho}  \right)^{l} e^{-\left(R_1 + R_2  \right)}.
\end{split}
\end{equation*}

Note that by Lemmas \ref{lem1},\ref{lem2} and \ref{lem3}
\begin{equation*}
\begin{split}
	&\lim_{\lambda \rightarrow \infty}\tilde{B_1} = \lim_{\lambda \rightarrow \infty} B_1 = \frac{\phi(\beta)\Phi(\eta)}{\beta},\\
	&\lim_{\lambda \rightarrow \infty}\tilde{B_2} = \lim_{\lambda \rightarrow \infty} B_2 = \frac{\phi(\sqrt{\eta^2+\beta^2})}{\beta}e^{\frac{1}{2}\eta_1^2}\Phi(\eta_1),\\
	&\lim_{\lambda \rightarrow \infty}\tilde{A} = \lim_{\lambda \rightarrow \infty} A = \int_{-\infty}^{\beta} \Phi\left( \eta + (\beta-t) \sqrt{\frac{\delta }{p\mu }} \right) d\Phi(t),	
\end{split}
\end{equation*}

\begin{equation*}
\begin{split}
	{\delta_1}&=  \sum_{i=0}^s \frac{1}{i!} R_1^i \frac{1}{(n-i)!} R_2^{n-i}  e^{-\left(R_1 + R_2  \right)} = \frac{1}{n!} e^{-\left(R_1 + R_2 \right)}\sum_{i=0}^s  \frac{n!}{i!(n-i)!} R_1^i R_2^{n-i} \\
	&= \frac{1}{n!} e^{-\left(R_1 + R_2 \right)}\left( R_1 +  R_2  \right)^{n} \sum_{i=0}^s  \frac{n!}{i!(n-i)!} \left( \frac{R_1}{ R_1 +  R_2 } \right)^i \left( \frac{R_2}{ R_1 +  R_2 } \right)^{n-i} \\
	&= \frac{1}{n!} e^{-\left(R_1 + R_2  \right)}\left( R_1 +  R_2  \right)^{n} \sum_{i=0}^s  \P(X_{\lambda}=(i,n-i))\\
	&= \P(Y_{\lambda}=n) \sum_{i=0}^s \P(X_{\lambda}=(i,n-i)) = \P(Y_{\lambda}=n) \P(X^1_{\lambda}\leq s)
\end{split}
\end{equation*}

where $X_{\lambda}$ is a random variable with Multinomial distribution with parameters $(n,p_i,p_j)$,  $p_i=\frac{R_1}{ R_1 +  R_2 }$, $p_j = \frac{R_2}{ R_1 +  R_2}$,  $Y_{\lambda}$ is a random variable with Poisson distribution with parameter $R_1 + R_2$, and
$X^1_{\lambda}$ is a random variable with Binomial distribution with parameters $(n,p_i)$.
By the CLT and the use of \ref{eq:lem2_app_n-s}
\begin{equation}
\label{eq:x_small_s}
\begin{split}
	\P(X^1_{\lambda}\leq s) &= \Phi \left( \frac{s-np_i}{\sqrt{np_i(1-p_i)}}\right)
	=\Phi \left( \frac{s-n \frac{R_1}{R_1+R_2}}{\sqrt{n\frac{R_1}{R_1+R_2}(1-\frac{R_1}{R_1+R_2})}}\right) \\
	&=\Phi \left( \frac{s-n \frac{R_1}{R_1+R_2}} {\sqrt{n\frac{R_1}{R_1+R_2}\frac{R_2}{R_1+R_2}}}\right) 
	=\Phi \left( \frac{s (R_1+R_2) - n R_1 } {\sqrt{n R_1 R_2}}\right)\\
	&=\Phi \left( \sqrt{\frac{R_1}{R_2}} \frac{s \frac{R_1+R_2}{R_1} - n } {\sqrt{n }}\right) 
	= \Phi \left( \sqrt{\frac{R_1}{R_2}} \frac{s (1+\frac{R_2}{R_1}) - n } {\sqrt{n }}\right)\\
	&=\Phi \left( \sqrt{\frac{R_1}{R_2}} \frac{s + \frac{R_2}{\rho} - n } {\sqrt{n }}\right) 
	= \Phi \left( - \sqrt{\frac{R_1}{R_2}} \sqrt{\frac{R_2}{n\rho}} \cdot \frac{n - s - \frac{R_2}{\rho} } {\sqrt{\frac{R_2}{\rho}}}\right)	\\
	&= \Phi \left( - \sqrt{\frac{s}{n}} \cdot \frac{n - s - \frac{R_2}{\rho} } {\sqrt{\frac{R_2}{\rho}}}\right)
	\approx \Phi \left( - \sqrt{\frac{R_1}{R_1+R_2}} \cdot \left(\eta-\beta \frac{R_2}{R_1} \right) \right)\\
	&= \Phi \left(    \frac{\beta \frac{R_2}{R_1}-\eta}{\sqrt{1+\frac{R_2}{R_1}}} \right)
	=  \Phi \left( \frac{ \beta \sqrt{\frac{p\mu}{\delta}} -  \eta } {\sqrt{  1+  \frac{p\mu}{\delta} }}\right). 
\end{split}
\end{equation}

By the normal approximation of the Poisson distribution:
\begin{equation}
\label{eq:y_equal_n}
\begin{split}
	\P(Y_{\lambda} = n) & \approx \frac{1}{\sqrt{R_1 + R_2}} \phi \left( \frac{n- \left( R_1 + R_2\right)}{\sqrt{R_1 + R_2}}\right) \approx \frac{1}{\sqrt{s}\sqrt{1+  \frac{p\mu}{\delta}}} \phi \left( \frac{\eta \sqrt{ \frac{p\mu}{\delta}} + \beta }{\sqrt{1+ \frac{p\mu}{\delta}}}\right). 
\end{split}
\end{equation}

We based on the following equivalences (as $\lambda$ tends to $\infty$) to develop Equations \ref{eq:x_small_s} and \ref{eq:y_equal_n}:
\begin{equation*}
\begin{split}
	&R_1 + R_2 = \frac{\lambda}{(1-p)\mu} + \frac{p\lambda}{(1-p)\delta} 
		\approx s +  (1-p)s \mu\left(\frac{p}{(1-p)\delta}\right) = s \left( 1+  \frac{ \mu p }{\delta }\right) = s \left( 1+  \frac{ R_2}{R_1}\right); \\
	& n-\left(R_1 + R_2\right)
	\approx s +  R_2  + \eta \sqrt{ R_2 } -\left(R_1 + R_2 \right) = s  + \eta \sqrt{ R_2} -R_1 
	\approx s  + \eta \sqrt{ R_2 } - s + \beta \sqrt{s}\\
	&\approx  \eta \sqrt{ \frac{ s\mu p }{\delta}} + \beta \sqrt{s};\\
	&\frac{n- \left( R_1 + R_2 \right)}{\sqrt{R_1 + R_2}} 
	\approx \frac{\eta \sqrt{ \frac{ \mu p }{\delta }} + \beta }{\sqrt{1+  \frac{ \mu p}{\delta}}}.
\end{split}
\end{equation*}

Following Equations \ref{eq:x_small_s} and \ref{eq:y_equal_n} we get  
\begin{equation}
\label{eq:tilda_delta_1}
\begin{split}
	{\delta_1}&= \P(Y_{\lambda}=n) \P(X^1_{\lambda}\leq s) \approx \frac{1}{\sqrt{s}\sqrt{1+  \frac{ \mu p}{\delta} }} \phi \left( \frac{\eta \sqrt{ \frac{ \mu p}{\delta}} + \beta }{\sqrt{1+  \frac{ \mu p}{\delta}}}\right) \Phi \left( \frac{ \beta \sqrt{\frac{ \mu p}{\delta}} -  \eta } {\sqrt{  1+  \frac{ \mu p}{\delta} }}\right). 
\end{split}
\end{equation}
Now lets find an approximation for $\tilde{\delta_2}$.
\begin{equation*}
\begin{split}
		{\delta_2}&=  \sum_{i=s+1}^n  \frac{1}{s!s^{i-s}} R_1^i \frac{1}{(n-i)!} R_2^{n-i} e^{-\left(R_1 + R_2 \right)} = \frac{ e^{-\left(R_1 + R_2\right)}}{s!s^{-s}} \sum_{i=s+1}^n \rho^i \frac{1}{(n-i)!} R_2^{n-i} \\
		&= \frac{ e^{-\left(R_1 + R_2 \right)}}{s!s^{-s}} R_2^{n} \sum_{j=0}^{n-s-1} \frac{1}{j!} \left( \frac{\rho}{R_2} \right)^{n-j} = \frac{ e^{-\left(R_1 + R_2 \right)}}{s!s^{-s}}  \rho^{n}\sum_{j=0}^{n-s-1} \frac{1}{j!} \left( \frac{R_2} {\rho}\right)^{j}.\\
\end{split}
\end{equation*}
When comparing ${\delta_2}$ to $B_2$ form Equation \ref{eq:B1B2}, we observe that
\begin{equation*}
\begin{split}
		{\delta_2}&=  (1-\rho) B_2 \approx \frac{\beta}{\sqrt{s}} B_2.\\
\end{split}
\end{equation*}
Therefore, based on the approximation of $B_2$ from Lemma \ref{lem2} we get
\begin{equation*}
\begin{split}
		&lim_{\lambda \rightarrow \infty} {\delta_2} = \frac{\phi(\sqrt{\eta^2 + \beta^2})}{\sqrt{s}} e^{\frac{1}{2} \omega ^2} \Phi(\omega) .\\
\end{split}
\end{equation*}
This proves Theorem \ref{thm5}.
\endproof

The next theorem gives the approximation for the case where $\beta=0$.
\begin{theorem} \label{thm6}
Let the variables $\lambda$, $s$ and $n$ tend to $\infty$ simultaneously and satisfy the QED scaling conditions in \eqref{eq:twofoldscaling} with $\beta = 0$. Define and $\xi=\frac{R_1}{R_2}=\frac{\delta }{p\mu }$, then
\begin{equation}
\lim_{\lambda \rightarrow \infty} \sqrt{s} \P(block)= \frac{\sqrt{r}\phi(\gamma)\Phi(-\frac{\gamma \sqrt{r}}{\sqrt{1- r}})+\frac{1}{\sqrt{2 \pi}}\Phi(\eta)}{\int_{-\infty}^{0} \Phi\left( \eta - t \sqrt{\xi} \right) d\Phi(t) + \sqrt{\frac{1-r}{r}} \frac{1}{\sqrt{2\pi}} \left(\eta \Phi(\eta) + \phi(\eta) \right)},
\end{equation}
where $\eta = \frac{\gamma-\beta \sqrt{r}}{\sqrt{1-r}}$.
\end{theorem}

\proof{Proof:}
It follows from \eqref{eq:delay_probability} that the probability of blocking is given by
\begin{equation*}
\P_n = \frac{{\delta_1}+{\delta_2}}{\tilde{A}+\tilde{B}},
\end{equation*}
where
\begin{equation*}
\begin{split}
	{\delta_1}&=  \sum_{i=0}^s \frac{1}{i!} R_1^i \frac{1}{i!} R_2^{n-i}  e^{-\left(R_1 + R_2  \right)}\\
	{\delta_2}&=  \sum_{i=s+1}^n  \frac{1}{s!s^{i-s}} R_1^i \frac{1}{(n-i)!} R_2^{n-i} e^{-\left(R_1 + R_2 \right)}\\
	\tilde{A} &=  \sum_{\substack{i,j| i \leq s, \\i+j \leq n }} \frac{1}{i!j!} R_1^i  R_2^j  e^{-\left(R_1 + R_2 \right)}\\
	\tilde{B}&=  \frac{1}{s!} R_1^{s} \frac{1}{1-\rho} \sum_{l=0}^{n-s} \frac{1}{l!} R_2^{l} e^{-\left(R_1 + R_2 \right)} - \frac{1}{s!} R_1^{s} \frac{\rho^{n-s+1}}{1-\rho} \sum_{l=0}^{n-s} \frac{1}{l!} \left(\frac{R_2}{\rho} \right)^{l} e^{-\left(R_1 + R_2 \right)}.
\end{split}
\end{equation*}

Note that by Lemmas \ref{lem3} and \ref{lem4}
\begin{equation*}
\begin{split}
	&\lim_{\lambda \rightarrow \infty}\tilde{A} = \lim_{\lambda \rightarrow \infty} A = \int_{-\infty}^{\beta} \Phi\left( \eta + (\beta-t) \sqrt{\frac{\delta }{\mu p}} \right) d\Phi(t)\\
	&\lim_{\lambda \rightarrow \infty}\tilde{B} = \lim_{\lambda \rightarrow \infty} B = \sqrt{\frac{\mu p}{\delta }} \frac{1}{\sqrt{2\pi}} \left(\eta \Phi(\eta) + \phi(\eta) \right).	
\end{split}
\end{equation*}

In addition, the approximations for $\delta_1$ and $\delta_2$ are the same as the proof of Theorem \ref{thm5}.
\begin{equation*}
\begin{split}
	&\lim_{\beta \rightarrow 0}\lim_{\lambda \rightarrow \infty} \sqrt{s} \delta_2 = \frac{1}{\sqrt{2 \pi}}\Phi(\eta)\\
	&\lim_{\beta \rightarrow 0}\lim_{\lambda \rightarrow \infty}  \sqrt{s} \delta_1 =  \frac{1}{\sqrt{1+  \frac{ \mu p}{\delta}}} \phi \left( \frac{\eta }{\sqrt{1+  \frac{\delta }{\mu p }}}\right) \Phi \left( - \frac{ \eta } {\sqrt{  1+  \frac{ \mu p }{\delta } }}\right).	
\end{split}
\end{equation*}

This proves Theorem \ref{thm6}.

\endproof

\section{Proof of Proposition \ref{prop:stability_convergence}}\label{app:proof_stability_convergence}
 
Define 
\[
A(s,n) = \sum_{k=0}^s \frac{k}{s} \, \binom{n}{k} b^k ,\qquad 
B(s,n) = \sum_{k=s+1}^n \frac{k!}{s!} \, \binom{n}{k} s^{s-k} b^k, \qquad
C(s,n) = \sum_{k=0}^s \binom{n}{k} \, b^k,
\]
\[
\]
where $b = \delta/p\mu = r/(1-r)$. Then
\[
\rho_J(s,n) = \frac{A(s,n)+B(s,n)}{C(s,n)+B(s,n)}.
\]
Proving that $\rho_{\rm max}(s,n) \to 1$ as $R\to\infty$ with $s$ and $n$ as in \eqref{eq:twofoldscaling} is equivalent to showing that
\begin{equation}\label{eq:proof_stab_1}
1-\rho_{\rm max}(s,n) = \frac{C(s,n)-A(s,n)}{C(s,n)+B(s,n)} = \frac{(1+b)^{-n}[C(s,n)-A(s,n)]}{(1+b)^{-n}[C(s,n)+B(s,n)]} \to 0.
\end{equation}
First, we rewrite
\begin{align*}
(1+b)^{-n} A(s,n) 
&= (1+b)^{-n} \sum_{k=1}^s \frac{n}{s} \binom{n-1}{k-1} b^k \\
&= \frac{n}{s}\left(\frac{b}{1+b}\right)\sum_{k=0}^{s-1} \binom{n-1}{k} \left(\frac{b}{1+b}\right)^k \left(\frac{1}{1+b}\right)^{n-1-k}\\
&= \frac{r n}{s}\sum_{k=0}^{s-1} \binom{n-1}{k} r^k (1-r)^{n-1-k}
= \frac{r n}{s} \P( {\rm Bin}(n-1,r) \leq s-1 ) \\
&= \frac{rn}{s} \P\left( \frac{{\rm Bin}(n-1,r) - (n-1)r}{\sqrt{nr(1-r)}} \leq \frac{s-1 - (n-1)r}{\sqrt{nr(1-r)}} \right)
\to \Phi\left(\frac{\beta-\gamma\sqrt{r}}{\sqrt{1-r}}\right),
\end{align*}
since $nr/s = 1 + O(1/\sqrt{R_1})$.
Also, 
\begin{align*}
(1+b)^{-n} C(s,n) 
&= \sum_{k=0}^s \binom{n}{k} \left(\frac{b}{1+b}\right)^k \left(\frac{1}{1+b}\right)^{n-k} = \sum_{k=0}^s \binom{n}{k} r^k (1-r)^{n-k}\\
&= \P( {\rm Bin}(n,r) \leq s) \to \Phi\left(\frac{\beta-\gamma\sqrt{r}}{\sqrt{1-r}}\right).
\end{align*}
Therefore, we have $(1+b)^{-n}[C(s,n)-A(s,n)] \to 0$ as $\lambda\to\infty$.
For the remaining term,
\begin{align*}
(1+b)^{-n} B(s,n) 
&= (1+b)^{-n}\sum_{k=s+1}^n \binom{n}{k}\,\frac{k!}{s!} s^{s-k} b^k  = (1+b)^{-n}\frac{n!}{s!}\, s^s\sum_{k=s+1}^n \frac{1}{(n-k)!} \left(\frac{s}{b}\right)^{-k}\\
&= (1+b)^{-n} \frac{n!}{s!}\, s^s\, \left(\frac{b}{s}\right)^n \sum_{k=s+1}^n \frac{1}{(n-k)!} \left(\frac{s}{b}\right)^{n-k}
=  r^n\, \frac{n!}{s!} s^{s-n} \sum_{m=0}^{n-s-1} \frac{1}{m!} \left(\frac{s}{b}\right)^m\\
&= \left(\frac{r}{s}\right)^n \frac{n!}{s!} s^s \,{\rm e}^{s/b} \, \P({\rm Pois}(s/b)\leq n-s-1),
\end{align*}
in which 
\[
\P({\rm Pois}(s/b)\leq n-s-1) 
= \P\left(\frac{{\rm Pois}(s/b)-s/b}{\sqrt{s/b}} \leq \frac{n-s-1-s/b}{\sqrt{s/b}}\right) \to \Phi\left(\frac{\gamma-\beta/\sqrt{r}}{\sqrt{1-r}}\right),
\]
as $\lambda\to\infty$.
By Stirling's approximation,
\begin{align*}
\left(\frac{r}{s}\right)^n \frac{n!}{s!} s^s \,{\rm e}^{s/b} 
&\sim \left(\frac{r}{s}\right)^n \sqrt{\frac{n}{s}} \,\frac{n^n {\rm e}^{-n}}{s^s {\rm e}^{-s}}\, s^s \,{\rm e}^{s/b} \\
&= \left(\frac{rn}{s}\right)^n \sqrt{\frac{n}{s}} {\rm e}^{-n+s+s/b} = \left(\frac{rn}{s}\right)^n \sqrt{\frac{n}{s}} {\rm e}^{-n+s/r}.
\end{align*}
Since, 
\[
\frac{rn}{s} = 1 + \frac{\gamma\sqrt{r}-\beta}{\sqrt{R_1}} + O(1/R_1),
\]
we find $\sqrt{n/s} = 1/\sqrt{r} + O(1/\sqrt{R_1})$ and 
\begin{align*}
\log\left[ \left(\frac{rn}{s}\right)^n \sqrt{\frac{n}{s}} {\rm e}^{-n+s/r} \right]
&= n \log\left[ \frac{rn}{s}\right] - n+\frac{s}{r}\\
&= -n \left[ \left(1-\frac{rn}{s}\right) + \frac{1}{2}\left(1-\frac{rn}{s}\right)^2 + O(1/R^{3/2}) \right] + \frac{s}{r}\left(1-\frac{rn}{s}\right)\\
&= \frac{s}{r}\left(1-\frac{rn}{s}\right)^2 - \frac{n}{2}\left(1-\frac{rn}{s}\right)^2 + O(1/\sqrt{R_1})\\
&= \frac{(\gamma\sqrt{r} - \beta)^2}{2r} + O(1/\sqrt{R_1}),
\end{align*}
as $R\to\infty$ and hence, 
\[
(1+b)^{-n} B(s,n) \to \phi\left(\frac{\gamma\sqrt{r}-\beta}{\sqrt{r}}\right)\Phi\left(\frac{\gamma-\beta/\sqrt{r}}{\sqrt{1-r}}\right).
\]
Hence, we conclude that the denominator of \eqref{eq:proof_stab_1} converges to a constant value as $\lambda$ grows, and hence the $1-\rho_{\rm max}(s,n)\to 0$ as $\lambda\to\infty$.

\section{Numerical results on accuracy of the asymptotic approximations}

In this section, we examine to accuracy of the asymptotic approximations of Theorem \ref{thm:limits_YT} and the heuristic method in Section \ref{sec:QED_limit_holding}.
We perform our numerical experiments for three case instances, with parameter settings as in Table \ref{tab:parameter_settings}.

\begin{table}[htb]
\centering
\begin{tabular}{|r|rrrr|}
\hline
		&	$\mu$	&	$\d$	&	$p$	&	$r$ \\
\hline
Case 1  &   1   & 0.10 & 0.90 & 0.10 \\
Case 2 	& 	1	& 0.25 &  0.75 &  0.25\\
Case 3  &   1	& 0.50  &  0.50  &  0.50 \\
\hline
\end{tabular}
\caption{Parameter settings for numerical experiments.}
\label{tab:parameter_settings}
\end{table}

\subsection{Restricted Erlang-R model with blocking}
\label{app:accuracy_blocking}
\begin{table}[h] \centering
\begin{tabular}{|r|rrr|rrr|}
\cline{2-7}\multicolumn{1}{r|}{} & \multicolumn{3}{c|}{$\beta = 1,\ \g = 1$} & \multicolumn{3}{c|}{$\beta = 1,\ \g = 2$} \bigstrut\\
\hline
$R_1$   & $\P({\rm d})$ & $\sqrt{R_1}\P({\rm b})$ & $\sqrt{R_1}\E[W]$ & $\P({\rm d})$ & $\sqrt{R_1}\P({\rm b})$ & $\sqrt{R_1}\E[W]$ \bigstrut\\
\hline
5     & 0.1270 & 0.0900 & 0.2283 & 0.1553 & 0.0212 & 0.1085 \bigstrut[t]\\
10    & 0.1340 & 0.0910 & 0.1919 & 0.1628 & 0.0206 & 0.1205 \\
25    & 0.1981 & 0.0945 & 0.1614 & 0.2356 & 0.0216 & 0.2145 \\
50    & 0.1513 & 0.0963 & 0.1588 & 0.1830 & 0.0205 & 0.1496 \\
100   & 0.1880 & 0.0956 & 0.1532 & 0.2231 & 0.0224 & 0.2055 \\
250   & 0.1797 & 0.0971 & 0.1399 & 0.2143 & 0.0219 & 0.2057 \\
\hline
\multicolumn{1}{r|}{} & \textit{0.1767} & \textit{0.0981} & \textit{0.1437} & \textit{0.2108} & \textit{0.0217} & \textit{0.1947} \bigstrut[b]\\
\cline{2-7}\end{tabular}%
\vspace{3mm}
\begin{tabular}{|r|rrr|rrr|}
\cline{2-7}\multicolumn{1}{r|}{} & \multicolumn{3}{c|}{$\beta = 2,\ \g = 1$} & \multicolumn{3}{c|}{$\beta = 2,\ \g = 2$} \bigstrut\\
\hline
$R_1$   & $\P({\rm d})$ & $\sqrt{R_1}\P({\rm b})$ & $\sqrt{R_1}\E[W]$ & $\P({\rm d})$ & $\sqrt{R_1}\P({\rm b})$ & $\sqrt{R_1}\E[W]$ \bigstrut\\
\hline
5     & 0.0237 & 0.0868 & 0.0282 & 0.0322 & 0.0192 & 0.0391 \bigstrut[t]\\
10    & 0.0206 & 0.0872 & 0.0188 & 0.0278 & 0.0183 & 0.0264 \\
25    & 0.0277 & 0.0876 & 0.0123 & 0.0363 & 0.0174 & 0.0174 \\
50    & 0.0185 & 0.0913 & 0.0116 & 0.0249 & 0.0175 & 0.0166 \\
100   & 0.0232 & 0.0888 & 0.0103 & 0.0303 & 0.0183 & 0.0145 \\
250   & 0.0203 & 0.0905 & 0.0079 & 0.0267 & 0.0179 & 0.0109 \bigstrut[b]\\
\hline
\multicolumn{1}{r|}{} & \textit{0.0188} & \textit{0.0914} & \textit{0.0084} & \textit{0.0247} & \textit{0.0177} & \textit{0.0118} \bigstrut\\
\cline{2-7}\end{tabular}%

\caption{Numerical results for Erlang-R model with blocking for Case 1.}
\label{tab:numerics_case1}
\end{table}
\begin{table}[h] \centering
\begin{tabular}{|r|rrr|rrr|}
\cline{2-7}\multicolumn{1}{r|}{} & \multicolumn{3}{c|}{$\beta = 1,\ \g = 1$} & \multicolumn{3}{c|}{$\beta = 1,\ \g = 2$} \bigstrut\\
\hline
$R_1$   & $\P({\rm d})$ & $\sqrt{R_1}\P({\rm b})$ & $\sqrt{R_1}\E[W]$ & $\P({\rm d})$ & $\sqrt{R_1}\P({\rm b})$ & $\sqrt{R_1}\E[W]$ \bigstrut\\
\hline
5     & 0.0911 & 0.1538 & 0.0479 & 0.1431 & 0.0345 & 0.0909 \bigstrut[t]\\
10    & 0.1010 & 0.1498 & 0.0560 & 0.1520 & 0.0326 & 0.1025 \\
25    & 0.1594 & 0.1509 & 0.1058 & 0.2192 & 0.0405 & 0.1785 \\
50    & 0.1201 & 0.1506 & 0.0726 & 0.1697 & 0.0381 & 0.1248 \\
100   & 0.1514 & 0.1539 & 0.1001 & 0.2088 & 0.0398 & 0.1704 \\
250   & 0.1459 & 0.1524 & 0.0957 & 0.2003 & 0.0397 & 0.1618 \\
\hline
\multicolumn{1}{r|}{} & \textit{0.1429} & \textit{0.1569} & \textit{0.0940} & \textit{0.1976} & \textit{0.0391} & \textit{0.1617} \bigstrut[b]\\
\cline{2-7}\end{tabular}%

\vspace{3mm}

\begin{tabular}{|r|rrr|rrr|}
\cline{2-7}\multicolumn{1}{r|}{} & \multicolumn{3}{c|}{$\beta = 2,\ \g = 1$} & \multicolumn{3}{c|}{$\beta = 2,\ \g = 2$} \bigstrut\\
\hline
$R_1$   & $\P({\rm d})$ & $\sqrt{R_1}\P({\rm b})$ & $\sqrt{R_1}\E[W]$ & $\P({\rm d})$ & $\sqrt{R_1}\P({\rm b})$ & $\sqrt{R_1}\E[W]$ \bigstrut\\
\hline
5     & 0.0130 & 0.1484 & 0.0044 & 0.0277 & 0.0294 & 0.0109 \bigstrut[t]\\
10    & 0.0121 & 0.1432 & 0.0042 & 0.0244 & 0.0267 & 0.0098 \\
25    & 0.0182 & 0.1383 & 0.0070 & 0.0319 & 0.0295 & 0.0141 \\
50    & 0.0119 & 0.1415 & 0.0043 & 0.0216 & 0.0301 & 0.0090 \\
100   & 0.0154 & 0.1413 & 0.0059 & 0.0270 & 0.0290 & 0.0119 \\
250   & 0.0136 & 0.1403 & 0.0051 & 0.0236 & 0.0291 & 0.0103 \bigstrut[b]\\
\hline
\multicolumn{1}{r|}{\textit{}} & \textit{0.0126} & \textit{0.1445} & \textit{0.0048} & \textit{0.0220} & \textit{0.0284} & \textit{0.0097} \bigstrut\\
\cline{2-7}\end{tabular}%
\caption{Numerical results for Erlang-R model with blocking for Case 2.}
\label{tab:numerics_case2}
\end{table}

\begin{table}[h] \centering
\begin{tabular}{|r|rrr|rrr|}
\cline{2-7}\multicolumn{1}{r|}{} & \multicolumn{3}{c|}{$\beta = 1,\ \g = 1$} & \multicolumn{3}{c|}{$\beta = 1,\ \g = 2$} \bigstrut\\
\hline
$R_1$   & $\P({\rm d})$ & $\sqrt{R_1}\P({\rm b})$ & $\sqrt{R_1}\E[W]$ & $\P({\rm d})$ & $\sqrt{R_1}\P({\rm b})$ & $\sqrt{R_1}\E[W]$ \bigstrut\\
\hline
5     & 0.0547 & 0.1945 & 0.0221 & 0.1181 & 0.0604 & 0.0617 \bigstrut[t]\\
10    & 0.0579 & 0.2158 & 0.0237 & 0.1325 & 0.0526 & 0.0746 \\
25    & 0.1113 & 0.2086 & 0.0544 & 0.1959 & 0.0641 & 0.1311 \\
50    & 0.0813 & 0.2050 & 0.0363 & 0.1523 & 0.0562 & 0.0933 \\
100   & 0.1060 & 0.2146 & 0.0509 & 0.1873 & 0.0632 & 0.1250 \\
250   & 0.1006 & 0.2179 & 0.0475 & 0.1820 & 0.0596 & 0.1214 \\
\hline
\multicolumn{1}{r|}{} & 0.1011 & \textit{0.2185} & \textit{0.0478} & \textit{0.1792} & \textit{0.0605} & \textit{0.1199} \bigstrut[b]\\
\cline{2-7}\end{tabular}%

\vspace{3mm}

\begin{tabular}{|r|rrr|rrr|}
\cline{2-7}\multicolumn{1}{r|}{} & \multicolumn{3}{c|}{$\beta = 2,\ \g = 1$} & \multicolumn{3}{c|}{$\beta = 2,\ \g = 2$} \bigstrut\\
\hline
$R_1$   & $\P({\rm d})$ & $\sqrt{R_1}\P({\rm b})$ & $\sqrt{R_1}\E[W]$ & $\P({\rm d})$ & $\sqrt{R_1}\P({\rm b})$ & $\sqrt{R_1}\E[W]$ \bigstrut\\
\hline
5     & 0.0034 & 0.1888 & 0.0009 & 0.0175 & 0.0510 & 0.0057 \bigstrut[t]\\
10    & 0.0030 & 0.2093 & 0.0008 & 0.0172 & 0.0416 & 0.0058 \\
25    & 0.0070 & 0.1937 & 0.0020 & 0.0243 & 0.0440 & 0.0089 \\
50    & 0.0043 & 0.1946 & 0.0011 & 0.0163 & 0.0414 & 0.0056 \\
100   & 0.0061 & 0.1999 & 0.0017 & 0.0207 & 0.0431 & 0.0076 \\
250   & 0.0052 & 0.2037 & 0.0014 & 0.0185 & 0.0401 & 0.0067 \bigstrut[b]\\
\hline
\multicolumn{1}{r|}{} & 0.0052 & \textit{0.2039} & \textit{0.0014} & \textit{0.0173} & \textit{0.0404} & \textit{0.0063} \bigstrut\\
\cline{2-7}\end{tabular}%

\caption{Numerical results for Erlang-R model with blocking for Case 3.}
\label{tab:numerics_case3}
\end{table}

\newpage
\subsection{Restricted Erlang-R model with holding}
\label{app:accuracy_holding}

\begin{table}[h] \centering
\begin{tabular}{|r|rr|rr|}
\cline{2-5}\multicolumn{1}{r|}{} & \multicolumn{2}{c|}{$\b=1,\ \g=1$} & \multicolumn{2}{c|}{$\b=1,\ \g=2$} \bigstrut\\
\hline
$R_1$     & $\P({\rm d})$ & $\sqrt{R_1}\E[W]$ & $\P({\rm d})$ & $\sqrt{R_1}\E[W]$ \bigstrut\\
\hline
5     & 0.1532 & 0.1031 & 0.1628 & 0.1216 \bigstrut[t]\\
10    & 0.1622 & 0.1272 & 0.1697 & 0.1331 \\
25    & 0.2340 & 0.2116 & 0.2413 & 0.2342 \\
50    & 0.1817 & 0.1468 & 0.1890 & 0.1678 \\
100   & 0.2199 & 0.1931 & 0.2304 & 0.2269 \\
250   & 0.2110 & 0.1852 & 0.2176 & 0.2230 \bigstrut[b]\\
\hline
\multicolumn{1}{r|}{} & \textit{0.2076} & \textit{0.1777} & \textit{0.2187} & \textit{0.2050} \bigstrut\\
\cline{2-5}\end{tabular}%

\vspace{3mm}

\begin{tabular}{|r|rr|rr|}
\cline{2-5}\multicolumn{1}{c|}{} & \multicolumn{2}{c|}{$\b=2,\ \g=1$} & \multicolumn{2}{c|}{$\b=2,\ \g=1$} \bigstrut\\
\hline
$R_1$     & $\P({\rm d})$ & $\sqrt{R_1}\E[W]$ & $\P({\rm d})$ & $\sqrt{R_1}\E[W]$ \bigstrut\\
\hline
5     & 0.0310 & 0.0121 & 0.0344 & 0.0148 \bigstrut[t]\\
10    & 0.0267 & 0.0123 & 0.0292 & 0.0128 \\
25    & 0.0348 & 0.0171 & 0.0373 & 0.0184 \\
50    & 0.0240 & 0.0108 & 0.0258 & 0.0125 \\
100   & 0.0293 & 0.0143 & 0.0317 & 0.0163 \\
250   & 0.0256 & 0.0120 & 0.0276 & 0.0145 \\
\hline
\multicolumn{1}{r|}{\textit{}} & \textit{0.0229} & \textit{0.0104} & \textit{0.0257} & \textit{0.0124} \bigstrut[b]\\
\cline{2-5}\end{tabular}%

\caption{Simulated and approximated probability of delay in Erlang-R model with holding for Case 1.}
\label{tab:heuristic_case1}
\end{table}
\begin{table}[h]\centering

\begin{tabular}{|r|rr|rr|}
\cline{2-5}\multicolumn{1}{r|}{} & \multicolumn{2}{c|}{$\b=1,\ \g=1$} & \multicolumn{2}{c|}{$\b=1,\ \g=2$} \bigstrut\\
\hline
$R_1$     & $\P({\rm d})$ & $\sqrt{R_1}\E[W]$ & $\P({\rm d})$ & $\sqrt{R_1}\E[W]$ \bigstrut\\
\hline
5     & 0.1327 & 0.0740 & 0.1620 & 0.1096 \bigstrut[t]\\
10    & 0.1446 & 0.0894 & 0.1683 & 0.1207 \\
25    & 0.2204 & 0.1631 & 0.2442 & 0.2203 \\
50    & 0.1694 & 0.1122 & 0.1888 & 0.1507 \\
100   & 0.2098 & 0.1524 & 0.2322 & 0.2111 \\
250   & 0.2033 & 0.1534 & 0.2190 & 0.1979 \bigstrut[b]\\
\hline
\multicolumn{1}{r|}{} & \textit{0.1840} & \textit{0.1277} & \textit{0.2109} & \textit{0.1759} \bigstrut\\
\cline{2-5}\end{tabular}%
\vspace{3mm}
\begin{tabular}{|r|rr|rr|}
\cline{2-5}\multicolumn{1}{c|}{} & \multicolumn{2}{c|}{$\b=2,\ \g=1$} & \multicolumn{2}{c|}{$\b=2,\ \g=1$} \bigstrut\\
\hline
$R_1$     & $\P({\rm d})$ & $\sqrt{R_1}\E[W]$ & $\P({\rm d})$ & $\sqrt{R_1}\E[W]$ \bigstrut\\
\hline
5     & 0.0219 & 0.0079 & 0.0322 & 0.0137 \bigstrut[t]\\
10    & 0.0199 & 0.0073 & 0.0284 & 0.0115 \\
25    & 0.0283 & 0.0128 & 0.0375 & 0.0163 \\
50    & 0.0190 & 0.0078 & 0.0255 & 0.0107 \\
100   & 0.0244 & 0.0097 & 0.0314 & 0.0151 \\
250   & 0.0214 & 0.0083 & 0.0272 & 0.0134 \bigstrut[b]\\
\hline
\multicolumn{1}{r|}{\textit{}} & \textit{0.0169} & \textit{0.0066} & \textit{0.0234} & \textit{0.0104} \bigstrut\\
\cline{2-5}\end{tabular}%

\caption{Simulated and approximated probability of delay in Erlang-R model with holding for Case 2.}
\label{tab:heuristic_case2}
\end{table}

\begin{table}
\centering
\begin{tabular}{|r|rr|rr|}
\cline{2-5}\multicolumn{1}{r|}{} & \multicolumn{2}{c|}{$\b=1,\ \g=1$} & \multicolumn{2}{c|}{$\b=1,\ \g=2$} \bigstrut\\
\hline
$R_1$     & $\P({\rm d})$ & $\sqrt{R_1}\E[W]$ & $\P({\rm d})$ & $\sqrt{R_1}\E[W]$ \bigstrut\\
\hline
5     & 0.0977 & 0.0413 & 0.1521 & 0.0851 \bigstrut[t]\\
10    & 0.1070 & 0.0469 & 0.1648 & 0.1028 \\
25    & 0.1926 & 0.1076 & 0.2421 & 0.1874 \\
50    & 0.1431 & 0.0727 & 0.1876 & 0.1342 \\
100   & 0.1855 & 0.1012 & 0.2282 & 0.1714 \\
250   & 0.1775 & 0.0963 & 0.2217 & 0.1765 \bigstrut[b]\\
\hline
\multicolumn{1}{r|}{} & \textit{0.1442} & \textit{0.0711} & \textit{0.1981} & \textit{0.1354} \bigstrut\\
\cline{2-5}\end{tabular}%

\vspace{3mm}

\begin{tabular}{|r|rr|rr|}
\cline{2-5}\multicolumn{1}{c|}{} & \multicolumn{2}{c|}{$\b=2,\ \g=1$} & \multicolumn{2}{c|}{$\b=2,\ \g=1$} \bigstrut\\
\hline
$R_1$     & $\P({\rm d})$ & $\sqrt{R_1}\E[W]$ & $\P({\rm d})$ & $\sqrt{R_1}\E[W]$ \bigstrut\\
\hline
5     & 0.0072 & 0.0019 & 0.0250 & 0.0081 \bigstrut[t]\\
10    & 0.0067 & 0.0018 & 0.0235 & 0.0082 \\
25    & 0.0148 & 0.0043 & 0.0325 & 0.0133 \\
50    & 0.0092 & 0.0025 & 0.0217 & 0.0081 \\
100   & 0.0132 & 0.0038 & 0.0277 & 0.0105 \\
250   & 0.0114 & 0.0033 & 0.0246 & 0.0099 \bigstrut[b]\\
\hline
\multicolumn{1}{r|}{\textit{}} & \textit{0.0078} & \textit{0.0022} & \textit{0.0188} & \textit{0.0069} \bigstrut\\
\cline{2-5}\end{tabular}

\caption{Simulated and approximated probability of delay in Erlang-R model with holding for Case 3.}
\label{tab:heuristic_case3}
\end{table}

\end{document}